\documentclass[letterpaper,11pt]{article}
\usepackage[margin=1in]{geometry}  

\usepackage{listings}
\usepackage{amsmath}
\usepackage{amsthm}
\usepackage{tikz}
\usepackage{caption}
\usepackage{array}
\usepackage{mdwmath}
\usepackage{multirow}
\usepackage{mdwtab}
\usepackage{eqparbox}
\usepackage{multicol}
\usepackage{amsfonts}
\usepackage{tikz}
\usepackage{multirow,bigstrut,threeparttable}
\usepackage{amsthm}
\usepackage{array}
\usepackage{bbm}
\usepackage{subfigure}
\usepackage{epstopdf}
\usepackage{mdwmath}
\usepackage{mdwtab}
\usepackage{eqparbox}
\usepackage{tikz}
\usepackage{latexsym}
\usepackage{cite}
\usepackage{amssymb}
\usepackage{bm}
\usepackage{amssymb}
\usepackage{graphicx}
\usepackage{mathrsfs}
\usepackage{epsfig}
\usepackage{psfrag}
\usepackage{setspace}
\usepackage[
            CJKbookmarks=true,
            bookmarksnumbered=true,
            bookmarksopen=true,
            colorlinks=true,
            citecolor=red,
            linkcolor=blue,
            anchorcolor=red,
            urlcolor=blue
            ]{hyperref}
\usepackage{algorithm}
\usepackage{algpseudocode}
\usepackage{stfloats}

\theoremstyle{plain}
\newtheorem{theorem}{Theorem}

\newtheorem{lemma}{Lemma}
\newtheorem{remark}{Remark}
\newtheorem{corollary}{Corollary}

\newtheorem*{question*}{Question}
\newtheorem{assumption}{Assumption}

\def \bP {\mathbb{P}}
\def \bE {\mathbb{E}}
\def \bR {\mathbb{R}}

\def \var {\mathsf{Var}}

\usepackage{xspace}

\newcommand{\torus}{\mathbb{T}}
\newcommand{\supp}{\mathsf{supp}}
\newcommand{\Lip}{\mathrm{Lip}}

\newcommand{\Lipstar}{\mathrm{Lip}^{\star}}

\newcommand{\stepa}[1]{\overset{\rm (a)}{#1}}
\newcommand{\stepb}[1]{\overset{\rm (b)}{#1}}

\newcommand{\floor}[1]{{\left\lfloor {#1} \right \rfloor}}

\newcommand{\reals}{\mathbb{R}}
\newcommand{\naturals}{\mathbb{N}}

\newcommand{\Expect}{\mathbb{E}}

\newcommand{\ie}{i.e.\xspace}
\newcommand{\iid}{i.i.d.\xspace}

\newcommand{\pth}[1]{\left( #1 \right)}
\newcommand{\qth}[1]{\left[ #1 \right]}
\newcommand{\sth}[1]{\left\{ #1 \right\}}

\newcommand{\iiddistr}{{\stackrel{\text{\iid}}{\sim}}}
\newcommand{\inddistr}{{\stackrel{\text{ind}}{\sim}}}

\newcommand{\Bern}{\text{Bern}}
\newcommand{\Poi}{\mathsf{Poi}}

\definecolor{myblue}{rgb}{.8, .8, 1}
\definecolor{mathblue}{rgb}{0.2472, 0.24, 0.6} 
\definecolor{mathred}{rgb}{0.6, 0.24, 0.442893}
\definecolor{mathyellow}{rgb}{0.6, 0.547014, 0.24}

\newcommand{\calP}{{\mathcal{P}}}

\newcommand{\diverge}{\to \infty}

\usepackage{cleveref}
\crefname{lemma}{Lemma}{Lemmas}
\Crefname{lemma}{Lemma}{Lemmas}
\crefname{thm}{Theorem}{Theorems}
\Crefname{thm}{Theorem}{Theorems}
\crefname{assumption}{Assumption}{Assumptions}
\Crefname{assumption}{Assumption}{Assumptions}
\crefformat{equation}{(#2#1#3)}

\newcommand{\prettyref}{\Cref}

\usepackage{enumerate}

\begin{document}

\title{Optimal rates of entropy estimation over Lipschitz balls}
\author{Yanjun Han, Jiantao Jiao, Tsachy Weissman, and Yihong Wu\thanks{Y.~Han and T.~Weissman are with the Department of Electrical Engineering, Stanford University, email: 
\url{{yjhan,tsachy}@stanford.edu}. 
J.~Jiao is with the Department of Electrical Engineering and Computer Sciences and Department of Statistics, University of California Berkeley, email: \url{jiantao@eecs.berkeley.edu}.
Y.~Wu is with the Department of Statistics and Data Science, Yale University, email: 
\url{yihong.wu@yale.edu}. Y.~Han and T.~Weissman are supported in part by the NSF grants CCF-0939370 and CCF-1527105. J.~Jiao is supported in part by the NSF grants CCF-0939270, CCF-1527105 and IIS-1901252. Y.~Wu is supported in part by the NSF grants CCF-1527105, CCF-1900507, an NSF CAREER award CCF-1651588, and an Alfred Sloan fellowship.
}}

\maketitle

%
%
%
%
%
%

\begin{abstract}
We consider the problem of minimax estimation of the entropy of a density over Lipschitz balls. Dropping the usual assumption that the density is bounded away from zero, we obtain the minimax rates $(n\ln n)^{-s/(s+d)} + n^{-1/2}$ for $0<s\leq 2$ for densities supported on $[0,1]^d$, where $s$ is the smoothness parameter and $n$ is the number of independent samples. We generalize the results to densities with unbounded support: given an Orlicz functions $\Psi$ of rapid growth (such as the sub-exponential and sub-Gaussian classes), the minimax rates for densities with bounded $\Psi$-Orlicz norm increase to $(n\ln n)^{-s/(s+d)} (\Psi^{-1}(n))^{d(1-d/p(s+d))} + n^{-1/2}$, where $p$ is the norm parameter in the Lipschitz ball. We also show that the integral-form plug-in estimators with kernel density estimates fail to achieve the minimax rates, and characterize their worst case performances over the Lipschitz ball. 

One of the key steps in analyzing the bias relies on a novel application of the Hardy-Littlewood maximal inequality, which also leads to a new inequality on the Fisher information that may be of independent interest.
\end{abstract}

\tableofcontents


\section{Introduction}

Estimation of functionals of data generating distributions is a fundamental problem in statistics. While this problem is relatively well-understood in finite dimensional parametric models \cite{bickel1993efficient,van2000asymptotic}, the corresponding nonparametric counterparts are often much more challenging and have attracted tremendous interest over the last two decades. Initial efforts have focused on inference of linear, quadratic, and cubic functionals in Gaussian white noise and density models and have laid the foundation for the ensuing research. We do not attempt to survey the extensive literature in this area, but instead refer to the interested reader to, e.g., \cite{hall1987estimation,bickel1988estimating,donoho1990minimax2,fan1991estimation,birge1995estimation,kerkyacharian1996estimating,laurent1996efficient,nemirovski2000topics,cai2003note,cai2005nonquadratic,tchetgen2008minimax} and the references therein. 

The monograph by~\cite{nemirovski2000topics} provides a general treatment of estimating smooth functionals and discusses cases where efficient parametric rate of estimation is possible. Recently, there has been progress toward the understanding of more complex nonparametric functionals over substantially more general observational models. These include causal effect functionals in observational studies and mean functionals in missing data models. For more details, we refer to \cite{robins2008higher,robins2017higher,mukherjee2017semiparametric}, which considers a general recipe to yield minimax estimation of a large class of nonparametric functionals common in statistical literature. 
However, among the class of nonparametric functionals considered in literature, most of the research endeavors, at least from the point of view of minimax optimality, have focused on ``smooth functionals" (see \cite{robins2008higher} for a discussion on general classes of ``smooth functionals"). 

In contrast, the results on optimal estimation of non-smooth functionals have been less comprehensive~\cite{hasminskii1980some,donoho1997renormalizing, korostelev2012minimax}. Notably, the seminal papers of \cite{lepski1999estimation} and \cite{Cai--Low2011} considered the estimating of $L_r$-norms in Gaussian mean models. Subsequently, significant progress has been made on testing and estimation of non-smooth functionals, such as the Shannon entropy, support size, total variation and Kullback-Leibler (KL) divergence, for discrete distributions on large domains~(see, e.g., \cite{Paninski2004,valiant2011power,Jiao--Venkat--Han--Weissman2015minimax,wu2016minimax,han2016minimax,jiao2018minimax,bu2018estimation,wu2019chebyshev}).

An important non-smooth functional of probability density function is the \emph{entropy}, which has been the subject of extensive studies. The main goal of this paper is to resolve the minimax rates of entropy estimation in the density model under smoothness constraints, specifically, over Lipschitz classes.
To this end, consider the following i.i.d.~sampling model:
\begin{align*}
X_1,\cdots,X_n \overset{ \text{i.i.d.} }{\sim} f
\end{align*}
where $f$ is a probability density function on $\reals^d$. The goal is to estimate the entropy (also known as the differential entropy in the information theory literature) of the density $f$:
\begin{align*}
H(f) \triangleq \int_{\reals^d} -f(x)\ln f(x)dx.
\end{align*}
This problem has extensive applications in various fields such as information theory, neuroscience, time series, and machine learning (c.f.~\cite{Kraskov2004,costa2004geodesic,hlavavckova2007causality} and the survey \cite{beirlant1997nonparametric,WKV09b}).

A prevalent assumption in nonparametric entropy estimation is that $f(x)\ge c$ everywhere for some constant $c>0$~\cite{hall1984limit,joe1989estimation,van1992estimating,sricharan2012estimation, kandasamy2015nonparametric}, while others impose various assumptions quantifying on average how close the density is to zero~\cite{hall1993estimation, levit1978asymptotically, tsybakov1996root,el2009entropy, gao2016breaking, singh2016finite,delattre2017kozachenko,gao2018demystifying}. Assuming the density is bounded away from zero makes entropy a \emph{smooth} functional, consequently, the general technique for estimating smooth nonparametric functionals~\cite{robins2008higher,robins2017higher,mukherjee2017semiparametric} can be directly applied to achieve the minimax rate $\Theta(n^{-4s/(4s+d)}+n^{-1/2})$. 

It is well-known that smoothness conditions or shape restrictions are often necessary for non-parametric problems. We allow the density to be arbitrarily close to zero and adopt Lipschitz ball $\Lip_{s,p,d}(L)$ smoothness assumptions. Assume smoothness parameter $s>0$,
 norm parameter $p\in [2,\infty)$ and dimensionality $d\in \mathbb{N} \triangleq \{1,2,\cdots\}$. The Lipschitz ball is defined as
\begin{align}
\Lip_{s,p,d}(L) \triangleq \{f: \|f\|_{\Lip_{s,p,d}} \le L\} \cap \{ f: \supp(f)\subseteq [0,1]^d \}, 
\label{eq:Lip}
\end{align}
where with $r\triangleq \lceil s\rceil$, the Lipschitz norm $\|\cdot\|_{\Lip_{s,p,d}}$ is defined as
\begin{align}
\|f\|_{\Lip_{s,p,d}} &\triangleq \|f\|_p + \sup_{t>0} t^{-s}\omega_r(f,t)_p, \label{eq:lipnorm}\\
\omega_r(f,t)_p &\triangleq \sup_{e\in \mathbb{R}^d, |e|\le 1} \|\Delta_{te}^rf(\cdot)\|_p, \label{eq.moduli_of_smoothness}\\
\Delta_{h}^r f(x) &\triangleq \sum_{k=0}^r (-1)^{r-k}\binom{r}{k}f\pth{x+\pth{k-\frac{r}{2}}h}, \qquad h\in \mathbb{R}^d. \label{eq:finitediff}
\end{align}
Here $|x|$ denotes the Euclidean norm of a vector $x\in\reals^d$ and  $\|\cdot\|_p$ denotes the $L_p$ norm of measurable functions on $\bR^d$. Note that the $L_p$ norm in \prettyref{eq.moduli_of_smoothness} is taken over the whole space $\mathbb{R}^d$ 
to ensure that the density $f$ vanishes smoothly at the boundary. 
For example, any density whose derivatives up to order $\lceil s\rceil - 1$ all vanish at the boundary of $[0,1]^d$ suffices.


We characterize the minimax rates of estimating $H(f)$ over the Lipschitz ball $\Lip_{s,p,d}(L)$ in the following theorem. 
\begin{theorem}[Compactly supported densities]\label{thm:main}
	For any $d\in \mathbb{N}, 0<s\le 2$ and $2\le p<\infty$, 
	there exist constants $L_0>1$ and $c,C>0$ depending on $s,p,d$, such that 
	for any $L_0 \le L \le (n\ln n)^{s/d}$ and any $n\in\naturals$,
	\begin{align}
	c( (n\ln n)^{-\frac{s}{s+d}}L^{\frac{d}{s+d}} + n^{-\frac{1}{2}}\ln L) 
	&\leq \left(\inf_{\hat{H}}\sup_{f\in \Lip_{s,p,d}(L)} \bE_f (\hat{H}-H(f))^2 \right)^{\frac{1}{2}} \label{eq:main}\\
	&\leq C((n\ln n)^{-\frac{s}{s+d}}L^{\frac{d}{s+d}} + n^{-\frac{1}{2}}\ln L). \nonumber 
	\end{align}
Moreover, the lower bound part of \prettyref{eq:main} holds for any $s>0, 1\leq p<\infty$.
\end{theorem}
\begin{remark}\label{remark:radius}
A careful inspection of the proof of Theorem \ref{thm:main} reveals that, for $s\in (0,2]$, $p\ge 2$ and $L_0\le L\le L'\le (n\ln n)^{s/d}$, the minimax $L_2$ risk for entropy estimation over densities supported on $[0,1]^d$ with $\|f\|_p \le L$ and $\sup_{t>0} t^{-s}\omega_{\lceil s\rceil}(f,t)_p\le L'$ is
\begin{align}
\Theta\left( (n\ln n)^{-\frac{s}{s+d}} (L')^{\frac{d}{s+d}} + n^{-\frac{1}{2}}\ln L \right).
\label{eq:rates-smoothnessL}
\end{align}
Hence, by scaling,\footnote{Indeed, let $\tilde f(x)\triangleq R^d f(Rx)$ denote the density of $X_i/R$. Then $H(\tilde f)=H(f)-d\log R$, $\|\tilde f\|_p=R^{d(1-1/p)} \|f\|_p$, $\Delta_h \tilde f(x) = R^d \Delta_{Rh}^r f(Rx)$ and hence $\sup_{t>0} t^{-s}\omega_{r}(\tilde f,t)_p = R^{d(1-1/p)+s} \sup_{t>0} t^{-s}\omega_{r}(f,t)_p$. } if the density is supported on $[0,R]^d$ with $R\ge 1$ and satisfies $\|f\|_{\text{\rm Lip}_{s,p,d}}\le L$ with $R^{d(1-1/p)}L\ge L_0$ and $R^{s+d(1-1/p)}L \le (n\ln n)^{s/d}$, the minimax $L_2$ risk is
\begin{align}
\Theta\left( (n\ln n)^{-\frac{s}{s+d}} \left(R^{s+d(1-1/p)}L \right)^{\frac{d}{s+d}} + n^{-\frac{1}{2}}\ln \left(R^{d(1-1/p)}L \right) \right).
\label{eq:rates-Rsupport}
\end{align}

\end{remark}
\begin{remark}\label{rmk:param}
A direct consequence of Theorem \ref{thm:main} is that, for fixed parameters $s>0, p\in [2,\infty)$ and $L>L_0$, when $d=1,2$, the parametric rate $\Theta(n^{-1/2})$ is attainable for entropy estimation over the Lipschitz ball $\Lip_{s,p,d}(L)$  if and only if $s\ge d$. Moreover, when $d\ge 3$, the parametric rate cannot be attained for all $s<d$. 
\end{remark}


To the best of our knowledge, \Cref{thm:main} is the first characterization of the minimax rate for nonparametric entropy estimation in arbitrary dimensions over Lipschitz balls (or even the simpler H\"{o}lder balls) without assuming the density is bounded away from zero. One observes that the exponents of $n$ or $L$ in the minimax rates \prettyref{eq:main} and \prettyref{eq:rates-smoothnessL} do not depend on the norm parameter $p$ under the assumption that $2\leq p<\infty$. Another observation from Remark \ref{rmk:param} is that the level of smoothness required for the parametric rate is $s\ge d$, which is more than $s\ge d/4$ that suffices for densities bounded away from zero on the support $[0,1]^d$~\cite{laurent1996efficient}, and also more than $s\ge d/2$ that suffices for densities satisfying a relative version of H\"{o}lder smoothness~\cite{berrett2019efficient}.

We construct the minimax rate-optimal estimator by first approximating the density $f$ by $f_h$ (a locally smoothed version of $f$), and then designing estimators to estimate $H(f_h)$. The key advantage of estimating $H(f_h)$ over estimating $H(f)$ is that for each $x\in [0,1]^d$ and positive integer $k\le \ln n$, the $k$-th power of $f_h(x)$ admits an unbiased estimator using a $U$-statistic, which enables us to employ the techniques of best polynomial approximation and Taylor expansion to reduce the bias in estimating $H(f_h)$. Moreover, our estimator is directly constructed and proved for the density model rather than the Poissonized model, unlike most prior work based on polynomial approximation \cite{Jiao--Venkat--Han--Weissman2015minimax,wu2016minimax}.

We improve the best known minimax lower bound for estimating non-smooth nonparametric functionals. The well-known lower bound $\Theta(n^{-4s/(4s+d)}+n^{-1/2})$~\cite{birge1995estimation}, which is optimal for smooth functionals such as quadratic functionals, is loose for entropy estimation. Instead, we reduce the nonparametric problem into a parametric submodel, and construct lower bound via the duality between moment matching and best approximation using {rational} functions. 

In addition to compactly supported densities, Theorem \ref{thm:main} can be extended to densities supported on $\reals^d$ with general tail conditions. Let $\Psi: [0,\infty] \to [0,\infty]$ be an Orlicz function, i.e., a continuous, increasing and convex function $\Psi$ satisfying $\Psi(0)=0$, $\Psi(u)>0$ for any $u>0$ and $\lim_{u\to\infty} \Psi(u)=\infty$. Moreover, we say $\Psi$ is of \emph{rapid growth} if there is a constant $\kappa=\kappa(\Psi)>1$ such that $\Psi(\kappa u)\ge \Psi(u)^2$ holds for all $u\ge 0$. Examples of rapidly growing Orlicz functions include $\Psi_q(u) = \exp(u^q)-1$ for any $q\geq 1$, with $\kappa(\Psi_q) = 2^{1/q}$; in particular, the cases of $q=1$ and $q = 2$ correspond to the sub-exponential and sub-Gaussian class, respectively.
 Consider the following class of densities:
\begin{align}\label{eq:Orlicz_tail}
\Lip_{s,p,d}^{\Psi}(L) \triangleq \{f: \|f\|_{\Lip_{s,p,d}} \le L\} \cap \left\{f:  \int_{\bR^d} \Psi\left(|x|\right)f(x)dx \le L \right\},
\end{align}
where $\|\cdot\|_{\Lip_{s,p,d}}$ is the Lipschitz norm defined in \eqref{eq:lipnorm}. Note that the second constraint of \eqref{eq:Orlicz_tail} implies that the $\Psi$-Orlicz norm of the random variable $|X|$ with $X\sim f$ is upper bounded by $L$. 

The following theorem presents the minimax rate for entropy estimation over $ \Lip_{s,p,d}^{\Psi}(L) $. 
\begin{theorem}[Densities with unbounded support]\label{thm:orlicz}
	Let $\Psi$ be an Orlicz function of rapid growth and $\Psi^{-1}$ its inverse function. 
	For any $d\in \mathbb{N}, 0<s\le 2, 2\le p<\infty$, 
	there exist constants $c,C,L_0>0$ depending on $s,p,d,\kappa(\Psi),\Psi(1)$ such that if $\Psi^{-1}(n)\ge 1$ and $L\ge L_0$, then
	\begin{align*}
	c\pth{ (n\ln n)^{-\frac{s}{s+d}} [\Psi^{-1}(n)]^{d\left(1-\frac{d}{p(s+d)}\right)} + n^{-\frac{1}{2}} } &\leq \left(\inf_{\hat{H}}\sup_{f\in \Lip_{s,p,d}^{\Psi}(L)} \bE_f (\hat{H}-H(f))^2 \right)^{\frac{1}{2}} \\
	&\leq C\pth{(n\ln n)^{-\frac{s}{s+d}} [\Psi^{-1}(n)]^{d\left(1-\frac{d}{p(s+d)}\right)} + n^{-\frac{1}{2}}}.
	\end{align*}
Moreover, the minimax lower bound works for any $s>0, 1\leq p<\infty$.
\end{theorem}

Comparing Theorem \ref{thm:orlicz} with Remark \ref{remark:radius}, we see that for general Orlicz function $\Psi$ with rapid growth, any density in $\Lip_{s,p,d}^{\Psi}(L)$ is effectively supported on $[-\Psi^{-1}(n), \Psi^{-1}(n)]^d$. There is also a subtle difference: the hidden constant in the parametric rate $\Theta(n^{-1/2})$ does not involve $\Psi^{-1}(n)$, thanks to the Orlicz norm constraint. Note that for simplicity we assume that $L$ is a constant and omit the dependence on $L$ in Theorem \ref{thm:orlicz}.

The estimator that achieves the minimax rates in Theorems~\ref{thm:main} and~\ref{thm:orlicz} relies on polynomial approximation. It is a natural question to ask whether an integral-form\footnote{Given a density estimate $\hat{f}$, an integral-form plug-in estimator for the entropy is $\int - \hat f(x)\ln \hat f(x)dx$, as opposed to $\frac{1}{n}\sum_{i=1}^n \log \hat f(X_i)$.} plug-in estimator using kernel density estimate can achieve the minimax rates. Recall that the kernel density estimator takes the form
\begin{align}\label{eq:KDE}
\hat{f}_h(x) = \frac{1}{nh^d}\sum_{i=1}^n K\left(\frac{x - X_i}{h}\right)
\end{align}
where $K(\cdot)$ is a kernel function, and $h$ is the bandwidth. The next result shows that the answer is negative for any sliding window kernel density estimator with a spatially invariant bandwidth (that is, the bandwidth $h$ can depend on the sample size $n$ but not on the location $x$):

\begin{theorem}[Suboptimality of integral-form plug-in estimators]
\label{lem.plug-in}
For $s\in (0,2], p\ge 2$, let $\hat{f}_h(x)$ be given in \eqref{eq:KDE} and define the integral-form plug-in estimator as 
$H(\hat f_h) = \int_{[0,1]^d} -\hat{f}_h(x)\ln \hat{f}_h(x)dx $. If the kernel $K(\cdot)$ satisfies \prettyref{ass:kernel} and $h\asymp (Ln)^{-1/(s+d)}$, then for $L\le n^{s/d}$, 
\begin{align*}
\left[\sup_{f\in \Lip_{s,p,d}(L)} \bE_f\left(H(\hat f_h)- H(f)\right)^2\right]^{\frac{1}{2}} \le C \left( n^{-\frac{s}{s+d}}L^{\frac{d}{s+d}} + n^{-\frac{1}{2}}\ln L \right), 
\end{align*}
where $C>0$ is a constant independent of $n,L$.

Conversely, for any kernel $K(\cdot)$ satisfying \prettyref{ass:kernel}
 and any bandwidth $h>0$, there exist constants $L_0>0, c>0$ independent of $n,L,h$, such that for any $L\ge L_0$, 
\begin{align*}
\left[\sup_{f\in \Lip_{s,p,d}(L)} \bE_f\left(H(\hat{f}_h) - H(f)\right)^2\right]^{\frac{1}{2}} \ge c\left(n^{-\frac{s}{s+d}}L^{\frac{d}{s+d}} + n^{-\frac{1}{2}}\ln L\right).
\end{align*}
\end{theorem}

Theorem~\ref{lem.plug-in} presents a tight characterization of the integral-form plug-in approach, and shows that the plug-in idea applied to the integral is strictly sub-optimal: the bias of the kernel-based plug-in estimator is $O(n^{-s/(s+d)}L^{d/(s+d)} )$, while for the optimal estimator it is $O((n\ln n)^{-s/(s+d)}L^{d/(s+d)} )$. 




Next we elaborate on the various assumptions in \Cref{thm:main}:

\vspace{1.5mm}
\noindent\emph{The Lipschitz ball $\Lip_{s,p,d}(L)$:} For $s=r+\alpha$ with $r$ integer and $\alpha\in (0,1]$, the H\"{o}lder ball $\text{H}_d^s(L)$ with smoothness $s$ and radius $L$ consists of functions $f$ with
$
\sup_{x\neq y} {|f^{(r)}(x) - f^{(r)}(y)|} / {|x-y|^{\alpha}} \le L.
$
The Lipschitz ball is a generalization of the H\"{o}lder ball by imposing the smoothness constraint \emph{on average} through the norm parameter $p$; for example, for $p=\infty$ the Lipschitz ball coincides with the H\"{o}lder $\Lip_{s,\infty,d}(L)=\text{H}_d^s(L)$ for any non-integer $s>0$.\footnote{As opposed to the definition of the Lipschitz ball in \prettyref{eq:lipnorm}, there is another slightly different definition using the modified Lipschitz norm $\Lip_{s,p,d}^{*}$, which coincides with a special case of the Besov ball $\text{B}_{p,\infty,d}^s$~\cite{Devore--Lorentz1993}. These two definitions are equivalent for non-integer $s$, while for integer $s$ the latter is strictly bigger. In this paper we adopt the former definition in \prettyref{eq:lipnorm} to avoid some technical subtleties.}

\vspace{1.5mm}
\noindent\emph{Radius of the Lipschitz ball:} The assumption $L\le (n\ln n)^{s/d}$ ensures that the minimax rate in Theorem \ref{thm:main} is $O(1)$, and $L\geq L_0$ is not superfluous as well. Indeed, if $sp \ge d$, then by standard embedding results of Lipschitz (or Besov) spaces \cite{jawerth1977some}, there exists $L_1=L_1(s,p,d)>0$ such that any $f \in \Lip_{s,p,d}(L_1)$ is bounded from below by a positive constant almost everywhere\footnote{In fact, \cite[Theorem 2.1]{jawerth1977some} states that if $sp\ge d$, $\|f(\cdot)-f(\cdot-t)\|_\infty \le C(s,p,d)\|f(\cdot)-f(\cdot-t)\|_{\Lip_{s,p,d}}$ for any $t\in [0,1]^d$. Since there must be some $x_0\in [0,1]^d$ such that $f(x_0)\ge 1$, we conclude that $f(x)\ge 1-2C(s,p,d)L$ for almost every $x\in [0,1]^d$, which is bounded from below by a constant if $L$ is sufficiently small.}, which, in view of the previous results \cite{robins2008higher,robins2017higher,mukherjee2017semiparametric}, implies that the entropy can be estimated at a faster rate $\Theta(n^{-{4s}/(4s+d)} + n^{-1/2})$ than that in \Cref{thm:main}. 

\vspace{1.5mm}
\noindent\emph{The smoothness condition $s\in (0,2]$:} Capturing high-order smoothness $s>2$ of a function is often challenging in nonparametric statistics, especially for density models. For example, if one would like to apply a kernel density estimator, for $s>2$ there does not exist a non-negative kernel to keep all polynomials with degree at most $\lfloor s\rfloor$. We will discuss this phenomenon in details in Section \ref{sec:highsmooth}. We note that the minimax lower bound $\Omega((n\ln n)^{-s/(s+d)}L^{d/(s+d)} + n^{-1/2}\ln L)$ only requires $0<s<\infty, 1\leq p<\infty$. 

\vspace{1.5mm}
\noindent\emph{The norm condition $p\in [2,\infty)$:} Our current upper bound requires $p \geq 2$, which ensures the difference between the entropy of the true density and its kernel-smoothed version is at the right order.
For the lower bound, the case of $p=\infty$ imposes a too strict constraint on the density (i.e., to be smooth everywhere), while $p<\infty$ only imposes an average-case smoothness constraint which can be handled by the current construction. When $p = \infty$ we prove a lower bound of $\Omega(n^{-{s/(s+d)}} (\ln n)^{-(s+2d)/(s+d)}L^{d/(s+d)} + n^{-1/2}\ln L)$ as shown in Theorem~\ref{thm.lower_bound}. 

\vspace{1.5mm}
\noindent\emph{The support of $f$:} For general nonparametric functional estimation problem, there are essentially three factors contributing to the minimax rates: the tail behavior if $f$ is supported on $\mathbb{R}^d$, the boundary behavior if $f$ is compactly supported, and the behavior of $f$ in the interior of its domain. In Theorem~\ref{thm:main}, we assume that $f$ is compactly supported and smoothly vanishing at the boundary so that sliding window kernel methods are applicable; this assumption is relaxed in Section \ref{sec.extension} to the so-called ``periodic boundary condition"~\cite{krishnamurthy2014nonparametric}. 
The effect of the tail behavior on the minimax rates is precisely quantified in Theorem~\ref{thm:orlicz} for densities with unbounded support.

%
%
%
%
%


\subsection{Related work}

%


The problem of estimating the entropy of a density has been investigated extensively in the literature. As discussed in the overview~\cite{beirlant1997nonparametric}, there exist two main approaches, based on either kernel density estimators, e.g.~\cite{joe1989estimation,gyorfi1991nonparametric,hall1993estimation,paninski2008undersmoothed,kandasamy2015nonparametric} or nearest neighbor methods, e.g.~\cite{tsybakov1996root, sricharan2012estimation,singh2016finite,delattre2017kozachenko,gao2018demystifying,berrett2019efficient}. Among these works, some focus on the consistency~\cite{gyorfi1991nonparametric,paninski2008undersmoothed}, $\sqrt{n}$-consistency~\cite{joe1989estimation,tsybakov1996root}, or the asymptotic efficiency~\cite{hall1993estimation,berrett2019efficient} of the proposed estimator, while others work on the minimax rate~\cite{sricharan2012estimation,kandasamy2015nonparametric,singh2016finite,delattre2017kozachenko,gao2018demystifying}. 

Similar estimator constructions have appeared in the literature. Asymptotic efficient estimators are obtained in \cite{Ibragimov--Nemirovskii--Khasminskii1987some,nemirovski2000topics} for smooth functionals by means of Taylor expansion; \cite{lepski1999estimation} and \cite{Cai--Low2011}
estimated the $L_1$ norm of the mean in Gaussian white noise model using trigonometric polynomial approximation. One related work~\cite{Han--Jiao--Mukherjee--Weissman2017adaptive} deserves special attention. 
Dealing with the Gaussian white noise model, \cite{Han--Jiao--Mukherjee--Weissman2017adaptive} analyzed the minimax rates of estimating the $L_r$ norms (for all $r\in [1,\infty)$) of the mean function over Besov spaces which was previously studied in \cite{lepski1999estimation}. Although both papers use the polynomial approximation technique for the upper and lower bound construction (which trace back to earlier work of \cite{lepski1999estimation,Cai--Low2011,Jiao--Venkat--Han--Weissman2015minimax,wu2016minimax}), there exist significant distinctions between this work and \cite{Han--Jiao--Mukherjee--Weissman2017adaptive}. First, here we analyze the density model as opposed to the location model, and it is crucial to design estimators to adapt to low-density regions. This specific problem has been investigated in \cite{patschkowski2016adaptation} for estimating linear functionals (density at a given point), where it was conjectured that the case of $s>2$ exhibits significantly different behavior from the case of $0<s\leq 2$; this is the underlying reason for the assumption $0<s\leq 2$ for our upper bound, which is discussed in more details in Section~\ref{sec:highsmooth}. In contract, in white noise models there is no need to adapt. Moreover, when $d=1$ these two models are asymptotically equivalent~\cite{nussbaum1996asymptotic} for $s>1/2$ provided that the density is bounded from below by a positive constant; however, they do \emph{not} imply the minimax rates of a given estimation problem for these two models must coincide, and for small densities the equivalence can break down~\cite{ray2018asymptotic}. In fact, in contrast to the conclusion of \prettyref{rmk:param}, it is shown in \cite{Han--Jiao--Mukherjee--Weissman2017adaptive} that the parametric rate is never achievable for ``entropy'' estimation in the white noise model. Second, the estimator construction in this paper requires more delicate analysis, and bounding the approximation error $H(f_h) - H(f)$ relies on a novel application of the Hardy-Littlewood maximal inequality in conjunction with the nonnegativity of the density function, which also leads to, as a by-product, a new inequality upper bounding the Fisher information in terms of the $L_p$ norm of the second derivative (\prettyref{thm:fisher}). Third, in the minimax lower bound, this work carefully chooses non-negative functions (not required in the Gaussian white noise model), and analyzes the total variation bound instead of the $\chi^2$-divergence bound which is simpler and more suitable for Gaussian models.

\subsection{Notation}
	\label{sec:notation}
	
	For a finite set $A$, let $|A|$ denote its cardinality. The norm $|\cdot|$ denotes the Euclidean norm of vectors in $\bR^d$, and $\|\cdot\|_p$ denotes the $L_p$ norm (with respect to the Lebesgue measure) of real-valued functions defined on $\bR^d$. 
		Let $\|\cdot\|_{\text{\rm op}}$ denotes the operator norm of matrices, \ie, the largest singular value.
		For $x\in\reals^d$, let $x_{\backslash i} \triangleq (x_j: j\neq i) \in \reals^{d-1}$. 
	For $n\in\naturals$, let $[n]\triangleq \{1,\ldots,n\}$. Denote by $\binom{[n]}{l} =\{	J\subseteq [n]: |J|=l\}$ the collection of all $l$-subsets of $[n]$.
	Throughout the paper,  for non-negative sequences $\{a_\gamma\}$ and $\{b_\gamma\}$, we write $a_\gamma \lesssim b_\gamma$ (or $a_n=O(b_n)$) if $a_\gamma \leq C b_\gamma$ for some positive constant $C$ that does \emph{not} depend on the sample size $n$, the bandwidth $h$, or the Lipschitz norm $L$. We use $a_\gamma \gtrsim b_\gamma$ (or $a_\gamma=\Omega(b_\gamma)$) to denote $b_\gamma \lesssim a_\gamma$, and $a_\gamma \asymp b_\gamma$ (or $a_\gamma=\Theta(b_\gamma)$) to denote both $a_\gamma \lesssim b_\gamma$ and $b_\gamma \lesssim a_\gamma$. We use $a_\gamma \ll b_\gamma$ (or $a_\gamma=o(b_\gamma)$) to denote $\lim_\gamma \frac{a_\gamma}{b_\gamma}=0$, and $a_\gamma \gg b_\gamma$ (or $a_\gamma=\omega(b_\gamma)$) to denote $b_\gamma \ll a_\gamma$. The support set of a probability measure $\mu$ is denoted by $\supp(\mu)$.
	Let $P_X$ denote the distribution of a random variable $X$.
	The KL (resp.~$\chi^2$) divergence  from distribution $\mu$ to $\nu$ is defined as 
	$D(\mu\|\nu) = \int d\mu\log \frac{d\mu}{d\nu}$ (resp.~$\chi^2(\mu\|\nu) = \int d\nu(\frac{d\mu}{d\nu}-1)^2$)
	if $\mu \ll \nu$ and $+\infty$ otherwise.
	
\subsection{Organization}
	\label{sec:org}
	
The rest of this paper is organized as follows. Section \ref{sec.construction} presents the construction of the minimax rate-optimal estimator. Section \ref{sec.analysis} proves the upper bound. In particular, the analysis of the bias incurred by the first-stage approximation relies on a novel application of the Hardy-Littlewood maximal inequality, and the same argument also leads to an inequality on Fisher information, which is presented at the end of \prettyref{sec:firststage} and might be of independent interest. Section \ref{sec.discussion} discusses generalizations and open problems.
In particular, Section \ref{sec.extension} extends the results to a broader class of densities that satisfy a periodic boundary conditions, and establishes the corresponding minimax rates of entropy estimation. 
 Remaining proofs are relegated to the appendices.

\section{Construction of the estimator}\label{sec.construction}

Define the smoothed density
\begin{align*}
f_h(x) \triangleq \int_{\reals^d} K_h(x-y)f(y)dy,
\end{align*}
where $K_h(\cdot)$ is some kernel function with bandwidth $h>0$. In the special case of $K_h(x) = \frac{1}{h^d} K(\frac{x}{h})$ for some kernel function $K: \mathbb{R}^d\to\reals$, 
we have
	\begin{align}
	f_h(x) = \int_{\bR^d} \frac{1}{h^d} K\pth{\frac{x-y}{h}} f(y)dy, 
	\label{eq.f_h}
	\end{align}
which admits the following natural unbiased estimator (kernel density estimate)
\begin{align}\label{eq.f_h_hat}
\hat{f}_h(x) \triangleq \frac{1}{n}\sum_{i=1}^n K_h(x-X_i), 
\end{align}
where $X_1,\cdots, X_n\iiddistr f$.

The optimal estimator for the entropy $H(f)$ is constructed in two steps:
First, 
by choosing a suitable (in particular, compactly-supported) kernel $K$, we approximate $f$ by $f_h$ and bound $|H(f_h) - H(f)|$ using functional-analytic properties of the density class. 
Next, we construct an estimator for $H(f_h)$ based on the kernel density estimator $\hat{f}_h$. 
The main insight is that  $\{\hat{f}_h(x): x \in [0,1]^d\}$ is essentially a \emph{finite-dimensional} parametric model, in the sense that $\hat{f}_h(x)$ roughly follows the binomial distribution $nh^d\hat{f}_h(x)\sim \mathsf{B}(n,h^df_h(x))$ (cf.~Lemma \ref{lemma.regime}). 
As a result, we essentially obtain a parametric binomial model with $h^{-d}$ parameters, so that the existing approximation-theoretic techniques for entropy estimation in parametric models \cite{wu2016minimax,Jiao--Venkat--Han--Weissman2015minimax} can be applied. 
We now describe the construction of the optimal entropy estimator for any $s\in (0,2]$ in any dimension.  For the first approximation stage, in order to find a suitable approximation $f_h$ for $f$, we recall the following property of Lipschitz spaces~\cite[Theorem 8.1]{hardle2012wavelets}:
\begin{lemma}
\label{lemma.lipschitz_approx}
	Fix any $s>0$ and any kernel $K: \mathbb{R}^d\to\reals$ which satisfies $\int_{\mathbb{R}^d} |x|^{\lceil s\rceil}|K(x)|dx<\infty$ and 
	maps any polynomial $q$ in $d$ variables of degree at most $\lceil s\rceil-1$ to themselves, i.e., 
	$\int_{\mathbb{R}^d} q(y) K(x-y)dy=q(x)$. 
	Then for any $f\in \Lip_{s,p,d}(L)$ and $f_h$ defined in \prettyref{eq.f_h}, 
	we have
	\begin{align*}
	\|f_h-f\|_p \triangleq \pth{	\int_{\bR^d} |f_h(x)-f(x)|^p dx}^{1/p} \lesssim Lh^s.
	\end{align*}
\end{lemma}

To apply Lemma \ref{lemma.lipschitz_approx}, we choose a kernel $K$ with the following properties:
\begin{assumption}
\label{ass:kernel}
Suppose $K:\reals^d\to\reals$ satisfies the following:
	\begin{enumerate}
\item Non-negativity: $K(t)\ge 0$ for any $t\in \mathbb{R}^d$;
\item Unit total mass: $\int_{\mathbb{R}^d} K(t)dt=1$;
\item Zero mean: $\int_{\mathbb{R}^d} tK(t)dt=0$;
\item Finite second moment: $\int_{\mathbb{R}^d} |t|^2K(t)dt < \infty$. 
\item Compact support: $\sup\{|t|: K(t)\neq 0\} <\infty$. 
\end{enumerate}
\end{assumption}

There are several kernel functions which fulfill \prettyref{ass:kernel}, e.g., the box kernel $K(t)=\mathbbm{1}(t\in [-1/2,1/2]^d)$. Note that the second and the third requirements ensure that it keeps all polynomials of degree at most one, i.e., 
$\int_{\reals^d} (a^\top x+b) K(x-y) dx = a^\top y+b$, 
and the first requirement (non-negativity) is crucial for proving the concentration result in Section \ref{sec.analysis}. In fact, the non-negativity requirement is the key reason why we need to impose the assumption $s\le 2$, and relaxing this requirement appears highly challenging (cf.~Section \ref{sec:highsmooth}). 

Since the kernel $K(\cdot)$ has a compact support, the approximation $f_h$ is compactly supported as well. By Lemma \ref{lemma.lipschitz_approx} and our assumption that $s\le 2$, we have
\begin{align}\label{eq:ffh}
\|f-f_h\|_p \lesssim Lh^s.
\end{align}
Later in \prettyref{sec:firststage} we will show that the entropy difference also satisfies $|H(f)-H(f_h)| \lesssim Lh^s$.



Fix an appropriate kernel $K$ that fulfills \Cref{ass:kernel} and define $f_h, \hat{f}_h$ as in \eqref{eq.f_h} and \eqref{eq.f_h_hat}. 
To construct an estimator $\hat H$ for $H(f_h) =  \int -f_h(x)\ln f_h(x) dx$, we let $\hat H = \int \hat H(x) dx$, where for each $x\in[0,1]^d$, $\hat H(x)$ is an estimator for $-f_h(x)\ln f_h(x)$ obtained as follows:

\begin{enumerate}
\item 
For notational convenience, let the sample size be $3n$ as opposed to $n$. Split the observations into three parts $X^{(1)}, X^{(2)}, X^{(3)}$, each consisting of $n$ observations.


\item For each part of observations, construct the kernel density estimators $\hat{f}_{h,1}(x), \hat{f}_{h,2}(x)$ and $\hat{f}_{h,3}(x)$ per \eqref{eq.f_h_hat}. The estimator $\hat{f}_{h,1}(x)$ will be used for classifying smooth versus non-smooth regime, and the other two for estimation. 
\item Regime classification and estimator construction:
\begin{itemize}
\item ``Non-smooth" regime: $\hat{f}_{h,1}(x)<\frac{c_1\ln n}{nh^d}$. Denote by $Q$ the best degree-$k$ polynomial approximation of $-t\ln t$ on $[0,\frac{2c_1\ln n}{nh^d}]$:
\begin{align}
\label{eq:Q}
Q = \sum_{l=0}^k a_lt^l = \arg\min_{P\in \mathsf{Poly}_k} \max_{t\in [0,\frac{2c_1\ln n}{nh^d}]} |-t\ln t-P(t)|,
\end{align}
where $\mathsf{Poly}_k$ denotes the collection of all polynomials of degree at most $k$.
Define the following unbiased estimator of $Q(f_h(x))$ in terms of $U$-statistics:
\begin{align}
\hat{H}_1(x) = \sum_{l=0}^k a_l\left(\binom{n}{l}^{-1}\sum_{J\in\binom{[n]}{l}}\prod_{j\in J} K_h(x-X_{j}^{(2)})\right).
\label{eq:H1x}
\end{align}

\item ``Smooth" regime: $\hat{f}_{h,1}(x)\ge\frac{c_1\ln n}{nh^d}$. Define the following bias-corrected plug-in estimator:
\begin{equation}\label{eq:H2x}
\begin{aligned}
\hat{H}_2(x) &= \mathbbm{1}(\hat{f}_{h,2}(x)\ge \frac{c_1\ln n}{4nh^d})\cdot
\Bigg\{-\hat{f}_{h,2}(x)\ln \hat{f}_{h,2}(x)  \\
&- (1+\ln\hat{f}_{h,2}(x))(\hat{f}_{h,3}(x)-\hat{f}_{h,2}(x))
-\frac{1}{2}\Bigg(\hat{f}_{h,2}(x) \\ 
&-2\hat{f}_{h,3}(x)+\frac{1}{\binom{n}{2}\hat{f}_{h,2}(x)}\sum_{i<j}K_h(x-X_{i}^{(3)})K_h(x-X_{j}^{(3)})\Bigg)\Bigg\}. 
\end{aligned}
\end{equation}
\item The final point estimate of $H(x)=-f_h(x)\ln f_h(x)$ is
\begin{equation}\label{eq:hatHx}
\begin{aligned}
\hat{H}(x) &\triangleq \min\sth{\hat{H}_1(x),\frac{1}{n^{1-2\varepsilon}h^d}}\mathbbm{1}\left(\hat{f}_{h,1}(x)<\frac{c_1\ln n}{nh^d}\right) \\
&\qquad+ \hat{H}_2(x)\mathbbm{1}\left(\hat{f}_{h,1}(x)\ge \frac{c_1\ln n}{nh^d}\right).
\end{aligned}
\end{equation}

\end{itemize}
\end{enumerate}

Finally, choose 
\begin{equation}
h = c_0(Ln\ln n)^{-\frac{1}{s+d}}, \quad k = \lceil{c_2 \ln n}\rceil,
\label{eq:parameters}
\end{equation}
where $c_0>0$ is any constant, $0<7c_2\ln 2< \varepsilon < \frac{s}{s+d}$ and 
$c_1>0$ is sufficiently large (per Lemma~\ref{lemma.non-smooth}--\ref{lemma.smooth}) 
 and output the estimator
\begin{align}
\label{eq:Ihat}
\hat{H} = \int_{\bR^d} \hat{H}(x)dx.
\end{align}
Note that the integration only need to be taken over the support of $f_h$, which is slightly larger than the unit cube $[0,1]^d$. 
This completes the construction of our estimator. A few remarks are in order:
\paragraph{Choice of the $U$-statistics} 
The following $U$-statistic
\begin{align*}
U_m = \frac{1}{\binom{n}{m}}\sum_{1\le i_1<i_2<\cdots<i_m\le n}\prod_{j=1}^m K_h(x-X_{i_j})
\end{align*}
has appeared several times in the estimator construction, which is the natural unbiased estimator for powers of $f_h(x)$:
\begin{align*}
\bE[U_m] 
&=\frac{1}{\binom{n}{m}}\sum_{1\le i_1<i_2<\cdots<i_m\le n}\prod_{j=1}^m \bE\left[K_h(x-X_{i_j})\right] =f_h(x)^m.
\end{align*}
The reason why we average over all possible subsets of size $m$ is to reduce the variance to the correct order (cf.~\Cref{lemma.U-statistic}). In practice, to compute the $k$-th order $U$-statistics, note that it is simply the (normalized) $k$-th elementary symmetric polynomial of $K_h(x-X_{i})$. Hence, it suffices to compute the power sum
$
\sum_{i=1}^n (K_h(x-X_{i}))^l
$
for all $l=1,\cdots,k$, and then invoke Newton's identity to compute elementary symmetric polynomials; this has overall time complexity $O(nk+k^2)=O(n\log n)$. For the special case of the box kernel $K(t) = \mathbbm{1}(t\in [-1/2,1/2]^d)$, which can be used to achieve the upper bound in \Cref{thm.achievability},
$\hat{H}_1(x)$ reduces to 
\begin{align}\label{eq:Zx}
	\hat{H}_1(x) = \sum_{l=0}^k a_l \cdot \frac{Z_x\cdot(Z_x-1)\cdot \ldots \cdot (Z_x-l+1)}{h^{ld} \cdot n\cdot(n-1)\cdot \ldots \cdot (n-l+1)},
\end{align}
where $
	Z_x = \sum_{i = 1}^n h^{d}K_h(x-X^{(2)}_i). 
$
Hence, the computational cost can be further reduced to $O(n+k^2)=O(n)$ in this simple example.

\paragraph{Polynomial approximation in the non-smooth regime} In the non-smooth regime (i.e., $\hat{f}_{h,1}(x)\le \frac{c_1\ln n}{nh^d}$) a suitable linear combination of the $U$-statistics is applied, where the coefficients come from the best approximating polynomial of our target functional $-x\ln x$. By the previous property of the $U$-statistic, in the non-smooth regime we estimate $Q(f_h(x))$ without any bias, and thus the bias in this regime becomes the polynomial approximation error. The coefficients of the polynomial $Q(\cdot)$ can be efficiently computed via the Remez algorithm, which converges double exponentially fast (see discussions in~\cite{Jiao--Venkat--Han--Weissman2015minimax}). The coefficients can also be pre-computed and stored so that there is no need to recompute the coefficients when applying the estimator.

\paragraph{Bias correction based on Taylor expansion in the smooth regime}
In the smooth regime (i.e., $\hat{f}_{h,1}(x)> \frac{c_1\ln n}{nh^d}$), we use the idea in \cite{han2016minimax} to correct the bias. Specifically, by Taylor expansion we can write
\begin{equation}\label{eq:debias-taylor}
\begin{aligned}
\phi(f_h(x)) &\approx \sum_{l=0}^R \frac{\phi^{(l)}(\hat{f}_{h,2}(x))}{l!}(f_h(x)-\hat{f}_{h,2}(x))^l\\
&=\sum_{l=0}^R\frac{\phi^{(l)}(\hat{f}_{h,2}(x))}{l!}\sum_{j=0}^l \binom{l}{j}f_h(x)^j(-\hat{f}_{h,2}(x))^{l-j} .
\end{aligned}
\end{equation}
A natural idea to de-bias is to find an unbiased estimator of the right-hand side in \cref{eq:debias-taylor}. Indeed, this can be done by sample splitting: we can split observations to obtain $\hat{f}_{h,3}(x)$, an independent copy of $\hat{f}_{h,2}(x)$, and then apply the previous $U$-statistics to $\hat{f}_{h,3}(x)$ to obtain an unbiased estimator $f_h(x)^j$. Our estimator construction uses this idea with $\phi(z)=-z\ln z$ and $R=2$ (which suffices for our de-biasing purposes).

\paragraph{Choice of bandwidth}
As will be clarified in Section \ref{sec.analysis} (cf.~\prettyref{eq:overall}), the bandwidth $h\asymp(n\ln n)^{-1/(s+d)}$ in \eqref{eq:parameters} is chosen in order to balance between two types of biases of our estimator. Compared with the optimal bandwidth $h\asymp n^{-1/(2s+d)}$ in estimating the density under $L_p$ risk for $p\in[1,\infty)$ \cite{nemirovski2000topics}, our choice of the bandwidth results in an ``undersmoothed" kernel estimator of the density, which is consistent with the findings in \cite{goldstein1992optimal,paninski2008undersmoothed} that an undersmoothed kernel estimator should be used in estimating nonparametric functionals. However, our specific choice of $h\asymp(n\ln n)^{-1/(s+d)}$ is different from the optimal bandwidths for other problems, such as $h\asymp n^{-2/(4s+d)}$ for estimating quadratic, cubic, and general smooth functionals \cite{bickel1988estimating,kerkyacharian1996estimating,robins2008higher,mukherjee2017semiparametric}, and the optimal bandwidth $h\asymp (n\ln n)^{-1/(2s+d)}$ for estimating the $L_r$ norm of the function in Gaussian white noise in one dimension with $r\in [1,\infty)$ not an even integer \cite{lepski1999estimation,Han--Jiao--Mukherjee--Weissman2017adaptive}.


\paragraph{Final integration}
The estimator $\hat{H}(x)$ in \prettyref{eq:hatHx} provides the pointwise estimation of $-f_h(x)\ln f_h(x)$ for all $x\in \supp(f_h)$, and an integration is required to produce the final entropy estimator. If the box kernel is used, notice that $n$ small cubes of equal size can partition the unit cube $[0,1]^d$ into $O(n^d)$ pieces, the mapping $x\mapsto Z_x$ in \eqref{eq:Zx} is piecewise constant on $O(n^d)$ pieces. Hence, for exact integration it suffices to evaluate $\hat{H}(x)$ at $O(n^d)$ points, which yields an overall $O(n^{d+1}\log n)$ time complexity of our estimator.
For practical implementation with general kernels, numerical integration methods and quadrature formulas can be used to evaluate the integral, and then we only need to evaluate $\hat{H}(x)$ at finitely many points.

In the next section we prove the following result, which completes the proof of the upper bound in Theorem \ref{thm:main}. 
\begin{theorem}\label{thm.achievability}
For $s\in (0,2]$, $p\ge 2$ and $L\le (n\ln n)^{s/d}$, the following holds for the estimator $\hat{H}$ defined in \cref{eq:Ihat}:
\begin{align*}
\left(\sup_{f\in \Lip_{s,p,d}(L)} \mathbb{E}_f(\hat{H}-H(f))^2\right)^{\frac{1}{2}} \le C\left((n\ln n)^{-\frac{s}{s+d}}L^{\frac{d}{s+d}} + n^{-\frac{1}{2}}\ln L\right), 
\end{align*}
where $C=C(s,p,d)>0$ is independent of $n, L$ (we omit the dependence of $C$ on the choice of parameters $c_0,c_1,c_2,\varepsilon$ and the kernel $K(\cdot)$).
\end{theorem}

\section{Proof of upper bound}\label{sec.analysis}
The error of our estimator $\hat{H}$ can be decomposed into three terms: the approximation error of $|H(f_h)-H(f)|$, the bias and the variance of the estimation error of $\hat{H}$ in estimating $H(f_h)$. Next we deal with these terms separately.

\subsection{First-stage approximation error}
\label{sec:firststage}
The approximation error between $H(f_h)$ and $H(f)$ is summarized in the following lemma, which is one of the key results in this paper. 
\begin{lemma}\label{lemma.first_approx_err}
Let $s\in (0,2]$, $p\ge 2$.
For any $f\in \Lip_{s,p,d}(L)$ and bandwidth $h$ with $0<Lh^s\le 1$, let $f_h$ be defined in \eqref{eq.f_h}. There exists a constant $C>0$ independent of $h$, such that
\begin{align}
\label{eq:HH}
|H(f)-H(f_h)| = \left|\int_{\bR^d} f(x)\ln f(x) dx - \int_{\bR^d} f_h(x)\ln f_h(x)dx\right|\le C\cdot Lh^s
\end{align}
whenever $0<h<h_0$, where $h_0$ is a constant depending only on $s$.
\end{lemma}

In view of the fact that $\|f-f_h\|_p\lesssim Lh^s$ (cf.~(\ref{eq:ffh})), Lemma \ref{lemma.first_approx_err} essentially says that the entropy functional $H(\cdot)$ is ``Lipschitz" with respect to convolution. The proof of Lemma \ref{lemma.first_approx_err} consists of three steps:
\begin{enumerate}
	\item By the convolution property, we first express the entropy difference $H(f_h)-H(f)$ as a mutual information term. Then using the variational representation of the mutual information and $\chi^2$-divergence, we reduce \prettyref{eq:HH} to an inequality that no longer involves the kernel;
\item By the equivalence between the $K$-functional and the modulus of continuity \cite{Devore--Lorentz1993}, we approximate $f$ by a non-negative $C^2$-function $g$ and further reduce the goal to an estimate 
of the form  
\[
\int_{[0,1]^d} \frac{|\nabla g(x)|^2}{f(x)+h^s}dx;
\]

\item We invoke the Hardy-Littlewood maximal inequality to control the above integral using the $L_2$-norm of the second-order derivative of $g$. This is the crux of the proof.
The same proof technology also leads to a new upper bound on Fisher information, which we summarize at the end of this subsection.
\end{enumerate}

\subsubsection{Mutual information and $\chi^2$-divergence}
	\label{sec:MIchi}
Recall the mutual information between random variables $A$ and $B$ is defined as 
the KL divergence between the joint distribution and product of the marginal distributions:
\[
I(A;B) = D(P_{AB}\|P_A \otimes P_B) = \mathbb{E}\left[ \ln \frac{dP_{AB}}{dP_A dP_B} \right].
\] 
Recall that, by \prettyref{ass:kernel}, the kernel satisfies $K\geq 0$ and $\int_{\bR^d} K(x)dx=1$.
Let $X$ and $U$ be independent random variables with density function $f$ and $K$, respectively. 
Then by the convolution property, the density of $X+hU$ is $f_h$, and, as a result, 
\begin{equation}
0 \leq H(f_h)-H(f) = I(U;X+h U).
\label{eq:HI1}
\end{equation}
Note that by the compact support of the kernel $K$, the density $f_h$ is supported on a cube slightly larger than $[0,1]^d$ (i.e., with edge size $1+O(h)$), and by a proper scaling we assume without loss of generality that both $f$ and $f_h$ are supported on $[0,1]^d$.

Next we reduce the desired inequality into a simpler one independent of the kernel $K(\cdot)$. 
Let $w$ be an arbitrary density supported on $[0,1]^d$. Then
\begin{align}
 I(U; X+hU) 
&= \mathbb{E}_U \qth{\int_{[0,1]^d} f(x-hU)\ln \frac{f(x-hU)}{f_h(x)}dx} \nonumber \\
&= \mathbb{E}_U \qth{\int_{[0,1]^d} f(x-hU)\ln \frac{f(x-hU)}{w(x)}dx} - D(f_h \|w) \nonumber \\
&\stepa{\le} \mathbb{E}_U \qth{\int_{[0,1]^d} f(x-hU)\ln \frac{f(x-hU)}{w(x)}dx} \label{eq:varMI} \\
&\stepb{\le} \mathbb{E}_U \qth{\int_{[0,1]^d} \frac{(f(x-hU)-w(x))^2}{w(x)} dx} \nonumber
\end{align}
where (a) follows from the non-negativity of the KL divergence, and (b) is due to the fact that the KL divergence is upper bounded by the $\chi^2$-divergence. 

Since $Lh^s\le 1$, there exists another density $w$ on $[0,1]^d$ such that $w(x)\ge \max\{f(x)/2, Lh^s/4\}$ for all $x\in [0,1]^d$ and $\|w-f\|_\infty\le Lh^s/4$. Such an existence may be constructed by adding $Lh^s/4$ to $f(x)$ for all $x$ with $f(x)\le Lh^s/4$, and then subtracting $Lh^s/4$ from $f(x)$ for a subset of $\{x\in [0,1]^d: f(x)\ge Lh^s/2\}$ to arrive at a density. As a result,
\begin{align}
\int_{[0,1]^d} \frac{(f(x-hU)-w(x))^2}{w(x)} dx &= \int_{[0,1]^d} \frac{(f(x-hU)-f(x) + f(x)-w(x))^2}{w(x)} dx \nonumber\\
&\le  2\int_{[0,1]^d} \frac{(f(x-hU)-f(x))^2}{w(x)} dx + 2\int_{[0,1]^d} \frac{(f(x)-w(x))^2}{w(x)} dx \nonumber\\
&\le 2\int_{[0,1]^d} \frac{(f(x-hU)-f(x))^2}{w(x)} dx + \frac{2\|f-w\|_\infty^2}{Lh^s/4} \nonumber\\
&\lesssim \int_{[0,1]^d} \frac{(f(x-hU)-f(x))^2}{w(x)} dx + Lh^s. \label{eq:HI3}
\end{align}
Combining \prettyref{eq:HI1}--\prettyref{eq:HI3}, we have
\[
0 \leq H(f_h)-H(f) \lesssim \Expect_U\qth{\int_{[0,1]^d} \frac{(f(x-hU)-f(x))^2}{ f(x) + Lh^s} dx }+ Lh^s.
\]
Recall that the finite second moment of the kernel $K$ ensures that $\bE |U|^2<\infty$.
Therefore, for Lemma \ref{lemma.first_approx_err} to hold, it suffices to prove that for any $u\in \reals$, 
\begin{align}\label{eq.target}
\int_{[0,1]^d} \frac{(f(x+hu)-f(x))^2}{f(x)+Lh^s} dx\lesssim Lh^s(1+|u|^2).
\end{align}
Note that $\eqref{eq.target}$ no longer involves the kernel $K(\cdot)$. 

We provide some insights why \eqref{eq.target} is expected to hold. When $s\le 1$ and $p=\infty$, the Lipschitz ball condition ensures that $|f(x+hu) - f(x)|\lesssim Lh^s|u|^s\le Lh^s(1+|u|)$, and \eqref{eq.target} clearly holds. However, when $1<s\le 2$, we will only have $|f(x+hu) - f(x)|\lesssim Lh|u|$ in general, and \eqref{eq.target} cannot be derived by this simple approach. The crux of proving \eqref{eq.target} for $1<s\le 2$ is that, when $f(x)$ is close to zero, the difference $|f(x+hu) - f(x)|$ also need to be small to maintain the non-negativity of $f(x-hu)$. In Section \ref{sec:HL}, we will essentially show that $|f(x+hu) - f(x)|\lesssim L\sqrt{f(x)h^s}  (1+|u|)$, which leads to \eqref{eq.target}. 

\subsubsection{Approximation by $C^2$ functions} 
We need the following lemma to replace $f$ with a smoother function $g$:
\begin{lemma}\label{lemma.smooth_approx_g}
	Let  $f\in \Lip_{s,p,d}(L)$ be a non-negative function, with $s\in (0,2]$ and $p\ge 1$. 
	Then there exists $C=C(s,p,d)$, such that for any $h>0$, there exists a non-negative function $g\in C^2(\mathbb{R}^d)$ such that
	\begin{align}
	\|f-g\|_p & \leq CL h^s, \label{eq:approxg1}\\
	\|\|\nabla^2 g(\cdot)\|_{\text{\rm op}}\|_p & \leq CL h^{s-2}.
	\label{eq:approxg2}
	\end{align}
\end{lemma}

Note that Lemma \ref{lemma.smooth_approx_g} is essentially the equivalence between the $K$-functional and the modulus of smoothness (\Cref{lmm:Kequiv} in Appendix \ref{app:aux}), with an extra constraint that $g$ being non-negative, which turns out to be crucial in proving the inequality \eqref{eq.new_target} below.

Let $g$ be given by \prettyref{lemma.smooth_approx_g}. Then
\begin{align}
& \int_{[0,1]^d} \frac{(f(x+hu)-f(x))^2}{f(x)+Lh^s}dx \nonumber \\
&\le 3\int_{[0,1]^d} \frac{(f(x+hu)-g(x+hu))^2 + (f(x)-g(x))^2 + (g(x+hu)-g(x))^2}{f(x)+Lh^s}dx \nonumber\\
&\le 3\int_{[0,1]^d} \frac{(f(x+hu)-g(x+hu))^2 + (f(x)-g(x))^2}{Lh^s}dx \nonumber\\
&\qquad + 6\int_{[0,1]^d} \frac{(h\nabla g(x)^\top u)^2 + (g(x+uh)-g(x)-h\nabla g(x)^\top u)^2}{f(x)+Lh^s}dx \nonumber\\
&\le \frac{6\|f-g\|_2^2}{Lh^s} + \frac{6}{Lh^s}\|\rho(\cdot,u)\|_2^2 + 6h^2|u|^2 \int_{[0,1]^d} \frac{|\nabla g(x)|^2}{f(x)+Lh^s}dx,
\label{eq:fish}
\end{align}
where 
\[
\rho_h(x,u) \triangleq g(x+uh)-g(x)-h\nabla g(x)^\top u.
\]

We bound the three terms in \cref{eq:fish} separately. By \eqref{eq:approxg1}, the first term is upper bounded by
\begin{align*}
\frac{6\|f-g\|_2^2}{Lh^s} \le \frac{6\|f-g\|_p^2}{Lh^s} \lesssim Lh^s.
\end{align*}
For the second term, by the integral representation of the Taylor remainder term, we have
\begin{align*}
|\rho_h(x,u)| &= \left|\int_0^1 (1-t)u^\top\nabla^2 g(x+t\cdot hu)u\cdot dt\right| 
\\
&\leq  h^2|u|^2 \cdot \int_0^1 (1-t)\|\nabla^2 g(x+t\cdot hu)\|_{\text{op}} dt
\end{align*}
and hence
\begin{align*}
\|\rho_h(\cdot,u) \|_2 
&\le  h^2|u|^2\cdot \left\|\int_0^1 (1-t)\|\nabla^2 g(x+t\cdot hu)\|_{\text{op}} dt \right\|_2 \\
&\stepa{\le}  h^2|u|^2\cdot \int_0^1 (1-t) \big\| \|\nabla^2 g(x+t\cdot hu)\|_{\text{op}} \big\|_2 dt \\
&\stepb{\lesssim}  h^2|u|^2\cdot Lh^{s-2} = Lh^s|u|^2,
\end{align*}
where (a) follows from the convexity of norms and (b) follows from \cref{eq:approxg2}.
Thus the first two terms in \cref{eq:fish} are both upper bounded by $O(Lh^s|u|^2)$. Hence, to show \eqref{eq.target}, it remains to prove that
\begin{align}\label{eq.new_target}
\int_{[0,1]^d} \frac{|\partial_i g(x)|^2}{f(x)+Lh^s}dx \lesssim Lh^{s-2}, \qquad \forall i\in [d]
\end{align}
where $\partial_i g=\frac{\partial g}{\partial x_i}$.

\subsubsection{Application of the Hardy-Littlewood maximal inequality}
	\label{sec:HL}
Finally, we use the non-negativity of $g$ and the Hardy-Littlewood maximal inequality \cite{stein2016singular} to prove \eqref{eq.new_target}. Fix any $\tau>0$ to be optimized later. Since $g$ is non-negative, we have
\begin{align*}
0 &\le g(x+\tau e_i) \\
&= g(x) + \tau\cdot \partial_i g(x) + (g(x+\tau e_i)-g(x) - \tau\cdot \partial_i g(x))
\end{align*}
and thus 
\begin{align*}
-\tau\cdot \partial_i g(x) \le g(x) + (g(x+\tau e_i)-g(x) - \tau\cdot \partial_i g(x)).
\end{align*}
Replacing $x+\tau e_i$ by $x-\tau e_i$, we also have
\begin{align*}
\tau\cdot \partial_i g(x)\le g(x) + (g(x-\tau e_i)-g(x) + \tau\cdot \partial_i g(x)).
\end{align*}
Combining these two inequalities, we arrive at the following pointwise bound:
\begin{align*}
\tau \cdot |\partial_i g(x)| &\le 2g(x) + |g(x+\tau e_i)-g(x) - \tau\cdot \partial_i g(x)| \\
&\qquad + |g(x-\tau e_i)-g(x) + \tau\cdot \partial_i g(x)| \\
&\le 2g(x) + \tau^2 \int_{-1}^{1} |\partial_{ii}g(x+t\cdot \tau e_i)| dt
\end{align*}
where for the second inequality we have used the integral representation of the Taylor remainder term again.

Since the previous inequality holds for any $\tau>0$, we choose $\tau=\tau_x=\sqrt{h^{2-s}f(x)/L+h^2}$ to obtain an upper bound on the derivative:
\begin{align*}
|\partial_i g(x)| \le \frac{2g(x)}{\sqrt{h^{2-s}f(x)/L+h^2}} + \tau_x \int_{-1}^{1} |\partial_{ii}g(x+t\cdot \tau_x e_i)| dt.
\end{align*}
Plugging this bound into \eqref{eq.new_target} and using the triangle inequality, we have
\begin{align*}
&\int_{[0,1]^d} \frac{|\partial_i g(x)|^2}{f(x)+Lh^s}dx \\
&\le 2\Big(
\underbrace{Lh^{s-2} \int_{[0,1]^d} \frac{4g(x)^2}{(f(x)+Lh^s)^2}dx}_{\triangleq A_1} + \underbrace{ (Lh^{s-2})^{-1}  \int_{[0,1]^d} \left(\int_{-1}^{1} |\partial_{ii}g(x+t\cdot \tau_x e_i)| dt\right)^2 dx}_{\triangleq A_2}\Big). 
\end{align*}

Next we upper bound $A_1$ and $A_2$ separately. For $A_1$, we use the triangle inequality again to obtain
\begin{align*}
A_1 &= Lh^{s-2}\cdot \int_{[0,1]^d} \frac{4g(x)^2}{(f(x)+Lh^s)^2}dx \\
&\le 8Lh^{s-2}\cdot \int_{[0,1]^d} \frac{(g(x)-f(x))^2+f(x)^2}{(f(x)+Lh^s)^2}dx \\
&\le 8Lh^{s-2}\cdot \left(\int_{[0,1]^d} \frac{(g(x)-f(x))^2}{L^2h^{2s}}dx + \int_{[0,1]^d} \frac{(f(x))^2}{(f(x))^2}dx\right) \\
&= 8Lh^{s-2}\cdot \left(\frac{\|g-f\|_2^2}{L^2h^{2s}}+1\right) \lesssim Lh^{s-2}.
\end{align*}

Hence, it remains to upper bound $A_2$, and it further suffices to prove that for any $x_{\backslash i}\in [0,1]^{d-1}$,
\begin{align}\label{eq.maximal}
\int_{0}^1 \left(\int_{-1}^{1} |\partial_{ii}g(x+t\cdot \tau_x e_i)| dt\right)^2 dx_i \le C\int_0^1 |\partial_{ii}g(x)|^2dx_i
\end{align}
for some  constant $C>0$. In fact, if \eqref{eq.maximal} holds, then integrating both sides over $x_{\backslash i}\in [0,1]^{d-1}$ together with the fact $\|\partial_{ii}g\|_2 \le \|\partial_{ii}g\|_p\lesssim Lh^{s-2}$ completes the proof of \eqref{eq.new_target}.

The proof of \eqref{eq.maximal} requires the introduction of the maximal inequality. Fix any $x_{\backslash i}\in [0,1]^{d-1}$ and define $h(y)\triangleq |\partial_{ii}g(x_{\backslash i},y)|$, \eqref{eq.maximal} is equivalent to 
\begin{align}\label{eq.maximal_red}
\int_0^1 \left(\frac{1}{2\tau_x}\int_{x-\tau_x}^{x+\tau_x} h(y) dy\right)^2 dx \le \frac{C}{4}\int_0^1 |h(x)|^2dx.
\end{align}

For any function $h$ on the real line, recall the Hardy--Littlewood maximal function $M[h]$ is defined as
\begin{align}
\label{eq:mf}
M[h](y) \triangleq \sup_{t>0} \frac{1}{2t}\int_{y-t}^{y+t} |h(z)|dz.
\end{align}
Next we recall the maximal inequality on the real line~\cite{stein2016singular}:
\begin{lemma}\label{lemma.maximal_inequality}
	For any non-negative real-valued measurable function $h$ on the real line $\bR$, the following tail bound holds: for any $t>0$, there exists a universal constant $C_1>0$ such that for $p>1$ we have	
	\begin{align*}
	\|M[h]\|_p \le C_1\left(\frac{p}{p-1}\right)^{\frac{1}{p}}\|h\|_p.
	\end{align*}
\end{lemma}

Applying this lemma with $p=2$ yields \eqref{eq.maximal_red} (and thus \eqref{eq.maximal}), as desired, completing the proof of \Cref{lemma.first_approx_err}.

We finish this subsection by noting that the proof technology developed based on the maximal inequality in fact leads to the following upper bound on Fisher information, which may be of independent interest.
\begin{theorem}
\label{thm:fisher}
	Let $f \in C^1(\bR^d)$ be a density function supported on $[0,1]^d$ with an absolute continuous gradient. Denote its Fisher information by
	\[
	J(f) \triangleq \int_{\bR^d} \frac{|\nabla f|^2}{f}.
	\]
	Then for any $p>1$, there exists a constant $C_p>0$, such that
	\[
	J(f) \leq C_p \sum_{i=1}^d \|\partial_{ii} f\|_p.
	\]
\end{theorem}
The connection between this result and the previous proof of \prettyref{lemma.first_approx_err} is the well-known fact that the local expansion of $\chi^2$-divergence is given by the Fisher information. Indeed, by Taylor expansion (assuming for simplicity that $d=1$ and $f=g$), the LHS of the main estimate \prettyref{eq.target} behaves as $h^2 J(f)$. Thanks to \prettyref{thm:fisher}, we can control the Fisher information by $J(f)=O(\|f''\|_2)$, which, by the smoothness assumption, is $O(Lh^{s-2})$, and leads to the desired \prettyref{eq.target}.

\subsection{Second-stage approximation error and variance}
In this subsection we analyze the performance of our pointwise estimator $\hat{H}(x)$ in estimating $-f_h(x)\ln f_h(x)$ for any $x\in \bR^d$.  Recall that
\begin{align*}
f_h(x) = \int_{\mathbb{R}^d} K_h(x-y)f(y)dy
\end{align*}
is the smoothed density, and 
\begin{align*}
\hat{f}_h(x) = \frac{1}{n}\sum_{i=1}^n K_h(x-X_i)
\end{align*}
is our estimator of $f_h(x)$, and $K_h(t)\triangleq h^{-d}K(h^{-1}t)$. In addition to the unbiasedness of $\hat{f}_h(x)$ in estimating $f_h(x)$, it satisfies some more properties: the random variable $h^d\hat{f}_h(x)$ roughly follows a binomial distribution $\mathsf{B}(n,h^df_h(x))$. This property confines to the one-dimensional simple example that $h^df_h(x)$ can be viewed as the discrete probability in the bin containing $x$, and $h^d\hat{f}_h(x)$ is the empirical frequency. Specifically, we prove the following lemma:
\begin{lemma}\label{lemma.regime}
If the kernel $K$ is non-negative everywhere, then there exists a constant $c_1>0$ depending on $d$ and $\|K\|_\infty$ only such that: 
\begin{enumerate}
\item If $f_h(x)\le\frac{c_1\ln n}{2nh^d}$, we have
\begin{align*}
\bP\left(\hat{f}_h(x)< \frac{c_1\ln n}{nh^d}\right) \ge 1-n^{-5d};
\end{align*}
\item If $f_h(x)\ge\frac{c_1\ln n}{2nh^d}$, we have
\begin{align*}
\bP\left(\frac{f_h(x)}{2}\le \hat{f}_h(x) \le 2f_h(x)\right) \ge 1-n^{-5d}.
\end{align*}
\end{enumerate}
\end{lemma}

Note that the concentration property of $\hat{f}_h(x)$ in Lemma \ref{lemma.regime} behaves as if $n \hat{f}_h(x)$ is distributed binomially as $\mathsf{B}(n,h^df_h(x))$. It shows that, based on our threshold to split the smooth and non-smooth regimes in \prettyref{eq:hatHx}, the probability of making an error in classification is negligible.

Now we prove that our estimators $\hat{H}_1$ and $\hat{H}_2$ in \prettyref{eq:H1x}--\prettyref{eq:H2x}  perform well in the corresponding regimes. 
To bound the variance, we invoke the well-known Efron--Stein--Steele inequality: 
\begin{lemma}[\cite{steele86}]\label{lemma:efron-stein}
	Let $X_1,X_2,\cdots,X_n$ be independent random variables, and for $i=1,\cdots,n$, let $X_i'$ be an independent copy of $X_i$. Then for any $f$,
	\begin{align*}
	\var(f(X_1,\cdots,X_n)) \le \frac{1}{2}\sum_{i=1}^n \bE(f(X_1,\cdots,X_n) - f(X_1,\cdots,X_{i-1},X_i',X_{i+1},\cdots,X_n))^2.
	\end{align*}
\end{lemma}
To apply \prettyref{lemma:efron-stein}, we need to bound the difference between the original estimator and the perturbed one where one observation is substituted by a fresh copy.
Recall that $\hat{H}_1$ depends only on the second part of observations $X^{(2)}$, and $\hat{H}_2$ depends on the last two parts $X^{(2)}\cup X^{(3)}$. Hence, we may define $\hat{H}_1', \hat{H}_2'$ to be the perturbed estimators where exactly one observation chosen uniformly at random from $X^{(2)}\cup X^{(3)}$ is replaced by its independent copy. The following lemmas summarize the upper bounds of the bias and the second moment of the perturbations: 
\begin{lemma}[Non-smooth regime]\label{lemma.non-smooth}
If $f_h(x)\le \frac{2c_1\ln n}{nh^d}$, $c_1\ge 2\|K\|_\infty c_2$, $0<7c_2\ln 2<\varepsilon$ and $n\ge 4c_2\ln n$, we have
\begin{align*}
|\bE\hat{H}_1(x) + f_h(x)\ln f_h(x)| &\lesssim \frac{1}{nh^d\ln n}, \\
\bE[(\hat{H}_1(x) - \hat{H}_1'(x))^2]&\lesssim \frac{1}{n^{3-\varepsilon}h^{2d}}.
\end{align*}
\end{lemma}

\begin{lemma}[Smooth regime]\label{lemma.smooth}
If $f_h(x)\ge \frac{c_1\ln n}{2nh^d}$ with sufficiently large constant $c_1>0$ (as in Lemma \ref{lemma.regime}), $h\le 1$ and $nh^d\ge 1$, we have
\begin{align*}
|\bE\hat{H}_2(x) + f_h(x)\ln f_h(x)| &\lesssim \frac{1}{nh^d\ln n}, \\
\bE[(\hat{H}_2(x) - \hat{H}_2'(x))^2] &\lesssim \frac{f_h(x)(1+(\ln f_h(x))^2)}{n^2h^d}. 
\end{align*}
\end{lemma}

For the final pointwise estimator $\hat{H}(x)$ in \prettyref{eq:hatHx}, let $\hat{H}'(x)$ be its perturbed version where exactly one observation chosen uniformly at random from $ X^{(2)}\cup X^{(3)}$ is replaced by its independent copy (note that $X^{(1)}$ is excluded). The following guarantee for $\hat{H}(x)$ follows from Lemma \ref{lemma.regime}--\ref{lemma.smooth}:
\begin{corollary}\label{cor.Ix}
Under the assumptions of Lemma \ref{lemma.regime}--\ref{lemma.smooth}, we have
\begin{align}
&\bE \left( \bE[\hat{H}(x)|X^{(1)}] + f_h(x)\ln f_h(x) \right)^2 \lesssim \frac{1}{(nh^d\ln n)^2} , \label{eq:Ix1} \\
&\bE[(\hat{H}(x) - \hat{H}'(x))^2]  \lesssim \frac{1}{n^{3-\varepsilon}h^{2d}}+\frac{f_h(x)(1+(\ln f_h(x))^2)}{n^2h^d} \label{eq:Ix2}.
\end{align}
\end{corollary}

\subsection{Overall performance}\label{subsec.overall}
Now we are ready to analyze the overall performance of the integrated estimator $\hat{H}=\int_{\bR^d} \hat{H}(x)dx$. As argued in Section \ref{sec:MIchi}, we assume without loss of generality that $f_h(\cdot)$ is supported on $[0,1]^d$, so that $\hat{H}=\int_{[0,1]^d} \hat{H}(x)dx$. By the triangle inequality, we have the following decomposition of the mean squared error (recall that $X^{(1)}$ is the first part of observations for regime classification): 
\begin{equation}\label{eq.decomp}
\begin{aligned}
\bE(\hat{H} - H(f))^2 &\le 2\left[ (H(f_h) - H(f))^2 + \bE(\hat{H} - H(f_h))^2 \right] \\
&= 2\left[ (H(f_h) - H(f))^2 + \bE(\bE[\hat{H}|X^{(1)}] - H(f_h))^2 + \bE[\var(\hat{H}|X^{(1)})] \right] . 
\end{aligned}
\end{equation}
We analyze different types of errors in \eqref{eq.decomp} separately: 
\begin{itemize}
\item First-stage approximation error: by Lemma \ref{lemma.first_approx_err} we know that
\begin{align*}
|H(f)-H(f_h)| \lesssim Lh^s.
\end{align*}
\item Conditional bias (second-stage approximation error): by Corollary \ref{cor.Ix} and Cauchy--Schwarz, 
\begin{align*}
\bE(\bE[\hat{H}|X^{(1)}] - H(f_h))^2 &= \bE\left( \int_{[0,1]^d} \left(\bE[\hat{H}(x)|X^{(1)}] + f_h(x)\ln f_h(x)\right) dx \right)^2  \\
&\le \int_{[0,1]^d} \bE\left( \bE[\hat{H}(x)|X^{(1)}] + f_h(x)\ln f_h(x)\right)^2 dx \\
&\lesssim \frac{1}{(nh^d\ln n)^2}. 
\end{align*}

\item Conditional variance: conditioned on $X^{(1)}$, the estimator $\hat{H}$ is a deterministic function of $(X^{(2)}, X^{(3)})$ consisting of mutually independent observations. We now apply \prettyref{lemma:efron-stein} to bound the variance. For $i=1,2,\cdots,2n$, 
define $\hat{H}_i(x)$ to be the pointwise estimator in \eqref{eq:hatHx} with $i$-th observation in $(X^{(2)}, X^{(3)})$ replaced by an independent copy, and let $\hat{H}_i = \int_{[0,1]^d} \hat{H}_i(x)dx$. Then by Lemma \ref{lemma:efron-stein}, we have
\begin{align*}
\var(\hat{H}|X^{(1)}) &\le \frac{1}{2}\sum_{i=1}^{2n} \bE[(\hat{H} - \hat{H}_i)^2 | X^{(1)} ] \\
&= \frac{1}{2}\sum_{i=1}^{2n} \bE\left[\left(\int_{[0,1]^d} (\hat{H}(x) - \hat{H}_i(x)) dx \right)^2 \bigg| X^{(1)} \right]. 
\end{align*}
Since $K$ has compact support (cf.~Assumption \ref{ass:kernel}), by our estimator construction we have $\mathsf{Leb}(\{x\in [0,1]^d: \hat{H}(x) \neq \hat{H}_i(x) \}) \lesssim h^d$. Hence, by Cauchy--Schwarz, we have
\begin{align*}
&\left(\int_{[0,1]^d} (\hat{H}(x) - \hat{H}_i(x)) dx \right)^2 \\
&\le\int_{[0,1]^d} (\hat{H}(x) - \hat{H}_i(x))^2 dx \cdot  \int_{[0,1]^d} \mathbbm{1}(\hat{H}(x)\neq\hat{H}_i(x))dx \\
&\lesssim h^d\int_{[0,1]^d} (\hat{H}(x) - \hat{H}_i(x))^2 dx.
\end{align*}
Combining the previous two displays, the conditional variance can be upper bounded as (recall the definition of $\hat{H}'(x)$ before Corollary \ref{cor.Ix})
\begin{align*}
\bE[\var(\hat{H}|X^{(1)})] &\lesssim nh^d \int_{[0,1]^d} \bE(\hat{H}(x) - \hat{H}'(x))^2dx \\
&\overset{\prettyref{eq:Ix2}}{\lesssim} \int_{[0,1]^d} \left( \frac{1}{n^{2-\varepsilon}h^d} + \frac{f_h(x)(1+(\ln f_h(x))^2)}{n} \right) dx \\
&\lesssim \frac{1}{n^{2-\varepsilon}h^d} + \frac{(\ln L)^2}{n}, 
\end{align*}
where the last inequality follows from Lemma \ref{lemma:varentropy} and $\|f_h\|_p\le \|f\|_p + \|f-f_h\|_p \lesssim L(1+h^s)$ (cf.~(\ref{eq:ffh})).

\end{itemize}

Substituting all three types of error bounds into \prettyref{eq.decomp}, we obtain
\begin{align}
\left(\sup_{f\in \Lip_{s,p,d}(L)}\bE_f\left(\hat{H}-H(f)\right)^2\right)^{\frac{1}{2}} 
\lesssim~& Lh^s + \frac{1}{nh^d\ln n} + \frac{1}{n^{1-\varepsilon/2}\sqrt{h^d}} + \frac{\ln L}{\sqrt{n}}. \label{eq:overall}
\end{align}
Finally, we choose 
$h = (Ln\ln n)^{-\frac{1}{s+d}}$ (note that the condition $Lh^s\le 1$ in Lemma \ref{lemma.first_approx_err} holds due to the assumption $L\le (n\ln n)^{s/d}$) and $\varepsilon<\frac{s}{s+d}$ in \eqref{eq:overall}
to obtain
\begin{align*}
\left(\sup_{f\in \Lip_{s,p,d}(L)}\bE_f\left(\hat{H}-H(f)\right)^2\right)^{\frac{1}{2}}  \lesssim (n\ln n)^{-\frac{s}{s+d}}L^{\frac{d}{s+d}} + n^{-\frac{1}{2}}\ln L,
\end{align*}
completing the proof of Theorem \ref{thm.achievability}.

\section{Further discussions}\label{sec.discussion}

\subsection{Extensions to densities satisfying periodic boundary conditions}\label{sec.extension}
In this section, we relax the assumptions that the underlying density $f$ is supported on $[0,1]^d$ and smoothly vanishing at the boundary, and establish the corresponding minimax rates in entropy estimation. 

Note that since the $L_p$ norm in \eqref{eq:lipnorm} is taken in $\bR^d$, the definition of the Lipschitz norm  requires that the density $f$ connects to zero smoothly at the boundary of $[0,1]^d$, which may exclude some well-known densities such as the uniform distribution. This assumption can be relaxed by considering the ``periodic boundary condition", which requires that the periodic extension of the density $f$ lies in the Lipschitz ball. Specifically, we define a new Lipschitz ball
\begin{align*}
\Lipstar_{s,p,d}(L) = \{f:  \|f\|_{\Lipstar_{s,p,d}} \le L \}\cap \{f: \supp(f)\subseteq [0,1]^d\}, 
\end{align*}
where the Lipschitz norm $\|\cdot\|_{\Lipstar_{s,p,d}}$ is defined in the same way as \eqref{eq:lipnorm} to \eqref{eq:finitediff}, with the only exceptions that the $L_p$ norm is taken over the unit cube $[0,1]^d$ instead of $\bR^d$, and $f$ is periodically extended to the entire space $\reals^d$ via $f(x \mod 1)\triangleq f(x_1-\lfloor x_1\rfloor, x_2-\lfloor x_2\rfloor, \cdots, x_d-\lfloor x_d\rfloor)$, for $x=(x_1,\cdots,x_d)\in \mathbb{R}^d$. Note that in this case we may also identify the unit cube $[0,1]^d$ as the $d$-dimensional torus $\torus^d$. 

The periodic boundary condition is weaker than the previous Lipschitz ball condition, in the sense of the norm comparison $\|f\|_{\Lipstar_{s,p,d}} \le C\|f\|_{\Lip_{s,p,d}}$ for some constant $C>0$.\footnote{This can be shown by applying the triangle inequality to $\|\sum_{i=1}^{3^d} f_i\|_{\Lip_{s,p,d}}$, where $f_i(x) = f(x-x_i)$ is the translation of $f$ with $\{x_1,\cdots,x_{3^d}\} = \{-1,0,1\}^d$. Consequently, one can choose $C=3^d$.} This assumption has already appeared in the literature~\cite{krishnamurthy2014nonparametric}, and the special case $s=2, d=1$ corresponds to $f(0)=f(1), f'(0)=f'(1)$. The next theorem shows that, the minimax rate remains unchanged in this weaker setting. 
\begin{theorem}\label{thm:periodic_boundary}
	For any $d\in \mathbb{N}, 0<s\le 2$ and $2\le p<\infty$, 
there exists $L_0>1$ depending on $s,p,d$, such that 
for any $L_0 \le L \le (n\ln n)^{s/d}$ and any $n\in\naturals$,
\begin{align*}
c( (n\ln n)^{-\frac{s}{s+d}}L^{\frac{d}{s+d}} + n^{-\frac{1}{2}}\ln L) 
&\leq \left(\inf_{\hat{H}}\sup_{f\in \Lipstar_{s,p,d}(L)} \bE_f (\hat{H}-H(f))^2 \right)^{\frac{1}{2}} \\
&\leq C((n\ln n)^{-\frac{s}{s+d}}L^{\frac{d}{s+d}} + n^{-\frac{1}{2}}\ln L)
\end{align*}
where $c, C>0$ are constants depending on $s,p,d$.
\end{theorem}

Theorem \ref{thm:periodic_boundary} is a straightforward extension of Theorem \ref{thm:main}. Since $\|\cdot\|_{\Lipstar_{s,p,d}}$ is a weaker norm than $\|\cdot\|_{\Lip_{s,p,d}}$, the minimax lower bound in Theorem \ref{thm:main} continues to hold for the new Lipschitz ball. As for the upper bound, we use the same estimator construction as in Section \ref{sec.construction}, with the understanding that the kernel convolution is taken with respect to the periodic extension of $f$ (or equivalently, is taken on the torus $\torus^d$). For  the analysis of this estimator, we apply a version of the maximal inequality on the torus $\torus^d$ in Section \ref{sec:firststage}, and the remaining arguments in Section \ref{sec.analysis} are essentially the same. We postpone the detailed proof to Section \ref{subsec.periodic_boundary_proof} in the appendix.

\subsection{The case of $s>2$}
\label{sec:highsmooth}

Note that the minimax lower bound in Theorem~\ref{thm.lower_bound} holds for all smoothness parameters $s>0$, but 
our current proof techniques of upper bound only work for the smoothness regime of $0<s\leq 2$. There are two main reasons:
\begin{enumerate}
\item Classifying smooth/nonsmooth regime (Lemma~\ref{lemma.regime}) fails when $s>2$;
\item Bias correction based on Taylor expansion (Lemma~\ref{lemma.smooth}) in the smooth regime does not extend to $s>2$. 
\end{enumerate}

The failures of Lemma~\ref{lemma.regime} and~\ref{lemma.smooth} are intrinsically related to the fact that one has to use kernels with negative parts to take advantage of smoothness $s>2$~\cite{Tsybakov2008}. Concretely, Lemma~\ref{lemma.regime} is closely related to the problem of adapting to the lowest values of density in density estimation~\cite{patschkowski2016adaptation}. It was conjectured in~\cite{patschkowski2016adaptation} that the case of $s>2$ exhibit significantly different behavior from the case of $0<s\leq 2$. Regarding bias correction, when the kernel is no longer non-negative, (\ref{eq:varKh}) can fail even when $f_h(x)>0$, which makes the proof of Lemma~\ref{lemma.smooth} break down. 
It is possible that bias correction based on Jackknife may achieve better performances when $s>2$. For the application of this approach in entropy estimation, we refer to~\cite{moon2016nonparametric, delattre2017kozachenko}.

Finally, we remark that the high smoothness regime of $s>2$ may not pose significant challenge for other problems of nonparametric statistics. For example, in the Gaussian white noise model, since it is a location model, the concentration of kernel estimators can be directly guaranteed using concentration inequality for sums of independent bounded random variables, which turn out to be sufficient for non-smooth functional estimation~\cite{Han--Jiao--Mukherjee--Weissman2017adaptive}. Even in the density model, 
 the case of $s>2$, which indeed calls for kernels with negative parts, can be easily handled. For example, to estimate density itself under $L_2$ risk, 
we can simply truncate the negative density estimates to obtain a better performance~\cite{Tsybakov2008}; in smooth functional estimation~\cite{birge1995estimation,tchetgen2008minimax}, the case of $s>2$ is also not special. It is mainly in estimating non-smooth functionals that designing procedures that can \emph{adapt} to low density regime becomes a crucial challenge~\cite{patschkowski2016adaptation}. 


\subsection{Connections to discrete entropy estimation}\label{sec:discreteentropyconnection}

Another intuitive idea for estimating the entropy of densities is to reduce it to a discrete entropy estimation problem. 
The motivation is that for a continuous random variable $X$ with density $f$, it is well-known \cite{Renyi59} that the Shannon entropy of its quantized version $[X]_k \triangleq \floor{kX}/k$ satisfies $H([X]_k)=d\log k + H(f)+o(1)$ as the quantization level $k\diverge$. Thus, to estimate $H(f)$, we can choose an appropriate $k$, quantize all the observations, and apply the optimal Shannon entropy estimator developed in \cite{Jiao--Venkat--Han--Weissman2015minimax,wu2016minimax}. Below we show that this approach achieves the minimax rate if $s \leq 1$.

For the ease of exposition, we consider $ s\in(0,1]$ and $p\ge 2$, with general dimension $d\in \mathbb{N}$. 
We split the unit cube $[0,1]^d$ into $S=h^{-d}$ sub-cubes $I_1,\ldots,I_S$ of size $h$, where $h$ is the ``bandwidth" we will choose later. For $i=1,\cdots,S$, 
define
\begin{align*}
p_i = \int_{I_i}f(t)dt, \quad \hat{p}_i = \frac{1}{n} \sum_{j = 1}^n \mathbbm{1}(X_j \in I_i). 
\end{align*}
as the ``probability" and the ``empirical frequency'' of the cube $I_i$, respectively. Since the entropy of the piecewise constant density
$
f_h(x) = \sum_{i=1}^S p_i\mathbbm{1}(x\in I_i)
$
is 
$
H(f_h) = \sum_{i  =1}^S p_i \ln \frac{1}{p_i} + \ln h^d, 
$
the problem of estimating $H(f_h)$ is reduced to a discrete Shannon entropy estimation problem. We can then use the minimax rate-optimal estimators $\hat{H}_{\mathsf{discrete}}$ \cite{wu2016minimax,Jiao--Venkat--Han--Weissman2015minimax} for the discrete entropy $\sum_{i=1}^S -p_i\ln p_i$
to define the estimator $\hat{H}$ for the entropy of the density $H(f)$ as
\begin{align}\label{eqn.discretereductionestimator}
\hat{H} =\hat{H}_{\mathsf{discrete}}+ \ln h^d.
\end{align}
One can show that $|H(f)-H(f_h)| =O(Lh^s)$. The optimal bandwidth is $h\asymp(Ln\ln n)^{-1/(s+d)}$, leading to the following risk bound, with an additional mild assumption that the density is bounded. 
\begin{lemma}\label{lemma.performancereducediscrete}
For $d \in \mathbb{N}, s\in (0,1], p\geq 2$, the performance of the estimator $\hat{H}$ in~(\ref{eqn.discretereductionestimator}) is given by
\begin{align*}
\left(\sup_{f\in\Lip_{s,p,d}(L), \|f\|_\infty \le L} \bE_f (\hat{H}-H(f))^2\right)^{\frac{1}{2}}\le C\left((n\ln n)^{-
\frac{s}{s+d}}L^{\frac{d}{s+d}}+n^{-\frac{1}{2}}\ln L\right), 
\end{align*}
where $C>0$ is a constant independent of $n, L$.
\end{lemma}


\appendix
\section{Auxiliary lemmas}
\label{app:aux}

The first lemma upper bounds the quantity $\int f(x)(\ln f(x))^2dx$ for any density $f\in \Lip_{s,p,d}(L)$. 
\begin{lemma}\label{lemma:varentropy}
Let $f$ be a density on $[0,1]^d$ with $\|f\|_p\le L$ for some $p>1$ and $L>0$. Then 
\begin{align*}
 \int_{[0,1]^d} f(x) (\ln f(x))^2 dx \le \frac{4}{e^2} + \frac{1+e^{1/(p-1)}}{(p-1)^2} + \left(\frac{p}{p-1}\right)^2(\ln L)^2.
\end{align*}
\end{lemma}

The following lemma presents the property of the best approximating polynomial of $-x\ln x$ on the interval $[0,\Delta]$, which follows from \cite[Section 7.5.4]{timan63}.
\begin{lemma}[Best approximating polynomial~\cite{Jiao--Venkat--Han--Weissman2015minimax,wu2016minimax}]\label{lemma.poly_approx}
Let $\Delta\le 1$. Denote by $Q(t)=\sum_{l=0}^k a_lt^l$ the best approximating polynomial of $-t\ln t$ on $[0,\Delta]$ in the uniform norm, we have
\begin{align*}
\sup_{t\in [0,\Delta]} |Q(t) + t\ln t| \leq \frac{C\Delta}{k^2}.
\end{align*}
where $C$ is a universal constant. Moreover, the coefficients satisfy:
\begin{align*}
|a_l| &\le 2^{3k}\Delta^{1-l}, \qquad l=0,2,3,\ldots,k\\
|a_1| &\le 2^{3k} - \ln \Delta.
\end{align*}
\end{lemma}

\begin{lemma}[Variance of second-order $U$-statistics]\label{lemma.U-statistic}
Let
\begin{align*}
U_2 \triangleq \frac{2}{n(n-1)} \sum_{1\le i<j\le n} X_iX_j,
\end{align*}
where $X_1,\ldots,X_n$ are i.i.d random variables with finite second moment. Then
\begin{align*}
\var(U_2) \le \frac{4(n-2)}{n(n-1)}(\bE X_1^2)(\bE X_1)^2 + \frac{2(\bE X_1^2)^2}{n(n-1)}.
\end{align*}
\end{lemma}

\begin{lemma}[Bennett's inequality \cite{Boucheron--Lugosi--Massart2013}]\label{lemma.bennett}
Let $X_1,\ldots,X_n\in [a,b]$ be independent random variables with 
\begin{align*}
\sigma^2 \triangleq \frac{1}{n}\sum_{i=1}^n \var(X_i).
\end{align*}
Then we have
\begin{align*}
\bP\left(\left|\frac{1}{n}\sum_{i=1}^n X_i - \frac{1}{n}\sum_{i=1}^n \bE[X_i]\right|\ge \varepsilon\right)\le 2\exp\left(-\frac{n\varepsilon^2}{2(\sigma^2+(b-a)\varepsilon/3)}\right).
\end{align*}
\end{lemma}

Finally, we need
the equivalence between Peetre's $K$-functional and modulus of smoothness on $\reals^d$. For a multi-index $\beta=(\beta_1,\ldots,\beta_d)$, we define $|\beta|\triangleq \sum_{i=1}^d|\beta_i|$, and write the differential operator $\prod_{i=1}^d(\partial/\partial x_i)^{\beta_i}$ as $D^\beta(\cdot)$. Now for $r\in\mathbb{N}$, the differential operator $D^r$ is defined as
\begin{align*}
D^rf \triangleq \sup_{|\beta|=r} |D^{\beta}f|.
\end{align*}
For $p\in [1,\infty]$ and $r\in\mathbb{N}$, the Peetre's $K$-functional for $f$ defined on $\bR^d$ is defined as
\begin{align*}
K_r(f,t^r)_p \triangleq \inf_{g} \|f-g\|_p + t^r\|D^rg\|_p
\end{align*}
where the infimum is taken over all functions $g$ defined on $\bR^d$ such that $D^{\beta}g$ for any $|\beta|=r-1$ is locally absolutely continuous. Also recall the definition of the modulus of smoothness $\omega_r(f,t)_p$ in \eqref{eq.moduli_of_smoothness}.
\begin{lemma}
\label{lmm:Kequiv}	
	For any $p\in [1,\infty]$ and $ r\in\mathbb{N}$, there exist universal constants $M=M(r,p)$ and $t_0=t_0(r,p)$ such that for any $0<t<t_0$ and $f$ defined on $\bR^d$,
	\begin{align}
	M^{-1}K_r(f,t^r)_p \le \omega_r(f,t)_p \le MK_r(f,t^r)_p.
	\label{eq:Kequiv}
	\end{align}
	Furthermore,
	\begin{equation}
	\omega_r(f,t)_p  \le M t^r \|D^r f\|_p, \quad 0<t<t_0
	\label{eq:Kequiv1}
	\end{equation}
	and
	\begin{equation}
	\omega_r(f,t)_p  \le 2^r \|f\|_p, \quad t>0.
	\label{eq:Kequiv2}
	\end{equation}	
\end{lemma}
For a proof this lemma, see \cite[Chapter 6, Theorem 2.4]{Devore--Lorentz1993} for the case of $d=1$, and \cite{johnen1977equivalence} for the general $\mathbb{R}^d$.
Finally, \prettyref{eq:Kequiv1} follows from \prettyref{eq:Kequiv} by choosing $g=f$, and 
\prettyref{eq:Kequiv2} follows from the definitions 
\prettyref{eq.moduli_of_smoothness}--\prettyref{eq:finitediff} and the triangle inequality.

\section{Proof of lower bound}\label{sec.lowerbound}
In this section we prove the following lower bound for entropy estimation: 
\begin{theorem}\label{thm.lower_bound}
	Let $s>0$ and $p\geq 1$. 
	Then there exist constants $L_0>1$ and $c_0>0$ depending on $s,p,d$ such that the following holds.
	For $L_0 \le L \le (n\ln n)^{s/d}$, let 
	\begin{align*}
	R_n \triangleq \left(\inf_{\hat{H}}\sup_{f\in \text{\rm Lip}_{s,p,d}(L)} \bE_f\left(\hat{H}-H(f)\right)^2\right)^{\frac{1}{2}}. 
	\end{align*}
	Then
	\begin{itemize}
		\item For any $p\in [1,\infty)$, $R_n \geq c_0 ((n\ln n)^{-\frac{s}{s+d}}L^{\frac{d}{s+d}} + n^{-\frac{1}{2}}\ln L)$.
		\item For $p=\infty$, $R_n \geq c_0 (n^{-\frac{s}{s+d}}(\ln n)^{-\frac{s+2d}{s+d}}L^{\frac{d}{s+d}}  + n^{-\frac{1}{2}}\ln L)$. 
	\end{itemize}
	
\end{theorem}

Note that Theorem \ref{thm.lower_bound} implies the lower bound part of Theorem \ref{thm:main}, which in fact holds for any $s>0$ and any $1 \leq p < \infty$. The second term $\Omega(n^{-1/2}\ln L)$ is essentially the parametric rate, where the only non-trivial part is to establish the proportionality to $\ln L$. This is shown below by means of a two-point argument provided that $L$ is at most polynomial in $n$:
\begin{lemma}\label{lemma:varentropy_lower}
Let $s>0$, $p\ge 1$. For any $\alpha > 0$, there exist constants $L_0,c_0,n_0>0$ depending on $s,p,\alpha$ such that for any $n\geq n_0$ and $L_0 \leq L \le n^{\alpha}$, 
\[
R_n \geq c_0 n^{-\frac{1}{2}}\ln L.
\]
\end{lemma}


Now our remaining goal is to show that for $p\in [1,\infty)$, 
\begin{equation}
R_n \gtrsim (n\ln n)^{-\frac{s}{s+d}}L^{\frac{d}{s+d}}.
\label{eq:lb-goal}
\end{equation}

The outline for the proof of \eqref{eq:lb-goal} is as follows:
In \prettyref{sec:reduction} we reduce the nonparametric problem to a parametric subproblem in the Poissonized sampling model. The minimax risk in the Poissonized sampling model is essentially no larger than\footnote{In fact, the minimax risks of the i.i.d.~sampling and Poissonized sampling models are closely related and it is simpler to show that the estimator constructed in \prettyref{sec.construction} satisfies the same risk bound in the Poisson model; see a previous version of the current paper \cite{HJWW17v2}.} that in the original i.i.d.~sampling model (cf.~Lemma \ref{lem.risk_multi_poisson}), and the central advantage of the Poissonized sampling model is that the sufficient statistics become independent (cf.~\eqref{eq:ss}), which simplifies the lower bound construction. The minimax lower bound for the parametric submodel is proved by a generalized version of Le Cam's method involving a pair of priors, also known as the method of two fuzzy hypotheses \cite{Tsybakov2008}. In \prettyref{subsec.measure} we construct the priors using duality to best approximation and \prettyref{sec:lb-complete} finishes the proof.

\subsection{Reduction to a parametric submodel}
\label{sec:reduction}

Let $d_0, d_1,d_2,d_3$ be positive constants to be specified later.
Set 
\begin{align}
\label{eq:Sh}
h= (d_0Ln\ln n)^{-\frac{1}{s+d}}, \quad S= (2h)^{-d}.
\end{align}
Without loss of generality we assume that $S$ is an integer. 
Fix  
\begin{equation}
\alpha \in \pth{0, \frac{1}{S}}.
\label{eq:alpha}
\end{equation}
Let $I_1,\ldots,I_S$ be the partition of the large cube $[\frac{1}{4},\frac{3}{4}]^d$ into $S$ smaller cubes of edge length $h$, and $t_i$ denote the leftmost corner of sub-cube $I_i$, i.e., $I_i=t_i+[0,h]^d$. 
Let $g$ and $w$ be some fixed smooth probability density functions on $\reals^d$ with finite entropy which vanish outside $[0,1]^d$ and $[0,1]^d\backslash[\frac{1}{4},\frac{3}{4}]^d$, respectively. The smoothness of $g$ and $w$ on $\mathbb{R}^d$ implies that their derivatives of any order vanish at their respective boundary.

To each vector $P=(p_1,\ldots,p_S)\in [0,\frac{d_3\ln n}{n}]^S$, we associate a function
\begin{align}\label{eq.f_P}
f_P(t) \triangleq \sum_{i=1}^S p_i\cdot \frac{1}{h^d}g\pth{\frac{t-t_i}{h}} + \left(1-S \alpha\right)w(t).
\end{align}
By the choice of $\alpha$, $f_P$ is non-negative everywhere. Moreover, since $w$ is smooth on $\mathbb{R}^d$ and vanishes outside $[0,1]^d$, the function $f_P$ connects to zero smoothly on the boundary.
The next result is a sufficient condition for $f_P$ to lie in the Lipschitz ball.
\begin{lemma}
	\label{lmm:fP}	
	If $L\ge L_0$ for some constant $L_0$ depending only on $w$, 
	\begin{align}
	\frac{1}{S}\sum_{i=1}^Sp_i^p \le \pth{\frac{2C_1}{n\ln n}}^p,
	\label{eq:PP}
	\end{align}
	and $d_0>0$ is sufficiently small depending only on $(g,L_0,C_1,s,p,d)$, then $f_P \in \Lip_{s,p,d}(L)$. 
\end{lemma}
\begin{proof}
	Let
	\begin{align*}
	\tilde{f}_P(t) \triangleq \sum_{i=1}^S p_i\cdot \frac{1}{h^d}g\pth{\frac{t-t_i}{h}}.
	\end{align*}
	Then $f_P = \tilde f_P + (1-S\alpha) w$.
	Since $\|f_P\|_{\Lip} \leq \|\tilde f_P\|_{\Lip} + (1-S\alpha) \|w\|_{\Lip}$ and 
	$S\alpha \leq 1$, after choosing $L_0 = 2\|w\|_{\Lip}$, it suffices to show $\tilde{f}_P\in \text{Lip}_{s,p,d}(L/2)$ for $L\ge L_0$. 
	Observe that
	\begin{align}
	\|\tilde{f}_P\|_p = \frac{\|g\|_p}{h^d}\left(\frac{1}{S}\sum_{i=1}^Sp_i^p\right)^{\frac{1}{p}} = h^s\|g\|_p\cdot d_0Ln\ln n\left(\frac{1}{S}\sum_{i=1}^Sp_i^p\right)^{\frac{1}{p}},
	\label{eq:tfP}
	\end{align}
	the condition \cref{eq:PP} ensures that $h^{-s}\|\tilde{f}_P\|_p$ is upper bounded by $2d_0C_1\|g\|_p\cdot L$. 
	By \cref{eq:Kequiv1} in \Cref{lmm:Kequiv}, \prettyref{eq:tfP} implies that $\omega_r(\tilde{f}_P,t)\le Lt^s/4$ for any $t\ge h$ and $d_0$ small.
	For $t<h$, further observe that for $r=\lceil s\rceil$, we have
	\begin{align}
	\|D^r\tilde{f}_P\|_p &= \frac{\|D^r g\|_p}{h^{d+r}}\left(\frac{1}{S}\sum_{i=1}^Sp_i^p\right)^{\frac{1}{p}} = \frac{\|D^r g\|_p}{h^{r-s}}\cdot d_0Ln\ln n\left(\frac{1}{S}\sum_{i=1}^Sp_i^p\right)^{\frac{1}{p}}
	\label{eq:rfP}
	\end{align}
	Hence, $h^{r-s}\|D^r\tilde{f}_P\|_p$ is upper bounded by $2d_0C_1\|D^r g\|_p\cdot L$.
	By \cref{eq:Kequiv2} in \Cref{lmm:Kequiv}, \prettyref{eq:rfP} implies that $\omega_r(\tilde{f}_P,t)\le Lt^s/4$ for any $t\le h$ and $d_0$ sufficiently small, completing the proof.
\end{proof}

Now we would like to reduce the original nonparametric model $f\in \text{Lip}_{s,p,d}(L)$ to some parametric submodel $f\in \{f_P:P\in \mathcal{P}\}$, with a properly chosen $\mathcal{P}$, and impose an iid prior on the coefficient vector $P$ with mean $\alpha$. However, $f_P$ need not be a valid density, since
\begin{align*}
\int_{[0,1]^d} f_P(t)dt = 1 + \sum_{i=1}^S(p_i-\alpha)
\end{align*}
need not normalize to one. 
To resolve this issue, we show that the minimax rate 
remains unchanged as long as $\|f\|_1$ is sufficiently close to one. To this end, let us
define a new model: instead of the original sampling model with a fixed sample size $n$, we draw $N\sim \mathsf{Poi}(n\|f\|_1)$ i.i.d. observations from the density $\frac{f}{\|f\|_1}$. In the new model, the parameter set of $f$ is a parametric one $\{f_P: P\in \mathcal{P}\}$, where
\begin{equation}
\begin{aligned}\label{eq:calP}
\mathcal{P}&\triangleq  \left\{P \in \qth{0,\frac{d_3\ln n}{n}}^S: \left|\sum_{i=1}^S(p_i-\alpha)\right|\le \frac{1}{nh^d(\ln n)^3\ln L}\right\}\\
&\qquad\cap \left\{P: \frac{1}{S}\sum_{i=1}^S |p_i|^p \le \pth{\frac{2C_1}{n\ln n}}^p\right\} 
\end{aligned}
\end{equation}
with $C_1$ given in (\ref{eq.bound_mean}). In view of \prettyref{lmm:fP}, we have $\{f_P:P\in\mathcal{P}\} \subseteq \text{Lip}_{s,p,d}(L)$.
We also extend the definition of entropy to positive functions verbatim: $H(f_P)=\int_{[0,1]^d} -f_P(t)\ln f_P(t)dt$, and denote the minimax risk restricted to the parametric submodel $f\in \{f_P: P\in\mathcal{P}\}$ as
\begin{align*}
R_n^P \triangleq \left(\inf_{\hat{H}}\sup_{P\in\mathcal{P}} \bE_{f_P}(\hat{H}-H(f_P))^2\right)^{\frac{1}{2}}.
\end{align*}

The following lemma relates the new minimax risk $R_n^P$ to the original risk $R_n$ in \prettyref{thm.lower_bound}. 
\begin{lemma}\label{lem.risk_multi_poisson}
	There exists a constant $C_2>0$ depending on $s,p$ only that
	\begin{align*}
	R_n^P \le C_2\left(R_\frac{n}{2} + \frac{1}{nh^d(\ln n)^3}\right).
	\end{align*}
\end{lemma}

Since our goal \prettyref{eq:lb-goal} is to show $R_n\gtrsim \frac{1}{nh^d\ln n}$, by Lemma \ref{lem.risk_multi_poisson}, it suffices to prove $R_n^P\gtrsim \frac{1}{nh^d\ln n}$. 
Note that
\begin{align*}
H(f_P) = C_0 + H(P) +(H(g)+d\ln h) \sum_{i=1}^Sp_i
\end{align*}
where the constant $C_0$ does not depend on $P$, and we denote, with a slight abuse of notation, the discrete entropy $H(P) \triangleq \sum_{i=1}^S -p_i\ln p_i$.
Hence, given the minimax optimal estimator $\hat{H}$ for $H(f_P)$, the following estimator
\begin{align*}
\tilde{H} \triangleq \hat{H}-C_0  - \left(H(g)+d\ln h\right)\cdot S\alpha
\end{align*}
for $H(P)$ satisfies
\begin{align*}
(\Expect[(\tilde{H}-H(P))])^{1/2} 
& \leq (\Expect[(\hat{H}-H(f_P))])^{1/2} + \left|H(g)+d\ln h\right|\cdot \left|\sum_{i=1}^S(p_i-\alpha)\right|\\
& \lesssim R_n^P + \frac{1}{\ln n}\cdot \frac{1}{nh^d\ln n}. 
\end{align*}
Therefore, to show \prettyref{eq:lb-goal}, it suffices to prove
\begin{align}\label{eq.para_lower}
\inf_{\hat{H}}\sup_{P\in\mathcal{P}} \bE_P\left( \hat{H} - H(P)\right)^2 \gtrsim \frac{1}{(nh^d\ln n)^2}.
\end{align}

Finally, we note that by the factorization theorem, for the Poisson sampling model, to estimate $H(P)$ with $P\in\mathcal{P}$, the histograms 
\begin{align*}
h_j = \sum_{i=1}^N \mathbbm{1}(X_i\in I_j), \qquad j=1,\ldots,S
\end{align*}
constitute a sufficient statistic. As a result, we may further assume that our observation model is
\begin{align}
\label{eq:ss}
h_j \inddistr~ \mathsf{Poi}(n\cdot p_j), \qquad j=1,\ldots,S
\end{align}
and the estimator $\hat H$ in \prettyref{eq.para_lower} is a function of $(h_1,\ldots,h_S)$.

\subsection{Construction of two priors}\label{subsec.measure}

The minimax lower bound \prettyref{eq.para_lower} follows from a generalized version of Le Cam's method involving two priors, also known as the method of two fuzzy hypotheses \cite{Tsybakov2008}. 
Given a collection of distributions $\{P_\theta: \theta\in\Theta'\}$, suppose the observation ${\bf Z}$ is distributed as $P_\theta$ with $\theta \in \Theta\subseteq \Theta'$. 
Let $\hat{T} = \hat{T}({\bf Z})$ be an arbitrary estimator of a function $T(\theta)$ based on $\bf Z$. 

Denote the total variation distance between two probability measures $P,Q$ by 
\begin{equation*}
V(P,Q) \triangleq \sup_{A\in \mathcal{A}} | P(A) - Q(A) | = \frac{1}{2} \int |p-q| d\nu,
\end{equation*}
where $p = \frac{dP}{d\nu}, q = \frac{dQ}{d\nu}$, and $\nu$ is a dominating measure so that $P \ll \nu, Q \ll \nu$.
The following general minimax lower bound follows from the same proof as \cite[Theorem 2.15]{Tsybakov2008}:
\begin{lemma}\label{lemma.tsybakov}
	Let $\sigma_0$ and $\sigma_1$ be prior distributions on $\Theta'$.
	Suppose there exist $\zeta\in \mathbb{R}, s>0, 0\leq \beta_0,\beta_1 <1$ such that
	\begin{align*}
	\sigma_0(\theta \in \Theta: T(\theta) \leq \zeta -s) & \geq 1-\beta_0, \\
	\sigma_1(\theta \in \Theta: T(\theta) \geq \zeta + s) & \geq 1-\beta_1.
	\end{align*}
	Then
	\begin{equation*}
	\inf_{\hat{T}} \sup_{\theta \in \Theta} \bP_\theta\left( |\hat{T} - T(\theta)| \geq s \right) \geq \frac{1-V(F_1,F_0)  - \beta_0 - \beta_1}{2},
	\end{equation*}
	where 
	$F_i=\int P_\theta \sigma_i(d \theta)$ is the marginal distribution of $\mathbf{Z}$ under the prior $\sigma_i$ for $i = 0,1$, respectively.
\end{lemma}



To apply \prettyref{lemma.tsybakov} to $\Theta'=[0,\frac{d_3\ln n}{n}]^S$ and $\Theta=\calP$, we first describe the construction of the two priors.
The following result is simply the duality between the problem of best uniform approximation and moment matching.\footnote{The proof of \prettyref{lem.measure} is identical to that of 
	\cite[Eqn.~(34)]{wu2016minimax}, since $\{x^{-q+1},x^{-q+2},\ldots,x^k\}$ forms a Haar system \cite[Section 3.3, Example 2]{Devore--Lorentz1993} and hence the Chebyshev alternating theorem holds.}


\begin{lemma}\label{lem.measure}
	Given a compact interval $I=[a,b]$ with $a>0$, integers $q,k>0$ and a continuous function $\phi$ on $I$, let 
	\begin{equation}
	E_{q-1,k}(\phi;I) \triangleq \inf_{\{a_i\}} \sup_{x\in I} \left| \sum_{i=-q+1}^k a_i x^i
	-\phi(x)\right|
	\label{eq:rationalapproxerror}
	\end{equation}
	denote the best uniform approximation error of $\phi$ by rational functions spanned by $\{x^{-q+1},x^{-q+2},\ldots,x^k\}$.
	Let $E_k(\phi;I) \triangleq E_{q-1,k}(\phi;I)$ denote the best uniform approximation error of $\phi$ by degree-$k$ polynomials.
		Then
	\begin{equation}
	\begin{aligned}
	2 E_{q-1,k}(\phi;I) = \max & ~ \int \phi(t) \nu_1(dt) - \int \phi(t) \nu_0(dt)   \\
	\text{\rm s.t.}     & ~ \int t^{l} \nu_1(dt) = \int t^{l} \nu_0(dt), \quad l=-q+1,\ldots,k
	\end{aligned}
	\label{eq:Rstar}
	\end{equation}
	where the maximum is taken over pairs of probability measures $\nu_0$ and $\nu_1$ supported on $I$.
\end{lemma}

Here we apply this lemma to $\phi_q(t)\triangleq t^{1-q}\ln t$ and 
\begin{align*}
\eta=\frac{d_1}{(\ln n)^2}, \quad I = \qth{\eta, 1}, \qquad k=\lceil d_2\ln n\rceil, \qquad q=\lceil p\rceil
\end{align*}
with constants $d_1,d_2>0$ to be specified later. The following lemma provides a lower bound for the approximation error $\phi_q(t) \triangleq t^{1-q}\ln t$:
\begin{lemma}\label{lem.approx}
	Fix $q\in \mathbb{N}$. There exists constants $c,c'>0$ depending on $q$ only such that
	\begin{align*}
	\liminf_{k\to\infty} k^{-2(q-1)}E_{q-1,k}\left(\phi_q;\qth{\frac{c}{k^2},1}\right) \ge c'.
	\end{align*}
\end{lemma}
Choosing $d_1=c/d_2^2$ with constant $c>0$ given in Lemma \ref{lem.approx} and in view of the definitions of $I$ and $k$, we conclude that
\begin{align*}
E_{q-1,k}(\phi_q;I) \gtrsim (\ln n)^{2q-2}.
\end{align*}
Let $\nu_0$ and $\nu_1$ be the maximizer of \prettyref{eq:Rstar}. We define probability measures $\tilde{\mu}_0, \tilde{\mu}_1$ on $I = \qth{\eta, 1}$ by
\begin{align*}
\tilde{\mu}_i(dt) = \pth{1- \int \pth{\frac{\eta}{t}}^{q} \nu_i(dt)  }\delta_0(dt) + \pth{\frac{\eta}{t}}^{q} \nu_i(dt)\qquad i=0,1.
\end{align*}
where $\delta_0$ is the delta measure at zero. It is straightforward to verify that $\tilde{\mu}_0$ and $\tilde{\mu}_1$ are probability measures, satisfying
\begin{enumerate}
	\item $\int t^l\tilde{\mu}_1(dt) = \int t^l\tilde{\mu}_0(dt)$, for all $l=0,1,\ldots,q+k$;
	\item $\int t\ln t\tilde{\mu}_1(dt) - \int t\ln t\tilde{\mu}_0(dt) \gtrsim (\ln n)^{-2}$;
	\item $\int t^{q}\tilde{\mu}_i(dt) = \eta^q= d_1^q(\ln n)^{-2q}$, for $i=0,1$.
\end{enumerate}

Finally, let $\mu_i$ be the dilation of $\tilde \mu_i$ by a factor of $\frac{d_3\ln n}{n}$, that is, if $Y_i \sim \tilde{\mu}_i$, then $\frac{d_3\ln n}{n} Y_i \sim \mu_i$. 
Then $\mathsf{supp}(\mu_i)\subseteq [0,\frac{d_3\ln n}{n}]$ and, furthermore,
\begin{enumerate}
	\item $\int t^l\mu_1(dt) = \int t^l {\mu}_0(dt)$, for all $l=0,1,\ldots,q+k$;
	\item $\int t\ln t{\mu}_1(dt) - \int t\ln t{\mu}_0(dt) \gtrsim (n\ln n)^{-1}$;
	\item $\int t^{q}{\mu}_i(dt) = (\frac{\eta d_3\ln n}{n})^q = (d_1d_3)^q (n\ln n)^{-q}$, for $i=0,1$.
\end{enumerate}

In the next subsection we will impose the prior where the parameters $(p_1,\ldots,p_S)$ are iid as either $\mu_0$ or $\mu_1$.
The utility of the above these properties are as follows:
The first condition ensures the priors $\mu_1$ and $\mu_2$ have matching first $q+k$ moments, and the induced Poisson mixtures are exponentially close in total variation; this is exactly the distribution of the sufficient statistics $(h_1,\ldots,h_S)$.
The second property ensures the separation of the functional values, while the third property ensures the smoothness of the functions $f_P$. Indeed, since $q\ge p$, H\"{o}lder's inequality yields
\begin{align}\label{eq.bound_mean}
\int t \mu_i(dt) \leq 
\left(\int t^p \mu_i(dt) \right)^{\frac{1}{p}} \le \left(\int t^{q} \mu_i(dt) \right)^{\frac{1}{q}}\le \frac{C_1}{n\ln n}
\end{align}
where $C_1$ is a constant depending on $d_2,d_3$ and $p$. It implies that controlling the $q$th moment of $\mu_i$ for $q\ge p$ also controls all of its lower-order moments.

\subsection{Minimax lower bound in the parametric submodel}
\label{sec:lb-complete}
In this subsection we invoke Lemma \ref{lemma.tsybakov} to finish the proof of (\ref{eq.para_lower}), thereby proving the lower bound in Theorem \ref{thm:main}. Consider the probability measures $\mu_0,\mu_1$ constructed in \Cref{subsec.measure}, and define
\begin{align}
\label{eq:Delta}
\Delta \triangleq \int t\ln t\mu_1(dt) - \int t\ln t\mu_0(dt) \gtrsim \frac{1}{n\ln n}. 
\end{align}
Put $\alpha = \int t \mu_0(dt) = \int t \mu_1(dt)$ and recall the parameter space $\calP$ defined in \prettyref{eq:calP}. By \prettyref{eq.bound_mean} and the assumption $L \le (n\ln n)^{s/d}$, for large $d_2$ (and therefore small $d_1$) we have
\begin{align*}
0\leq S \alpha \lesssim \frac{d_1S}{n\ln n}\lesssim \frac{d_1}{nh^d\ln n}\le 1, 
\end{align*}
which fulfills \prettyref{eq:alpha}.

Let $\mu_i^{\otimes S}$ denote the $S$-fold product of $\mu_i$. 
Consider the following event:
\begin{align*}
E_i \triangleq \{P: P\in\mathcal{P}\} \cap \left\{P: \left|\frac{1}{S}\sum_{j=1}^Sp_j\ln p_j-\bE_{\mu_i}p\ln p\right| \le \frac{\Delta}{4} \right\}, \qquad i=0,1. 
\end{align*}
We first show that $\mu_i^{\otimes S}(E_i)\to 0$ as $n\to\infty$ for any $i=0,1$. Recall the definitions of $S$ and $h$ from \prettyref{eq:Sh}.
By Chebyshev's inequality and the fact that $\mu_i$ is supported on $[0,\frac{d_3 \log n}{n}]$, we have
\begin{align*}
\mu_i^{\otimes S}\left\{ P:\left|\sum_{i=1}^S(p_i-\bE_{\mu_i}p_i)\right|> \frac{1}{nh^d(\ln n)^3\ln L} \right\}&\le (nh^d(\ln n)^3\ln L)^2\cdot \var_{\mu_i^{\otimes S}}\left[\sum_{i=1}^S p_i\right]\\
&\le (nh^d(\ln n)^3\ln L)^2\cdot S\pth{\frac{d_3\ln n}{n}}^2\\
&\lesssim h^d(\ln n)^8(\ln L)^2 \to 0.
\end{align*}
Similarly, 
\begin{align*}
&\mu_i^{\otimes S}\left\{P: \frac{1}{S}\sum_{i=1}^S |p_i|^p > \pth{\frac{2C_1}{n\ln n}}^p\right\} \\
& \overset{\prettyref{eq.bound_mean}}{\le} \mu_i^{\otimes S}\left\{P: \frac{1}{S}\sum_{i=1}^S (|p_i|^p - \bE_{\mu_i}|p_i|^p) > \pth{\frac{C_1}{n\ln n}}^p\right\}\\
&\le \pth{\frac{C_1}{n\ln n}}^{-2p}\cdot \var_{\mu_i^{\otimes S}}\left[\frac{1}{S}\sum_{i=1}^S |p_i|^p \right]\\
&\lesssim (n\ln n)^{2p}\cdot \frac{1}{S} \pth{\frac{\ln n}{n}}^{2p}\\
&\asymp h^d(\ln n)^{4p} \to 0
\end{align*}
and
\begin{align*}
&\mu_i^{\otimes S}\left\{P: \left|\frac{1}{S}\sum_{j=1}^Sp_j\ln p_j-\bE_{\mu_i}[p\ln p]\right| > \frac{\Delta}{4} \right\} \\
&\le \frac{16}{\Delta^2}\mathsf{Var}_{\mu_i^{\otimes S}}\left[\frac{1}{S}\sum_{j=1}^S p_j\ln p_j\right]\\
&\le \frac{16}{S\Delta^2}\cdot \pth{\frac{(\ln n)^2}{n}}^{2} \\
&\overset{\prettyref{eq:Delta}}{\lesssim} h^d(\ln n)^6\to 0.
\end{align*}
Hence, by the union bound, we have $\beta_i\triangleq \mu_i^{\otimes S}(E_i^c) \to 0$
for $i=0,1$.

Now we are ready to apply Lemma \ref{lemma.tsybakov} to
\begin{align*}
T(\theta) &= H(P) = \sum_{i=1}^S -p_i\ln p_i\\
\zeta &= \frac{\bE_{\mu_1^{\otimes S}} H(P) + \bE_{\mu_0^{\otimes S}} H(P)}{2}, \quad s = \frac{S\Delta}{4}
\end{align*}
Under the prior $\sigma_i \triangleq \mu_i^{\otimes S}$, the sufficient statistics $(h_1,\ldots,h_S)$ are i.i.d.~with distribution $F_i \triangleq \gamma_i^{\otimes S}$, where $\gamma_i \triangleq \int \Poi(n \lambda) \mu_i(d\lambda)$ is a Poisson mixture, in view of \prettyref{eq:ss}.
Note that $\mu_0$ and $\mu_1$ have matching first $q+k$ moments. 
By \cite[Lemma 3]{wu2016minimax}, we have
\[
V(\pi_0,\pi_1) \leq \pth{\frac{2e d_3 \ln n}{q+k} }^{q+k} \leq n^{-d_2\log \frac{d_2}{2e d_3}}
\]
provided that $d_2\log \frac{d_2}{2e d_3} \geq 2$. Therefore
$V(F_0,F_1) \leq S \cdot V(\pi_0,\pi_1)\to 0$ by choosing constant $d_2>0$ large enough.
Applying Lemma \ref{lemma.tsybakov} together with Markov's inequality yields
\begin{align*}
\inf_{\hat{H}}\sup_{P\in\mathcal{P}} \bE_P\left( \hat{H} - H(P)\right)^2 \gtrsim s^2\cdot \inf_{\hat{H}}\sup_{P\in\mathcal{P}} \bP\left(|\hat{H}-H(P)|\ge s\right) \gtrsim \frac{1}{(nh^d\ln n)^2}, 
\end{align*}
which is (\ref{eq.para_lower}), as desired.

Finally, the lower bound for $p=\infty$ follows from the same argument, except that we will set $q=1$, choose the bandwidth differently
\[
h=\left(\frac{\ln n}{d_0Ln}\right)^{\frac{1}{s+d}},
\]
and replace the RHS of \eqref{eq.bound_mean} by $\frac{d_3\ln n}{n}$.

\section{Proof of Theorems~\ref{thm:orlicz}, \ref{lem.plug-in} and \ref{thm:periodic_boundary}}

\subsection{Proof of Theorem \ref{thm:orlicz}}\label{subsec.subGaussian_proof}
We prove the upper and lower bounds separately. 
\subsubsection{Proof of upper bound}
First we propose the entropy estimator $\hat{H}$ for the unbounded support case. We begin with a lemma for the tail of density $f\in \Lip_{s,p,d}^{\Psi}(L)$. 

\begin{lemma}\label{lemma:orlicz_tail}
Let $f\in \Lip_{s,p,d}^{\Psi}(L)$ with $p\in [2,\infty)$ and Orlicz function $\Psi$ of rapid growth. There exist constants $C_0, C_1, C_2>0$ depending on $p,d,\kappa(\Psi), \Psi(1),L$ only, such that 
\begin{align}
\left| \int_{|x|\ge C_0\Psi^{-1}(n) } f(x)\ln \frac{1}{f(x)}dx\right| &\le \frac{C_1}{n^4}, \label{eq:orlicz_tail} \\
\int_{\reals^d} f(x)\ln^2 f(x)dx &\le C_2. \label{eq:orlicz_varentropy}
\end{align}
\end{lemma}
\begin{proof}
We first prove the second inequality \eqref{eq:orlicz_varentropy}. Since $t\ln^2 t\le t^{1/2} + t^{3/2}$ for all $t\ge 0$, 
\begin{align}\label{eq:decomposition}
\int_{\reals^d} f(x)\ln^2 f(x)dx \le \int_{\reals^d} f(x)^{\frac{3}{2}}dx + \int_{\reals^d} f(x)^{\frac{1}{2}}dx. 
\end{align}
By Cauchy--Schwartz and $\|f\|_2 \le \|f\|_p \le L$, we have
\begin{align}\label{eq.one_and_half}
\int_{\reals^d} f(x)^{\frac{3}{2}}dx \le \left( \int_{\reals^d} f(x)^2dx\right)^{\frac{1}{2}}\left(\int_{\reals^d} f(x)dx\right)^{\frac{1}{2}} \le L. 
\end{align}
For the second integral, by the convexity of $\Psi$ and $\Psi(0)=0$ we know that there is a constant $a>0$ depending only on $\Psi(1)$ such that $\Psi(a)\ge 2$. Then by Cauchy--Schwartz again, 
\begin{align}\label{eq.one_half_inner}
\int_{|x|\le a} f(x)^{\frac{1}{2}}dx \le (2a)^{\frac{d}{2}} \left(\int_{|x|\le a} f(x)dx\right)^{\frac{1}{2}} \le a^{\frac{d}{2}}. 
\end{align}
Since for any $x\in \reals^d$, we have
\begin{align*}
2f(x)^{\frac{1}{2}} \le \Psi(|x|)f(x)+ \frac{1}{\Psi(|x|)}. 
\end{align*}
Integrating over $|x|> a$ gives
\begin{align}\label{eq.one_half_outer}
2\int_{|x|> a} f(x)^{\frac{1}{2}}dx &\le \int_{|x|>a} \Psi(|x|)f(x)dx + \int_{|x|>a} \frac{1}{\Psi(|x|)}dx \nonumber \\
&\le L + \sum_{k=0}^\infty \int_{\kappa^{k}a < |x| \le \kappa^{k+1} a} \frac{1}{\Psi(|x|)}dx \nonumber \\
&\le L + \sum_{k=0}^\infty \frac{(2\kappa^{k+1} a)^d}{\Psi(\kappa^{k} a)} < \infty, 
\end{align}
where the finiteness follows from $\Psi(\kappa^{k}a ) \ge \Psi(a)^{2^k} \ge 2^{2^k}$ by the rapid growth assumption with parameter $\kappa=\kappa(\Psi)$. Combining \eqref{eq:decomposition}--\eqref{eq.one_half_outer} gives \eqref{eq:orlicz_varentropy}. 

To establish the first inequality \eqref{eq:orlicz_tail}, let $X\sim f$. Then by Cauchy--Schwartz, 
\begin{align}
\left| \int_{|x|\ge C_0\Psi^{-1}(n) } f(x)\ln \frac{1}{f(x)}dx\right| &= \left| \bE\left[\ln \frac{1}{f(X)} \mathbbm{1}(|X|\ge C_0\Psi^{-1}(n)) \right]\right| \nonumber \\
&\le \sqrt{\bE[(\ln f(X))^2]\cdot \bP(|X|\ge C_0\Psi^{-1}(n))}.  \label{eq.cauchy--schwartz}
\end{align}
The probability $\bP(|X|\ge C_0\Psi^{-1}(n))$ can be upper bounded easily: 
\begin{align*}
\bP(|X|\ge C_0\Psi^{-1}(n)) \le \frac{1}{\Psi(C_0\Psi^{-1}(n))} \int_{\reals^d} \Psi(|x|)f(x)dx \le \frac{L}{\Psi(C_0\Psi^{-1}(n))}. 
\end{align*}
By the rapid growth property of $\Psi$, choosing $C_0 = \kappa^3$ we have $\Psi(C_0\Psi^{-1}(n))\ge n^8$, and therefore
\begin{align}\label{eq.tail_probability}
\bP(|X|\ge C_0\Psi^{-1}(n)) \le \frac{L}{n^8}. 
\end{align}

Using the upper bound of $\bE[(\ln f(X))^2]$ in \eqref{eq:orlicz_varentropy}, a combination of \eqref{eq.cauchy--schwartz} and \eqref{eq.tail_probability} gives \eqref{eq:orlicz_tail} and completes the proof of the lemma. 
\end{proof}

Let $C_0$ be given in Lemma \ref{lemma:orlicz_tail}, and $R \triangleq C_0\Psi^{-1}(n)\ge 1$. We define the entropy estimator $\hat{H}$ by
\begin{align*}
\hat{H} = \int_{|x|\le R} \hat{H}(x)dx,
\end{align*}
where the pointwise estimator $\hat{H}(x)$ is defined in \eqref{eq:hatHx}, with the same parameters except that the bandwidth $h$ is chosen to be
\begin{align}\label{eq:orlicz_bandwidth}
h = c_0(n\ln n)^{-\frac{1}{s+d}}\cdot R^{\frac{d}{p(s+d)}}. 
\end{align}
To show the upper bound, first note that if the density $f$ is indeed supported on $[-R,R]^d$, then the scaled density $g(x) = f(R(2x-1))$ is supported on $[0,1]^d$, with $g \in \Lip_{s,p,d}(L(2R)^{s+d(1-1/p)})$ and $\|g\|_p \le L(2R)^{d(1-1/p)}$. Moreover, the rapid growth property of $\Psi$ ensures that $R= o(n^\varepsilon)$ for any $\varepsilon>0$, where the hidden constant depends on $\Psi(1), \kappa(\Psi)$ and $\varepsilon$. Following the same analysis in Section \ref{sec.analysis}, and noting that the scaled bandwidth becomes $h/(2R)$, we have
\begin{align*}
\left( \bE_f \left(\hat{H} - H(f)\right)^2 \right)^{\frac{1}{2}} \le C\left( R^{s+d(1-1/p)}\cdot (h/R)^s + \frac{1}{n(h/R)^d\ln n} + \frac{1}{\sqrt{n}} \right)
\end{align*}
for some $C>0$ depending on $s,p,d,\kappa(\Psi),\Psi(1),L$, where the hidden constant in the $O(n^{-1/2})$ term is thanks to \eqref{eq:orlicz_varentropy} in Lemma \ref{lemma:orlicz_tail}. Hence, plugging in the bandwidth in \eqref{eq:orlicz_bandwidth} gives the desired upper bound. In the general case, we may simply truncate the density $f$ onto the set $\{x: |x|\le R\}$ and apply the same analysis\footnote{The density after truncation is approximately a density, which integrates to $1+O(n^{-4})$ following the same line of proofs in Lemma \ref{lemma:orlicz_tail}.}, with the additional error upper bounded in Lemma \ref{lemma:orlicz_tail} which is negligible.

\subsubsection{Proof of lower bound}
The proof of the lower bound is similar to Section \ref{sec.lowerbound}, with a slightly different construction. Since the parametric rate $\Omega(n^{-1/2})$ is trivial, we only need to show the first term. Set
\begin{align}\label{eq:h_and_S}
h = (n\ln n)^{-\frac{1}{s+d}}\cdot [\Psi^{-1}(n)]^{\frac{d}{p(s+d)}}, \qquad R=c_0\Psi^{-1}(n), \qquad S = \frac{R^d-1}{h^d}
\end{align}
where $c_0>0$ is a small numerical constant to be specified later, and without loss of generality we assume that $S$ is an integer. Fix $\alpha\in (0,S^{-1})$, and constants $d_1,d_2,d_3$ as in Section \ref{sec.lowerbound}.

Let $I_1,\ldots,I_S$ be the partition of $[0,R]^d - [0,1]^d$ into $S$ smaller cubes of edge length $h$, and $t_i$ denote the leftmost corner of sub-cube $I_i$, i.e., $I_i=t_i+[0,h]^d$. 
Let $g$ and $w$ be some fixed smooth probability density functions on $\reals^d$ with finite entropy both of which vanish outside $[0,1]^d$. Note that the smoothness of $g$ and $w$ on $\mathbb{R}^d$ implies that their derivatives of any order vanish at their respective boundary.

To each vector $P=(p_1,\ldots,p_S)\in [0,\frac{d_3\ln n}{n}]^S$, we associate a function (also appeared in \eqref{eq.f_P})
\begin{align*}
f_P(t) \triangleq \sum_{i=1}^S p_i\cdot \frac{1}{h^d}g\pth{\frac{t-t_i}{h}} + \left(1-S \alpha\right)w(t).
\end{align*}
By the choice of $\alpha$, $f_P$ is non-negative everywhere. We show that with a proper choice of the vector $P$, the density $f_P$ lies in the Lipschitz ball $\Lip_{s,p,d}^{\Psi}(L_0)$ for a proper constant $L_0$. 
\begin{lemma}\label{lemma.subGaussian_Lipschitz}
Let $\rho>0$. There exists a numerical constant $c_0>0$ depending on $\kappa(\Psi)$ such that if 
\begin{align*}
\frac{1}{S}\sum_{i=1}^S p_i^p \le \left(\frac{2C_1}{n\ln n}\right)^p, 
\end{align*}
then $f_P\in \Lip_{s,p,d}^{\Psi}(L_0)$ for some constant $L_0>0$ independent of $n$. 
\end{lemma}
\begin{proof}
	We first verify the Lipschitz norm condition, whose proof is similar to Lemma \ref{lmm:fP}. Define $\tilde{f}_P$ as in Lemma \ref{lmm:fP}, it suffices to show $\|\tilde f_P\|_{\Lip_{s,p,d}}\le L_0$ for some constant $L_0$ independent of $n$.
	Observe that
	\begin{equation}\label{eq:tfP_new}
		\begin{aligned}
	\|\tilde{f}_P\|_p &= \frac{\|g\|_p}{h^d}\left(h^d\sum_{i=1}^Sp_i^p\right)^{\frac{1}{p}} \\
	&\le \frac{\|g\|_pR^{\frac{d}{p}}}{2h^d}\left(\frac{1}{S}\sum_{i=1}^Sp_i^p\right)^{\frac{1}{p}} \\
	&= \frac{c_0^{\frac{d}{p}}h^s\|g\|_p}{2}\cdot n\ln n\left(\frac{1}{S}\sum_{i=1}^Sp_i^p\right)^{\frac{1}{p}},
	\end{aligned}
	\end{equation}
	the condition ensures that $h^{-s}\|\tilde{f}_P\|_p$ is upper bounded by a constant depending only on $c_0, C_1$ and $g$. 
	By \cref{eq:Kequiv1} in \Cref{lmm:Kequiv}, \prettyref{eq:tfP_new} implies that there exists a constant $L_0$ independent of $n$ such that $\omega_r(\tilde{f}_P,t)\le L_0t^s$ for any $t\ge h$.
	For $t<h$, further observe that for $r=\lceil s\rceil$, we have
	\begin{equation}	\label{eq:rfP_new}
	\begin{aligned}
	\|D^r\tilde{f}_P\|_p &= \frac{\|D^r g\|_p}{h^{d+r}}\left(h^d\sum_{i=1}^Sp_i^p\right)^{\frac{1}{p}} \\
	&\le \frac{\|D^r g\|_pR^{\frac{d}{p}}}{2h^{d+r}}\left(\frac{1}{S}\sum_{i=1}^Sp_i^p\right)^{\frac{1}{p}} \\
	&= \frac{c_0^{\frac{d}{p}}\|D^r g\|_p}{2h^{r-s}}\cdot n\ln n\left(\frac{1}{S}\sum_{i=1}^Sp_i^p\right)^{\frac{1}{p}}
	\end{aligned}
	\end{equation}
	Hence, $h^{r-s}\|D^r\tilde{f}_P\|_p$ is upper bounded by a constant depending on $C_1,g$. 
	By \cref{eq:Kequiv2} in \Cref{lmm:Kequiv}, \prettyref{eq:rfP_new} implies that $\omega_r(\tilde{f}_P,t)\le L_0t^s$ for any $t\le h$, completing the proof of $\|\tilde f_P\|_{\Lip_{s,p,d}}\le L_0$. 
	
	Next we show that $f_P$ satisfies the Orlicz norm condition. Note that the density $f_P$ is supported on $[0,R]^d$, and
	\begin{align*}
	  f_P(x) \le \frac{d_3\ln n}{nh^d}\cdot \|g\|_\infty + \|w\|_\infty\cdot \mathbbm{1}(x\in [0,1]^d).
	\end{align*}
	As a result, there exists a numerical constant $C$ depending on $d,d_3,\Psi(1), \|g\|_\infty, \|w\|_\infty$ such that 
	\begin{align}\label{eq:subGaussian_norm}
	\int_{\bR^d} \Psi(|x|) f_P(x)dx \le C\left(1+ \Psi\left(c_0\Psi^{-1}(n) \right)\cdot \frac{R^d\ln n}{nh^d}\right). 
	\end{align}
	Since $\Psi$ is of rapid growth, for $c_0 = \kappa(\Psi)^{-m}$ with $m\in\naturals$ we have $\Psi(c_0\Psi^{-1}(n))\le n^{2^{-m}}$. Hence, choosing $m$ large enough such that $2^{-m} <  \frac{s}{s+d}$, the RHS of \eqref{eq:subGaussian_norm} is upper bounded by a numerical constant $L_0$ independent of $n$, completing the proof.  
\end{proof}

By Lemma \ref{lemma.subGaussian_Lipschitz}, the construction of two measures in Section \ref{subsec.measure} still ensures that $f_P$ lies in the Lipschitz ball $\Lip_{s,p,d}^{\Psi}(L_0)$. Consequently, the remaining arguments (reduction to parametric submodel, fuzzy hypothesis testing, etc.) in Section \ref{sec.lowerbound} can be applied to obtain
\begin{align*}
\inf_{\hat{H}} \sup_{f\in \Lip_{s,p,d}^{\Psi}(L_0)} \bE_f \left(\hat{H} - H(f)\right)^{\frac{1}{2}} \gtrsim \frac{S}{n\ln n} \asymp (n\ln n)^{-\frac{s}{s+d}}\cdot [\Psi^{-1}(n)]^{d\left(1-\frac{d}{p(s+d)}\right)},
\end{align*}
as desired. \qed

\subsection{Proof of Theorem~\ref{lem.plug-in}}
\subsubsection{Proof of upper bound}
We first prove the achievability part of the plug-in approach. For kernel $K(\cdot)$ satisfying \prettyref{ass:kernel} and $\hat{f}_h(x)$ given by \eqref{eq:KDE}, Lemma \ref{lemma.first_approx_err} shows that the first approximation error is at most
\begin{align}\label{eq:plugin_error_1}
|H(f) - H(f_h)| \lesssim Lh^s.
\end{align}

Next we upper bound the bias $|\bE H(\hat{f}_h) - H(f_h)|$ of the plug-in estimator. Let $\phi(x)\triangleq x\ln x$, and $g(x)$ be any twice continuously differentiable function on $[0,\infty)$. By \cite[Eqn. (2.5)]{Totik1994approximation}, for any $u,v\ge 0$,
\begin{align*}
|g(v) - g(u) - g'(u)(v-u)| \le \frac{4(v-u)^2}{u} \|wg''(w)\|_\infty.
\end{align*}
Choosing $u=f_h(x), v=\hat{f}_h(x)$ and taking expectation at both sides, and noting that
\begin{align*}
\bE(\hat{f}_h(x)-f_h(x))^2 \le \frac{\|K\|_{\infty}f_h(x)}{nh^d}
\end{align*} 
we have 
\begin{align*}
|\bE \phi(\hat{f}_h(x)) - \phi(f_h(x))| &\le 2\|\phi-g\|_\infty + |\bE g(\hat{f}_h(x))-g(f_h(x))| \\
&\le 2\|\phi-g\|_\infty + \frac{4\|wg''(w)\|_\infty}{f_h(x)}\cdot \bE(\hat{f}_h(x)-f_h(x))^2\\
&\lesssim \|\phi-g\|_\infty + \frac{1}{nh^d}\cdot \|wg''(w)\|_\infty.
\end{align*}

The previous inequality holds for any $g$, so we may take infimum over $g$ and reduce the problem to the following lemma:
\begin{lemma}[{\cite[Theorem 2.1.1]{Ditzian--Totik1987}}]
	\label{lemma.equivalence_K-func_moduli}
	For any real-valued function $\phi$ defined on $[0,\infty)$ and $t>0$, there exists a universal constant $C>0$ independent of $\phi,t$ such that
	\begin{align*}
	\inf_{g\in C^2[0,\infty)} \left(\|\phi-g\|_\infty + t^2\|wg''(w)\|_\infty \right) \le C\omega_\varphi^2(f,t)_\infty 
	\end{align*}
	where $\omega_\varphi^2(f,t)_\infty$ is the second-order modulus of smoothness defined in \eqref{eq.DT_modulus} with $\varphi(x)=\sqrt{x}$.
\end{lemma}
It follows from \cite[Example 3.4.2]{Ditzian--Totik1987} that for $\phi(x)=x\ln x$, we have
\begin{align*}
\omega_\varphi^2(\phi,t)_\infty \asymp t^2.
\end{align*}
As a result, by the previous arguments and Lemma \ref{lemma.equivalence_K-func_moduli} we obtain
\begin{align*}
|\bE \phi(\hat{f}_h(x))-\phi(f_h(x))| \lesssim \frac{1}{nh^d},
\end{align*}
which integrates into
\begin{align}\label{eq:plugin_error_2}
|\bE H(\hat{f}_h) - H(f_h)| \le \int_{[0,1]^d} |\bE \phi(\hat{f}_h(x))-\phi(f_h(x))| dx \lesssim \frac{1}{nh^d}. 
\end{align}

To upper bound the variance $\var(H(\hat{f}_h))$, we apply the same arguments in Section \ref{subsec.overall} to obtain
\begin{align*}
\var(H(\hat{f}_h)) \lesssim nh^d\cdot \int_{[0,1]^d} \bE[(\hat{f}_h(x)\ln \hat{f}_h(x) - \hat{f}_h'(x)\ln \hat{f}_h'(x) )^2]dx, 
\end{align*}
where $\hat{f}_h'(x)$ is the density estimate based on observations $X_1', X_2,\ldots,X_n$, with $X_1'$ being an independent copy of $X_1$. For any $x\in [0,1]^d$, we split into two cases: 
\begin{enumerate}
	\item Case I: $f_h(x) \le \frac{c_1\ln n}{nh^d}$, where $c_1>0$ is the constant appearing in Lemma \ref{lemma.regime}. By Lemma \ref{lemma.regime} and the union bound, with probability at least $1-2n^{-5d}$ we have $\hat{f}_h(x), \hat{f}_h'(x) \le \frac{2c_1\ln n}{nh^d}$. Note that for any $\varepsilon\in (0,1)$, there exists some constant $C_\varepsilon>0$ such that $|x\ln x-y\ln y|\le C_{\varepsilon}|x-y|^{1-\varepsilon}$ whenever $x,y\in [0,\frac{1}{2}]$. Hence,
	\begin{align*}
	\bE(\hat{f}_h(x)\ln \hat{f}_h(x) - \hat{f}_h'(x)\ln \hat{f}_h'(x) )^2 &\le C_{\varepsilon}^2\bE|\hat{f}_h(x) - \hat{f}_h'(x)|^{2(1-\varepsilon)} \\
	&\quad + 2n^{-5d}\cdot \left(2\cdot \frac{\|K\|_\infty}{h^d}\ln \frac{\|K\|_\infty}{h^d}\right)^2. 
	\end{align*}
	Note that for $\varepsilon \in(0,\frac{1}{2})$, we have
	\begin{align*}
	\bE|\hat{f}_h(x) - \hat{f}_h'(x)|^{2(1-\varepsilon)} &= (nh^d)^{-2(1-\varepsilon)}\cdot \bE\left|K\left(\frac{x-X_1}{h}\right) - K\left(\frac{x-X_1'}{h}\right)\right|^{2(1-\varepsilon)} \\
	&\lesssim (nh^d)^{-2(1-\varepsilon)}\cdot \bE K^{2(1-\varepsilon)}\left(\frac{x-X_1}{h}\right) \\
	&\lesssim (nh^d)^{-2(1-\varepsilon)}\cdot h^df_h(x) \lesssim \frac{\ln n}{n^3h^{2d}}\cdot (nh^d)^{2\varepsilon}, 
	\end{align*}
	where the last inequality follows from $f_h(x)\le \frac{c_1\ln n}{nh^d}$. 
	\item Case II: $f_h(x) > \frac{c_1\ln n}{nh^d}$. By Lemma \ref{lemma.regime} and the union bound, with probability at least $1-2n^{-5d}$ we have $\hat{f}_h(x), \hat{f}_h'(x) \in [\frac{f_h(x)}{2},2f_h(x)]$. Then the mean value theorem gives
	\begin{align*}
	\bE(\hat{f}_h(x)\ln \hat{f}_h(x) - \hat{f}_h'(x)\ln \hat{f}_h'(x) )^2 & \lesssim (1+(\ln f_h(x))^2)\cdot \bE(\hat{f}_h(x) - \hat{f}_h'(x))^2 \\
	&\qquad + 2n^{-5d}\cdot \left(2\cdot \frac{\|K\|_\infty}{h^d}\ln \frac{\|K\|_\infty}{h^d}\right)^2 \\
	&\lesssim \frac{f_h(x)(1+(\ln f_h(x))^2)}{n^2h^d},
	\end{align*}
	where in the last step we have used that 
	\begin{align*}
	\bE(\hat{f}_h(x) - \hat{f}_h'(x))^2 &= \frac{1}{n^2h^{2d}}\cdot \bE\left( K\left(\frac{x-X_1}{h}\right) - K\left(\frac{x-X_1'}{h}\right) \right)^2 \\
	& \le \frac{4}{n^2h^{2d}}\cdot \bE K^2\left(\frac{x-X_1}{h}\right) \le \frac{4\|K\|_\infty}{n^2h^d}\cdot f_h(x). 
	\end{align*}
\end{enumerate}

Combining the previous two cases, for any $x\in [0,1]^d$ and $\varepsilon\in(0,\frac{1}{2})$ we have
\begin{align*}
\bE(\hat{f}_h(x)\ln \hat{f}_h(x) - \hat{f}_h'(x)\ln \hat{f}_h'(x) )^2 \lesssim \frac{\ln n}{n^3h^{2d}}\cdot (nh^d)^{2\varepsilon} + \frac{f_h(x)(1+(\ln f_h(x))^2)}{n^2h^d}. 
\end{align*}
Now by Lemma \ref{lemma:varentropy}, an integration yields to
\begin{equation}\label{eq:plugin_error_3}
\begin{aligned}
\var(H(\hat{f}_h)) &\lesssim \frac{\ln n}{n^2h^{d}}\cdot (nh^d)^{2\varepsilon} + \int_{[0,1]^d}\frac{f_h(x)(1+(\ln f_h(x))^2)}{n}dx \\
&\lesssim \frac{\ln n}{n^2h^{d}}\cdot (nh^d)^{2\varepsilon} + \frac{(\ln L)^2}{n}. 
\end{aligned}
\end{equation}

In summary, combining \eqref{eq:plugin_error_1}, \eqref{eq:plugin_error_2} and \eqref{eq:plugin_error_3}, we conclude that
\begin{align*}
\left[\sup_{f\in \text{\rm Lip}_{s,p,d}(L)} \bE_f\left(H(\hat{f}_h) - H(f)\right)^2\right]^{\frac{1}{2}} \lesssim Lh^s + \frac{1}{nh^d} + \frac{\ln L}{\sqrt{n}} + \frac{\sqrt{\ln n}}{n\sqrt{h^d}}\cdot (nh^d)^{\varepsilon}.
\end{align*}
Choosing $h\asymp (Ln)^{-\frac{1}{s+d}}$ and $\varepsilon>0$ sufficiently small, we arrive at the desired upper bound. \qed

%

\subsubsection{Proof of lower bound}
Next we establish the lower bound. The parametric rate $\Omega(n^{-1/2}\ln L)$ has been established in Theorem \ref{thm.lower_bound}, so we focus only on the $\Omega(n^{-s/(s+d)}L^{d/(s+d)})$ part. Suppose that some plug-in approach with kernel $K(\cdot)$ and bandwidth $h$ attains a worst-case $L_2$ risk $o(n^{-s/(s+d)}L^{d/(s+d)})$. We make the following claims:
\begin{enumerate}
\item We must have $h\gg (Ln)^{-\frac{1}{s+d}}$. Otherwise, consider $f\equiv 1$, we have $f_h\equiv 1$ and
\begin{align*}
\bE (\hat{f}_h(x) - f_h(x))^2 = \frac{1}{n}\left(\frac{1}{h^d}\int_{\mathbb{R}^d} K^2(t)dt -1\right) \gtrsim \frac{1}{nh^d}.
\end{align*}
By Lemma \ref{lemma.regime}, $\hat{f}_h(x)\le 2$ with probability at least $1-n^{-5d}$, Taylor expansion with Lagrange remainder term give
\begin{align*}
&\bE[\hat{f}_h(x)\ln \hat{f}_h(x)] - f_h(x)\ln f_h(x) \nonumber\\
&\ge \frac{1}{2\cdot 2} \bE(\hat{f}_h(x)-f_h(x))^2 - n^{-5d}\cdot \left(\frac{\|K\|_\infty f_h(x)}{nh}\ln \frac{\|K\|_\infty f_h(x)}{nh} + |f_h(x)\ln f_h(x)|\right)\\
&\gtrsim \frac{1}{nh^d}.
\end{align*}
As a result, since $h\lesssim (Ln)^{-1/(s+d)}$, the overall bias is at least
\begin{align*}
H(f) - \bE H(\hat{f}_h) &= \int_{[0,1]^d} \left(\bE \hat{f}_h(x)\ln \hat{f}_h(x) - f_h(x)\ln f_h(x)\right) dx \\
&\gtrsim \frac{1}{nh^d} \gtrsim n^{-\frac{s}{s+d}}L^{\frac{d}{s+d}}, 
\end{align*}
a contradiction to the assumed performance of $H(\hat{f}_h)$.
\item We must have $h\lesssim (Ln)^{-1/(s+d)}$. Otherwise, $h\gg (Ln)^{-1/(s+d)}$, and by \eqref{eq:plugin_error_2} and triangle inequality, for any $f\in \Lip_{s,p,d}(L)$ we must have
\begin{align*}
|H(f)-H(f_h)| &\le |H(f) - \bE H(\hat{f}_h)| + |\bE H(\hat{f}_h) - H(f_h)| \\
&\lesssim |H(f) - \bE H(\hat{f}_h)| + \frac{1}{nh^d} \\
&\ll n^{-\frac{s}{s+d}}L^{\frac{d}{s+d}} \\
&\lesssim Lh^s \wedge 1.
\end{align*}
In the sequel we will construct some $f\in \Lip_{s,p,d}(L)$ such that $|H(f)-H(f_h)|\gtrsim Lh^s\wedge 1$, establishing the desired contradiction. Let $f$ be the density $f_P$ which was constructed in \eqref{eq.f_P} with $P$ consisting of all identical elements (denoted by $p$), or specifically,
\begin{align*}
f(t) = p\cdot \sum_{i=1}^S Lh^s g\pth{\frac{t-t_i}{h}} + \left(1-pLh^s\right)w(t) \triangleq f_1(t) + f_2(t).
\end{align*} 
We choose the parameter $p = \Theta( (Lh^s)^{-1}\wedge 1)$ so that $f$ is a density, and choose the smooth function $g$ to satisfy $H(g)\neq H(g_1)$ for $g_1(x)\triangleq \int_{[0,1]^d} K(x-y)g(y)dy$. Note that the density $f\in \text{H}_{d}^s(L)$ lies in the H\"{o}lder ball for small $p$. By the structure of the kernel-smoothed density, we may write
\begin{align*}
f_h(t) = p\cdot \sum_{i=1}^S Lh^s g_1\pth{\frac{t-t_i}{h}} + f_{2h}(t) + \varepsilon(t)
\end{align*}
where $f_{2h}(t)$ is the kernel-smoothed version of $f_2(t)$, and $\varepsilon(t)$ can be non-zero only in cubes which intersect the boundary of $[\frac{1}{4},\frac{3}{4}]^d$. Moreover, since $f$ lies in the H\"{o}lder ball, we have $\|f-f_h\|_\infty \lesssim Lh^s$ and $\|\varepsilon\|_\infty\lesssim Lh^s$. Noting that the total volume of all cubes which intersect the boundary is $O(h)$, we have
\begin{align*}
H(f_h) - H(f) \gtrsim (H(f_{2h}) - H(f_2)) + p\cdot Lh^s (H(g_1)-H(g)) - Lh^s\ln(1/h)\cdot h.
\end{align*}
Finally, since convolution increases the entropy, we have $H(f_{2h})\ge H(f_2), H(g_1)\ge H(g)$. Now our choice of $p=\Theta((Lh^s)^{-1}\wedge 1)$ gives
\begin{align*}
H(f_h) - H(f) \gtrsim p\cdot Lh^s - Lh^{s+1}\ln (1/h) \gtrsim Lh^s \wedge 1, 
\end{align*}
as desired.
\end{enumerate}
The proof of the lower bound is complete noting that these two claims contradict each other.
\qed


\subsection{Proof of Theorem \ref{thm:periodic_boundary}}\label{subsec.periodic_boundary_proof}
The lower bound directly follows from Theorem \ref{thm:main} since we are considering a larger density class. For the upper bound, we define the estimator $\hat{H}$ as in Section \ref{sec.construction} except that the kernel convolution (i.e., the definition of $f_h$ in \eqref{eq.f_h}) is understood on the torus $\torus^d$. The analyses parallel those in Section \ref{sec.analysis}, where Lemmas \ref{lemma.smooth_approx_g} and \ref{lemma.maximal_inequality} now need some caution when dealing with the torus $\torus^d$. For Lemma \ref{lemma.smooth_approx_g}, same result holds simply by changing $\bR^d$ into $\torus^d$ in the proof of Lemma \ref{lemma.smooth_approx_g}. As for Lemma \ref{lemma.maximal_inequality}, we remark that the Hardy--Littlewood maximal inequality also holds on the torus\footnote{For example, see \cite{tao-torus}; the proof follows from that of standard maximal inequality on Euclidean space via the Vitali covering lemma \cite[Chapter 1, Theorem 1]{stein2016singular}.}. \qed

\section{Proof of main lemmas}

\subsection{Proof of Lemma \ref{lemma.smooth_approx_g}}
	Let $K_d$ be the following non-negative kernel on $\mathbb{R}^d$:
	\begin{align*}
	K_d(x) = \prod_{i=1}^d\left(1-|x_i|\right)_+, \qquad x=(x_1,\ldots,x_d)
	\end{align*}
	where we note that $\int_{\mathbb{R}^d} K_d(x)dx=1$. Now define
	\begin{align*}
	K_{d,h}(x) = \frac{1}{h^d} K_d(\frac{x}{h}),
	\end{align*}
	and
	\begin{align*}
	g(x) = \int_{\bR^d} f(y)K_{d,h}(x-y)dy.
	\end{align*}
	
	Since both $f$ and $K_{d,h}$ are non-negative, so is $g$. Also, by the symmetry of $K_{d}(\cdot)$, 
	\begin{align*}
	f(x) - g(x) &= \int_{\bR^d} (f(x)-f(y))K_{d,h}(x-y)dy \\
	&= \int_{\bR^d} (f(x)-f(x-hu))K_d(u)du \\
	&= \frac{1}{2} \int_{\bR^d} (2f(x)-f(x-hu)-f(x+hu))K_d(u)du
	\end{align*}
	As a result,
	\begin{align*}
	\|f-g\|_p &\le \frac{1}{2}\int_{\bR^d} \|2f(x)-f(x-hu)-f(x+hu)\|_p K_d(u)du\\
	&\lesssim Lh^s\cdot \frac{1}{2}\int_{\bR^d} K_d(u)du \asymp Lh^s.
	\end{align*}
	
	As for the Hessian of $g$, we have
	\begin{align*}
	[\nabla^2 g(x)]_{ii} &= \frac{\partial^2}{\partial x_i^2} \int_{\bR^d} f(y)K_{d,h}(x-y)dy \\
	&= \int_{\bR^d} f(y)\frac{\partial^2}{\partial x_i^2} K_{d,h}(x-y)dy \\
	&= \frac{1}{h^2}\int_{\mathbb{R}^d} f(x-hu) \frac{\partial^2}{\partial u_i^2} K_{d}(u)du \\
	&= \frac{1}{h^2}\int_{\mathbb{R}^d} f(x-hu) \cdot (\delta(u_i+1)-2\delta(u_i)+\delta(u_i-1))K_{d-1}(u_{\backslash i})du \\
	&= \frac{1}{h^2}\int_{\mathbb{R}^{d-1}} \Delta_{he_i}^2 f(x_{\backslash i}-hu_{\backslash i}, x_i)K_{d-1}(u_{\backslash i})du_{\backslash i}
	\end{align*}
	where $\Delta_{he_i}^2 f(x) \triangleq f(x+he_i)-2f(x)+f(x-he_i)$, with $e_i$ being the $i$-th coordinate vector, and $\delta(\cdot)$ is the delta function.
	As a result, we have
	\begin{align*}
	\|[\nabla^2 g(x)]_{ii}\|_p &\le \frac{1}{h^2}\int_{ \mathbb{R}^{d-1}} \|\Delta_{he_i}^2 f(x_{\backslash i}-hu_{\backslash i}, x_i)\|_pK_{d-1}(u_{\backslash i})du_{\backslash i} \\
	&= \frac{1}{h^2}\int_{ \mathbb{R}^{d-1}} \|\Delta_{he_i}^2 f(x)\|_pK_{d-1}(u_{\backslash i})du_{\backslash i}  \\
	&\lesssim \frac{1}{h^2}\int_{ \mathbb{R}^{d-1}} Lh^s\cdot K_{d-1}(u_{\backslash i})du_{\backslash i} \asymp Lh^{s-2}.
	\end{align*}
	
	Similarly, for $i\neq j$, we have
	\begin{align*}
	[\nabla^2 g(x)]_{ij} &= \frac{1}{h^2}\int_{\mathbb{R}^d} f(x-hu) \frac{\partial^2}{\partial u_i\partial u_j} K_{d}(u)du \\
	&= \frac{1}{h^2}\int_{\mathbb{R}^d} f(x-hu) \cdot (\mathbbm{1}_{u_i\in (0,1)}-\mathbbm{1}_{u_i\in (-1,0)})(\mathbbm{1}_{u_j\in (0,1)}-\mathbbm{1}_{u_j\in (-1,0)})K_{d-2}(u_{\backslash (i,j)})du \\
	&= \frac{1}{h^2}\int_{\mathbb{R}^{d-2}}\int_0^1\int_0^1 \Delta_{he_i}\Delta_{he_j} f(x-hu)K_{d-2}(u_{\backslash (i,j)})du_idu_jdu_{\backslash (i,j)}
	\end{align*}
	where $\Delta_{he_i}f(x)= f(x+he_i)-f(x)$ denotes the forward difference. Note that
	\begin{align*}
	\Delta_{he_i}\Delta_{he_j}f(x) &= f(x+he_i+he_j) - f(x+he_i) - f(x+he_j) + f(x) \\
	&= \frac{f(x+he_i+he_j)-2f(x+he_i)+f(x-he_i+he_j)}{2} \nonumber \\
	&\qquad + \frac{f(x+he_i+he_j)-2f(x+he_j)+f(x-he_j+he_i)}{2} \nonumber\\
	&\qquad - \frac{f(x-he_i+he_j) - 2f(x) + f(x-he_j+he_i)}{2} \\
	&= \frac{\Delta_{he_i}^2f(x+he_i) + \Delta_{he_j}^2f(x+he_j)- \Delta_{h(e_i+e_j)}^2 f(x) }{2}
	\end{align*}
	we conclude that
	\begin{align*}
	\|\Delta_{he_i}\Delta_{he_j}f \|_p &\le \frac{1}{2}\left( \|\Delta_{he_i}^2 f\|_p + \|\Delta_{he_j}^2 f\|_p + \|\Delta_{h(e_i+e_j)}^2 f\|_p\right) 
	\lesssim Lh^s.
	\end{align*}
	Then using the same technique, we also have
	\begin{align*}
	\|[\nabla^2 g(x)]_{ij}\|_p \lesssim Lh^{s-2}.
	\end{align*}
	
	Finally, note that $\|A\|_{\text{op}}\le \sqrt{\text{Tr}(A^TA)}\le \sum_{i,j} |A_{ij}|$, we conclude that
	\begin{align*}
	\|\|\nabla^2 g\|_{\text{op}}\|_p &\le \left\| \sum_{i,j=1}^d |[\nabla^2 g(x)]_{ij}| \right\|_p
	\le \sum_{i,j=1}^d \|[\nabla^2 g(x)]_{ij}\|_p 
	\lesssim Lh^{s-2}, 
	\end{align*}
	as desired.\qed

\subsection{Proof of Lemma \ref{lemma.regime}}
Applying the Bennett inequality (Lemma \ref{lemma.bennett}) to independent random variables $Y_i=K_h(x-X_i)$, we obtain
\begin{align*}
&\bP\left(\left|\hat{f}_h(x)-f_h(x)\right|\ge\varepsilon \right)\\
&\le 2\exp\left(-\frac{n\varepsilon^2}{2\var(K_h(x-X_1))+2\|K\|_\infty\varepsilon/3h^d}\right).
\end{align*}
For the variance, by the non-negativity of the kernel $K$ and the density $f$, we have
\begin{equation}\label{eq:varKh}
\begin{aligned}
\var(K_h(x-X_1)) &\le \bE[K_h(x-X_1)^2]\\
&= \int K_h(x-y)^2f(y)dy\\
&\le \frac{\|K\|_\infty}{h^d}\int K_h(x-y)f(y)dy\\
&= \frac{\|K\|_\infty}{h^d}f_h(x). 
\end{aligned}
\end{equation}
Combining these two inequalities yields the desired result. \qed

\subsection{Proof of Lemma \ref{lemma.non-smooth}}
Recall that $\Expect[\hat{H}_1(x)]=Q(f_h(x))$, where 
$Q(t)$ defined in \cref{eq:Q} is the best degree-$k$ polynomial that uniformly approximates $-t \ln t$ on the interval $[0,\frac{2c_1\ln n}{nh^d}]$. Thus the bias is upper bounded by
\begin{align*}
|\bE[\hat{H}_1(x)] + f_h(x)\ln f_h(x)| &= \left|\sum_{l=0}^k a_lf_h(x)^l +  f_h(x)\ln f_h(x)\right| \\
&\lesssim \frac{1}{k^2}\cdot \frac{2c_1\ln n}{nh^d} \lesssim \frac{1}{nh^d\ln n}
\end{align*}
provided that $f_h(x) \leq \frac{2c_1\ln n}{nh^d}$.

To upper bound the second moment of the perturbation, first note that $\hat{H}_1(x)$ does not depend on observations from $X^{(3)}$ and is symmetric in observations from $X^{(2)}$. Consequently, we may assume that the original estimator $\hat H_1(x)$ uses the observations $X^{(2)}=(X_1,\ldots,X_n)$ and the perturbed estimator $\hat{H}_1'(x)$ uses observations $(X_1',\ldots,X_n')$, where $X_j'=X_j$ for $j=1,\ldots,n-1$ and $X_n'$ is an independent copy of $X_n$. Then
\begin{align*}
\hat{H}_1(x)  - \hat{H}_1'(x)
= & ~ \sum_{l=0}^k a_l \frac{1}{\binom{n}{l}} \sum_{J\in \binom{[n]}{l}} \left(\prod_{j\in J} K_h(x-X_{j}) -\prod_{j\in J} K_h(x-X_{j}') \right)  \\
= & ~ \frac{1}{n} (K_h(x-X_{n})-K_h(x-X_{n}')) \\
&\qquad \cdot \underbrace{\sum_{l=1}^k l a_l \frac{1}{\binom{n-1}{l-1}} \sum_{J'\in \binom{[n-1]}{l-1}} \prod_{j\in J'} K_h(x-X_{j})}_{\triangleq \tilde{H}_1(x)}.
\end{align*}
By the independence of $(X_n, X_n')$ and $(X_1,\ldots,X_{n-1})$, 
\begin{align}\label{eq.perturbation_H1}
\bE[(\hat{H}_1(x)  - \hat{H}_1'(x))^2] = \frac{1}{n^2}\bE[(K_h(x-X_{n})-K_h(x-X_{n}'))^2]\cdot \bE[\tilde{H}_1(x)^2]. 
\end{align}
By \eqref{eq:varKh} and the assumption on $f_h(x)$, we have
\begin{equation}\label{eq.H1_part1}
\begin{aligned}
\bE[(K_h(x-X_{n})-K_h(x-X_{n}'))^2] &= 2\var(K_h(x-X_n)) \\
&\le \frac{2\|K\|_\infty f_h(x)}{h^d}\le \frac{4c_1\|K\|_\infty \ln n}{nh^{2d}}.
\end{aligned}
\end{equation}
It remains to upper bound $\bE[\tilde{H}_1(x)^2]$. By Cauchy-Schwarz, 
\begin{align}
\bE[\tilde{H}_1(x)^2]
\le k^3\cdot \sum_{l=0}^{k-1} \frac{a_{l+1}^2}{\binom{n-1}{l}^{2}}\sum_{|J|=l}\sum_{|T|=l} \bE\left[\prod_{j\in J}\prod_{t\in T} K_h(x-X_{j})K_h(x-X_{t})\right] \label{eq:varH11}.
\end{align}
Recall from \ref{eq:varKh} that
\begin{align*}
\bE K_h(x-X_i) &= f_h(x), \\
\bE K_h(x-X_i)^2 & \le \frac{\|K\|_\infty f_h(x)}{h^d}.
\end{align*}
Therefore, if $|J\cap T|=r$, we have
\begin{align*}
\bE\left[\prod_{j\in J} \prod_{t\in T}  K_h(x-X_j)K_h(x-X_t)\right] &\le \left(\frac{\|K\|_\infty f_h(x)}{h^d}\right)^r \cdot f_h(x)^{2(l-r)} \nonumber \\
&= \|K\|_\infty^rh^{-dr} f_h(x)^{2l-r}.
\end{align*}
Moreover, the number of pairs of subsets $J, T\in\binom{[n-1]}{l}$ such that $|J\cap T|=r$ is $\frac{(n-1)!}{r!((l-r)!)^2(n-1-2l+r)!}$.
Hence, for $0\leq l\le k-1$,
\begin{equation}\label{eq:varH12}
\begin{aligned}
&\binom{n-1}{l}^{-2}\sum_{|J|=l, |T|=l} \bE\left[\prod_{j\in J} \prod_{t\in T}  K_h(x-X_j)K_h(x-X_t)\right]  \\
&\le \frac{(l!)^2((n-1-l)!)^2}{(n-1)!} \sum_{r=0}^l \frac{\|K\|_\infty^rh^{-dr}f_h(x)^{2l-r} }{r!((l-r)!)^2(n-1-2l+r)!}  \\
&= \frac{((n-1-l)!)^2}{(n-1)!(n-1-2l)!}\sum_{r=0}^l \binom{l}{r}^2\cdot \frac{r!  |K\|_\infty^rh^{-dr}f_h(x)^{2l-r} }{(n-1-2l+r)(n-2-2l+r)\cdots(n-2l)} \\
&\le \sum_{r=0}^l 2^{2l} \left(\frac{r}{n-2l}\right)^r
\|K\|_\infty^rh^{-dr}f_h(x)^{2l-r}  \\
&\le 2^{2l}\sum_{r=0}^l \left(\frac{2k}{nh^d}\right)^r\|K\|_\infty^rf_h(x)^{2l-r}  \\
&\le 2^{2l}\sum_{r=0}^l \left(\frac{2\|K\|_\infty c_2\ln n}{nh^d}\right)^r\left(\frac{2c_1\ln n}{nh^d}\right)^{2l-r} \\
&= \left(\frac{2c_1\ln n}{nh^d}\right)^{2l} \sum_{r=0}^l \left(\frac{\|K\|_\infty c_2}{c_1}\right)^r  \\
&\le 2\left(\frac{2c_1\ln n}{nh^d}\right)^{2l}, 
\end{aligned}
\end{equation}
where we have used our assumptions that $n\ge 4k= 4c_2\ln n$, and $c_1\ge 2\|K\|_\infty c_2$. By Lemma \ref{lemma.poly_approx}, we have the following coefficient bound:
\begin{equation}\label{eq:varH13}
\begin{aligned}
|a_{l+1}| &\le \left(\frac{2c_1\ln n}{nh^d}\right)^{-l}\ln\left(\frac{nh^d}{2c_1\ln n}\right)\cdot 2^{3c_2\ln n} 
\\
&\le \left(\frac{2c_1\ln n}{nh^d}\right)^{-l}\cdot n^{3c_2\ln 2}\ln n. 
\end{aligned}
\end{equation}
Combining \cref{eq:varH11}--\cref{eq:varH13} yields
\begin{equation}\label{eq.H1_part2}
\begin{aligned}
\bE[\tilde{H}_1(x)^2] &\le k^3\cdot \sum_{l=0}^{k-1} \left(\frac{2c_1\ln n}{nh^d}\right)^{-2l}(n^{3c_2\ln 2}\ln n)^2\cdot 2\left(\frac{2c_1\ln n}{nh^d}\right)^{2l} \\
&= 2k^4\cdot (n^{3c_2\ln 2}\ln n)^2 \lesssim n^\varepsilon, 
\end{aligned}
\end{equation}
provided that $7c_2\ln 2<\varepsilon$. Now combining \eqref{eq.perturbation_H1}, \eqref{eq.H1_part1} and \eqref{eq.H1_part2} completes the proof. 
\qed

\subsection{Proof of Lemma \ref{lemma.smooth}}
Recall that by Lemma \ref{lemma.regime}, with $c_1>0$ large enough we have
\begin{align}\label{eq.concentration}
\bP\left(\hat{f}_{h,2}(x) \notin \Big(\frac{f_h(x)}{2},2f_h(x)\Big)\right)\le n^{-5d}.
\end{align}
Moreover, since both $f_h$ and $\hat f_h$ are convolutions of probability measures with $K_h$, we have the deterministic upper bound: for all $x$,
\begin{align}
\max\{f_h(x),\hat{f}_{h,2}(x), \hat{f}_{h,3}(x)\} \le \|K_h\|_\infty = \frac{\|K\|_\infty}{h^d}.
\label{eq:fhsup}
\end{align}
Consequently, the following deterministic upper bound of $|\hat{H}_2(x)|$ (in \eqref{eq:H2x}) holds for all $x$: 
\begin{align}
|\hat{H}_2(x)| \lesssim 1 + nh^d + \frac{1}{h^d}\ln \frac{n}{h^d}. \label{eq:H2_sup}
\end{align}

Recall that throughout this proof we have by assumption
\begin{equation}
f_h(x)\ge \frac{c_1\ln n}{2nh^d}.
\label{eq:fh-smooth}
\end{equation}
For $\phi(z)=-z\ln z$ and $\bE_3$ denoting the partial expectation taken only with respect to the third group of observations $X^{(3)}$, we have
\begin{equation}\label{eq.bias}
\begin{aligned}
&\bE_3\left[\hat{H}_2(x)\bigg|\hat{f}_{h,2}(x)\ge \frac{f_h(x)}{2}\right] \\
&= \sum_{k=0}^2 \frac{\phi^{(k)}(\hat{f}_{h,2}(x))}{k!}(f_h(x)-\hat{f}_{h,2}(x))^k  \\
&= \phi(f_h(x)) - \frac{\phi^{(3)}(\xi_1)}{6}(f_h(x)-\hat{f}_{h,2}(x))^3 \\
&= \phi(f_h(x)) - \frac{\phi^{(3)}(f_h(x))}{6}(f_h(x)-\hat{f}_{h,2}(x))^3 \\
&\qquad + \frac{\phi^{(4)}(\xi_2)}{6}(f_h(x)-\hat{f}_{h,2}(x))^3(f_h(x)-\xi_1),
\end{aligned}
\end{equation}
where the last two steps follow from the Taylor expansion, with $\xi_1$ lying between $f_h(x)$ and $\hat{f}_{h,2}(x)$, and $\xi_2$ lying between $f_h(x)$ and $\xi_1$. The central moments of $\hat{f}_{h,2}(x)$ are computed as follows:
\begin{align*}
\Big|\bE(\hat{f}_{h,2}(x)-f_h(x))^3\Big| &= \left|\bE\left(\frac{1}{n}\sum_{i=1}^n K_h(x-X_i) - f_h(x)\right)^3\right| \\
&= \frac{1}{n^2} |\bE(K_h(x-X_1)-f_h(x))^3|\\
&\le \frac{1}{n^2}\left(\bE K_h(x-X_1)^3 + f_h(x)^3\right)\\
&\le \frac{1}{n^2}\left(\frac{\|K\|_\infty^2}{h^{2d}}\bE K_h(x-X_1) + f_h(x)^3\right)\\
&= \frac{f_h(x)}{n^2}\left(\frac{\|K\|_\infty^2}{h^{2d}} + f_h(x)^2\right) \overset{\prettyref{eq:fhsup}}{\le} \frac{2f_h(x)\|K\|_\infty^2}{n^2h^{2d}}, 
\end{align*}
and
\begin{align*}
\bE(\hat{f}_{h,2}(x)-f_h(x))^4 &= \bE\left(\frac{1}{n}\sum_{i=1}^n K_h(x-X_i) - f_h(x)\right)^4 \\
&\le \frac{3}{n^2}\bE^2[K_h(x-X_1)^2] + \frac{1}{n^3} \bE[K_h(x-X_1)^4]\\
&\le \frac{3}{n^2}\left(\frac{\|K\|_\infty f_h(x)}{h^d}\right)^2 + \frac{1}{n^3} \frac{\|K\|_\infty^3f_h(x)}{h^{3d}}.
\end{align*}

Consequently, since $\phi^{(3)}(z) = z^{-2}$, we have
\begin{equation}\label{eq.third_moment}
\begin{aligned}
\left| \bE \left[ \frac{\phi^{(3)}(f_h(x))}{6}(f_h(x)-\hat{f}_{h,2}(x))^3 \right]\right| &\le \frac{1}{6f_h(x)^2}\cdot \frac{2f_h(x)\|K\|_\infty^2}{n^2h^{2d}} \\
&\overset{\prettyref{eq:fh-smooth}}{\lesssim} \frac{1}{nh^d\ln n}. 
\end{aligned}
\end{equation}
Similarly, since $\phi^{(4)}(z) = -2z^{-3}, \xi_2\ge f_h(x)/2$ and $|f_h(x)-\xi_1|\le | f_h(x)-\hat{f}_{h,2}(x)|$, we have
\begin{equation}\label{eq.fourth_moment}
\begin{aligned}
& \left| \bE \left[ \frac{\phi^{(4)}(\xi_2)}{6}(f_h(x)-\hat{f}_{h,2}(x))^3(f_h(x)-\xi_1) \right]\right| \\
&\le \frac{1}{3(f_h(x)/2)^3}\cdot \bE(\hat{f}_{h,2}(x) - f_h(x))^4 \\
&\le \frac{1}{3(f_h(x)/2)^3}\cdot \left[\frac{3}{n^2}\left(\frac{\|K\|_\infty f_h(x)}{h^d}\right)^2 + \frac{1}{n^3}\cdot \frac{\|K\|_\infty^3f_h(x)}{h^{3d}}\right] \\
&\lesssim \frac{1}{n^2h^{2d}f_h(x)} + \frac{1}{n^3h^{3d}f_h(x)^2} \lesssim \frac{1}{nh^d\ln n}. 
\end{aligned}
\end{equation}
Then the desired bias bound is a direct consequence of \eqref{eq.concentration}--\eqref{eq.fourth_moment}. 

Next we upper bound the second moment of the perturbation. First we consider the case where one observation in $X^{(3)}$, say $X_n$, is replaced by an independent copy $X_n'$. In this case, by definition \prettyref{eq:H2x}, we have
\begin{align}
\hat{H}_2(x) - \hat{H}_2'(x) &= \mathbbm{1}\left(\hat{f}_{h,2}(x)\ge \frac{c_1\ln n}{4nh^d}\right)\left[\frac{\ln \hat{f}_{h,2}(x)}{n}(K_h(x-X_n') - K_h(x-X_n))\right. \nonumber\\
& \qquad + \left.  \frac{K_h(x-X_n') - K_h(x-X_n)}{n(n-1)\hat{f}_{h,2}(x)}\sum_{i=1}^{n-1}K_h(x-X_i) \right]. \label{eq.H2_diff1}
\end{align}
Note that $(X_n, X_n')$, $(X_1,\ldots,X_{n-1})$ and $\hat{f}_{h,2}(x)$ are mutually independent. Recall from \eqref{eq.H1_part1} that
\begin{align}\label{eq.H2_diff1_1}
\bE[(K_h(x-X_n') - K_h(x-X_n))^2] \le \frac{2\|K\|_\infty f_h(x)}{h^d}. 
\end{align}
Moreover, 
\begin{align}
\bE\left[\left(\frac{1}{n-1}\sum_{i=1}^{n-1}K_h(x-X_i) \right)^2 \right] 
& = ~ f_h(x)^2 + \frac{1}{n-1}\var\left(K_h(x-X_i) \right)  \nonumber \\
& \overset{\eqref{eq:varKh}}{\le} ~ f_h(x)^2 + \frac{\|K\|_\infty f_h(x)}{(n-1)h^d} \overset{\prettyref{eq:fh-smooth}}{\asymp} f_h(x)^2. \label{eq.H2_diff1_2}
\end{align}
To bound the remaining two expectations, let $E$ be the event $\hat{f}_{h,2}(x) \in (f_h(x)/2, 2f_h(x))$, which satisfies $\bP(E)\ge 1-n^{-5d}$ by \eqref{eq.concentration}. Then
\begin{align}
\ln^2 \hat{f}_{h,2}(x) \mathbbm{1}(E) &\lesssim 1 + (\ln f_h(x))^2, \label{eq.H2_diff1_3} \\
\hat{f}_{h,2}(x)^{-2} \mathbbm{1}(E) &\lesssim f_h(x)^{-2}. \label{eq.H2_diff1_4}
\end{align}

Therefore, by \eqref{eq.H2_diff1}--\eqref{eq.H2_diff1_4}, we conclude that
\begin{align}
\bE[(\hat{H}_2(x) - \hat{H}_2'(x))^2 \mathbbm{1}(E)] &\lesssim \frac{f_h(x)}{n^2h^d}\left[1+(\ln f_h(x))^2 + \frac{1}{f_h(x)^2}\left(f_h(x)^2 + \frac{f_h(x)}{nh^d} \right)\right] \nonumber\\
&\lesssim \frac{f_h(x)(1+(\ln f_h(x))^2)  }{n^2h^d} + \frac{1}{n^3h^{2d}} \nonumber\\
&\overset{\prettyref{eq:fh-smooth}}{\asymp} \frac{f_h(x)(1+(\ln f_h(x))^2)  }{n^2h^d}. \nonumber
\end{align}
Moreover, by \eqref{eq.concentration} and \eqref{eq:H2_sup}, 
\begin{align*}
\bE[(\hat{H}_2(x) - \hat{H}_2'(x))^2 \mathbbm{1}(E^c) ] \lesssim \left(1 + nh^d + \frac{1}{h^d}\ln\frac{n}{h^d} \right)^2 \bP(E^c) \lesssim n^{2-5d}(\ln n)^2,
\end{align*}
in view of the assumptions $h\le 1$ and $nh^d\ge 1$. Combining the previous two inequalities with \eqref{eq:fh-smooth} gives
\begin{align}\label{eq.H2_diff1_bound}
\bE[(\hat{H}_2(x) - \hat{H}_2'(x))^2] \lesssim  \frac{f_h(x)(1+(\ln f_h(x))^2)  }{n^2h^d}. 
\end{align}

Now we consider the case where one observation in $X^{(2)}$ is replaced by an independent copy, and $\hat{f}_{h,2}'(x)$ is the perturbed density estimate. Similar to the event $E$, let $E'$ be the event $\hat{f}_{h,2}'(x) \in (f_h(x)/2, 2f_h(x))$. Then on the event $E\cap E'$, by the assumption \prettyref{eq:fh-smooth},
we have $\hat{f}_{h,2}(x) \geq f_h(x)/2 \geq \frac{c_1\ln n}{4nh^d}$ and similarly for $\hat{f}_{h,2}'(x)$, and hence
\begin{align*}
\hat{H}_2(x) - \hat{H}_2'(x) = \frac{1}{2}D_1 + C_2D_2 + \frac{1}{2}C_3D_3,
\end{align*}
where
\begin{align*}
D_1 &\triangleq \hat{f}_{h,2}(x) - \hat{f}_{h,2}'(x), \\
D_2 &\triangleq \ln \hat{f}_{h,2}'(x - \ln \hat{f}_{h,2}(x)), \quad C_2 \triangleq \hat{f}_{h,3}(x) \\
D_3 &\triangleq \hat{f}_{h,2}'(x)^{-1} - \hat{f}_{h,2}(x)^{-1}, \quad C_3 \triangleq \frac{2}{n(n-1)}\sum_{i<j} K_h(x-X_{i}^{(3)})K_h(x-X_{j}^{(3)}).  
\end{align*}
By independence and the triangle inequality, it remains to upper bound the second moments of $D_1, D_2, D_3, C_2, C_3$ separately. By \eqref{eq.H2_diff1_1}, 
\begin{align}\label{eq.D1}
\bE[D_1^2] \le \frac{2\|K\|_\infty f_h(x)}{n^2h^d}. 
\end{align}
By the mean value theorem, 
\begin{align}
\bE[D_2^2\mathbbm{1}(E\cap E')] &\le \sup_{\xi\in (f_h(x)/2, 2f_h(x))}\left|\frac{d\ln \xi}{d\xi} \right|^2 \cdot \bE[D_1^2]  \lesssim \frac{1}{n^2h^df_h(x)} \label{eq.D2}, \\
\bE[D_3^2\mathbbm{1}(E\cap E')] &\le \sup_{\xi\in (f_h(x)/2, 2f_h(x))}\left|\frac{d(\xi^{-1})}{d\xi} \right|^2 \cdot \bE[D_1^2]  \lesssim \frac{1}{n^2h^df_h(x)^3} \label{eq.D3}.
\end{align}
Using $\bE[X^2] = (\bE X)^2 + \var(X)$, we also have
\begin{align}\label{eq.C2}
\bE[C_2^2] = f_h(x)^2 + \var(\hat{f}_{h,3}(x))\overset{\eqref{eq:varKh}}{\lesssim} f_h(x)^2 + \frac{f_h(x)}{nh^d}  \overset{\prettyref{eq:fh-smooth}}{\asymp} f_h(x)^2, 
\end{align}
\begin{equation} \label{eq.C3}
\begin{aligned}
\bE[C_3^2] &\le f_h(x)^4 + \frac{4(n-2)}{n(n-1)}f_h(x)^2\left(f_h(x)^2+ \frac{f_h(x)}{nh^d}\right) \\
&\qquad + \frac{2}{n(n-1)} \left(f_h(x)^2+ \frac{f_h(x)}{nh^d}\right)^2 \\
&\lesssim f_h(x)^4 + \frac{f_h(x)^3}{n^2h^d} + \frac{f_h(x)^2}{n^4h^{2d}} \overset{\prettyref{eq:fh-smooth}}{\asymp} f_h(x)^4 ,
\end{aligned}
\end{equation}
where \prettyref{eq.C3} is due to Lemma \ref{lemma.U-statistic}. Now combining \eqref{eq.D1}--\eqref{eq.C3}, the triangle inequality gives
\begin{align*}
&\bE[(\hat{H}_2(x) - \hat{H}_2'(x))^2\mathbbm{1}(E\cap E')]  \\
&\lesssim 
\Expect[D_1^2] + \bE[C_2^2] \Expect[D_2^2\mathbbm{1}(E\cap E')] + \bE[C_3^2] \Expect[D_3^2\mathbbm{1}(E\cap E')]  \\
&\lesssim \frac{f_h(x)}{n^2h^d}. 
\end{align*}

Moreover, by \eqref{eq.concentration} and \eqref{eq:H2_sup}, 
\begin{align*}
\bE[(\hat{H}_2(x) - \hat{H}_2'(x))^2 \mathbbm{1}(E^c \cup E'^c) ] &\lesssim \left(1 + nh^d + \frac{1}{h^d}\ln\frac{n}{h^d} \right)^2 \bP(E^c \cup E'^c) \\
&\lesssim n^{2-5d}(\ln n)^2,
\end{align*}
in view of the assumptions $h\le 1$ and $nh^d\ge 1$. Combining the previous two inequalities with \eqref{eq:fh-smooth} gives
\begin{align}\label{eq.H2_diff2_bound}
\bE[(\hat{H}_2(x) - \hat{H}_2'(x))^2] \lesssim  \frac{f_h(x) }{n^2h^d}. 
\end{align}
The proof is then completed by combining \eqref{eq.H2_diff1_bound} and \eqref{eq.H2_diff2_bound}.
\qed

\subsection{Proof of Lemma~\ref{lemma.performancereducediscrete}}

To analyze the performance of $\hat{H}$, we consider the following decomposition:
\begin{equation}\label{eq.discretization}
\begin{aligned}
\hat{H} - H(f) &= H(f_h) - H(f) + \hat{H} - H(f_h) \\
&= H(f_h) - H(f) + \left(\hat{H}_{\mathsf{discrete}} + \sum_{i=1}^{S}p_i\ln p_i\right) \\
&=h^d\sum_{i=1}^{S} (-\frac{p_i}{h^d}\ln \frac{p_i}{h^d}+\frac{1}{h^d}\int_{I_i}f(t)\ln f(t)dt) + \left(\hat{H}_{\mathsf{discrete}} + \sum_{i=1}^{S}p_i\ln p_i\right). 
\end{aligned}
\end{equation}
The first term deals with the approximation error in the first approximation stage. Note that by our Lipschitz ball assumption, we have
\begin{align*}
\|f-f_h\|_2^2 = \sum_{i=1}^{S} \int_{I_i} \left|f(t)-\frac{p_i}{h^d}\right|^2dt \lesssim L^2h^{2s}.
\end{align*}

Now we deal with the difference $H(f_h)-H(f)$. For $Z_i\sim\mathsf{Unif}(I_i)$, $Y=f(Z_i)$ and $\phi(x)=x\ln x$, we have
\begin{align*}
\frac{1}{h^d}\int_{I_i}f(t)\ln f(t)dt-\frac{p_i}{h^d}\ln \frac{p_i}{h^d} = \bE \phi(Y_i) - \phi(\bE Y_i)
\end{align*}
which is the gap in the Jensen's inequality. Using the equivalence between the $K$-functional and the modulus of smoothness~\cite{jiao2017bias}\cite[Chapter 6, Theorem 2.4]{Devore--Lorentz1993}, we have the following lemma. \begin{lemma}\label{lemma.K-functional}\cite{Strukov--Timan1977mathematical}
There exists a universal constant $C>0$ such that for any real-valued function $\phi$ defined on an interval $I\subset \mathbb{R}$ and any random variable $X$ supported on $I$ with a finite variance, 
\begin{align*}
|\bE\phi(X) - \phi(\bE X)| \le C\omega_I^2(\phi,\sqrt{\var(X)})
\end{align*}
where $\omega_I^2(\phi,t)$ is the second-order modulus of smoothness defined in \cite{Devore--Lorentz1993}.
\end{lemma}

To apply Lemma \ref{lemma.K-functional}, we note that~\cite[Chapter 2.9, Example 1]{Devore--Lorentz1993}
\begin{align*}
\omega_{\mathbb{R}_+}^2(-x\ln x, t) \asymp t
\end{align*}
and
\begin{align*}
\var(Y_i) = \frac{1}{h^d}\int_{I_i} \left|f(t)-\frac{p_i}{h^d}\right|^2dt.
\end{align*}
As a result,
\begin{align*}
|H(f_h)-H(f)|&\le \frac{1}{S}\sum_{i=1}^{S} |\bE \phi(Y_i) - \phi(\bE Y_i)|\lesssim \frac{1}{S}\sum_{i=1}^S \omega_{\mathbb{R}_+}^2(\phi,\sqrt{\var(Y_i)})\\
&\lesssim \frac{1}{S}\sum_{i=1}^{S}\sqrt{\var(Y_i)}\le \sqrt{\frac{1}{S}\sum_{i=1}^{S} \var(Y_i)} \\
&= \sqrt{\sum_{i=1}^{S} \int_{I_i} \left|f(t)-\frac{p_i}{h^d}\right|^2dt} \lesssim Lh^s.
\end{align*}

As for the second term in \eqref{eq.discretization}, we have $p_i\le \frac{L}{S}$ for any $i=1,\ldots,S$ since $\|f\|_\infty \le L$. Note that the discrete entropy can be related to the KL divergence as 
\begin{align*}
\sum_{i=1}^S -p_i\ln p_i = - \sum_{i=1}^S p_i\ln \frac{p_i}{S^{-1}} + \ln S, 
\end{align*}
where the likelihood ratio between $(p_1,\ldots,p_S)$ and the uniform distribution is upper bounded by $L$, the result in \cite{han2016minimax} yields the following risk bound:
\begin{align*}
\left(\sup_{f} \bE_f\left(\hat{H}_{\mathsf{discrete}}+\sum_{i=1}^{S}p_i\ln p_i\right)^2\right)^{\frac{1}{2}} \lesssim \frac{1}{nh^d\ln n} + \frac{\ln L}{\sqrt{n}}.
\end{align*}

Combining the previous inequalities together, we arrive at
\begin{align*}
\left(\sup_{f\in\text{Lip}_{s,p,d}(L), \|f\|_\infty\le L} \bE_f (\hat{H}-H(f))^2\right)^{\frac{1}{2}}\lesssim Lh^s + \frac{1}{nh^d\ln n}+\frac{\ln L}{\sqrt{n}}.
\end{align*}
Now choosing $h\asymp(Ln\ln n)^{-\frac{1}{s+d}}$ completes the proof. \qed

%

\subsection{Proof of Lemma \ref{lemma:varentropy}}
Let $X\sim f$, then
\begin{align*}
\int_{[0,1]^d} f(x)(\ln f(x))^2 dx &= \bE[(\ln f(X))^2] \\
&= \underbrace{\bE[(\ln f(X))^2\mathbbm{1}(f(X)\le e^{1/(p-1)}) ]}_{\triangleq A_1} \\
&\qquad + \underbrace{\bE[(\ln f(X))^2\mathbbm{1}(f(X)> e^{1/(p-1)}) ]}_{\triangleq A_2}. 
\end{align*}

Since $\max_{t\in [0,e^{1/(p-1)}]} t\ln^2 t = \max\{4e^{-2}, e^{1/(p-1)}/(p-1)^2\}$, we have
\begin{align*}
A_1 = \int_{x\in [0,1]^d: f(x)\le e^{p-1}} f(x)(\ln f(x))^2 dx \le \frac{4}{e^2} + \frac{e^{1/(p-1)}}{(p-1)^2}. 
\end{align*}
As for $A_2$, note that whenever $t>e$, 
\begin{align*}
\frac{d^2}{dt^2}\left[(\ln t)^2 \right] = \frac{2-2\ln t}{t^2} < 0. 
\end{align*}
Hence, by conditional Jensen's inequality, we have
\begin{align*}
A_2 &= \frac{1}{(p-1)^2} \bE[(\ln f(X)^{p-1})^2 \mathbbm{1}(f(X)^{p-1}>e) ] \\
&\le \frac{\bP(f(X)^{p-1}>e)}{(p-1)^2} \cdot \left(\ln \bE[f(X)^{p-1} | f(X)^{p-1} > e] \right)^2 \\
&\le \frac{\bP(f(X)^{p-1}>e)}{(p-1)^2} \cdot \left(\ln \frac{\int_{[0,1]^d} f(x)^pdx }{\bP(f(X)^{p-1}>e)} \right)^2 \\
&\le \frac{\bP(f(X)^{p-1}>e)}{(p-1)^2} \cdot \left(\ln \frac{L^p}{\bP(f(X)^{p-1}>e)} \right)^2 \\
&\le \frac{1 + p^2(\ln L)^2}{(p-1)^2}, 
\end{align*}
where the last inequality follows from
\begin{align*}
\max_{t\in [0,1]} t\left(\ln \frac{a}{t}\right)^2 = \begin{cases}
a/e &\text{if } 1\le a\le e \\
(\ln a)^2 & \text{if }a > e 
\end{cases} \le 1 + (\ln a)^2 
\end{align*}
whenever $a\ge 1$. Combining the upper bounds of $A_1$ and $A_2$ completes the proof. \qed

\subsection{Proof of Lemma \ref{lemma.U-statistic}}
It is straightforward to see $\bE U=(\bE X_1)^2$. As a result,
\begin{align*}
\var(U_2) = \frac{4}{n^2(n-1)^2}\sum_{1\le i<j\le n,1\le i'<j'\le n} \bE[X_iX_jX_{i'}X_{j'}] - (\bE X_1)^4.
\end{align*}
Denote by $I_0, I_1, I_2$ the set of indices $(i,j,i',j')$ where $(i,j)$ and $(i',j')$ have $0, 1$ and $2$ elements in common, respectively. Then
\begin{align*}
|I_0| &= \frac{n!}{2!2!(n-4)!} = \frac{n(n-1)(n-2)(n-3)}{4}, \\
|I_1| &= \frac{n!}{1!1!1!(n-3)!} = n(n-1)(n-2), \\
|I_2| &= \binom{n}{2} = \frac{n(n-1)}{2}.
\end{align*}

As a result,
\begin{align*}
\var(U_2) &= \frac{4}{n^2(n-1)^2}\left(|I_0|\cdot (\bE X_1)^4 + |I_1|\cdot (\bE X_1^2)(\bE X_1)^2 + |I_2|\cdot (\bE X_1^2)^2\right) - (\bE X_1)^4 \\
&= -\frac{4n-6}{n(n-1)}(\bE X_1)^4 + \frac{4(n-2)}{n(n-1)}(\bE X_1^2)(\bE X_1)^2 + \frac{2(\bE X_1^2)^2}{n(n-1)},
\end{align*}
and the result follows.
\qed

\subsection{Proof of Corollary \ref{cor.Ix}}
Recall that
\begin{align*}
\hat{H}(x) &= \min\sth{\hat{H}_1(x), \frac{1}{n^{1-2\varepsilon}h^d}}\mathbbm{1}\left(\hat{f}_{h,1}(x)<\frac{c_1\ln n}{nh^d}\right) \\
&\qquad + \hat{H}_2(x)\mathbbm{1}\left(\hat{f}_{h,1}(x)\ge \frac{c_1\ln n}{nh^d}\right).
\end{align*}

To establish \eqref{eq:Ix1}, we first show that for $f_h(x)\le \frac{2c_1\ln n}{nh^d}$, 
\begin{align}\label{eq.min}
\left| \bE \min\sth{\hat{H}_1(x), \frac{1}{n^{1-2\varepsilon}h^d}} + f_h(x)\ln f_h(x) \right| \lesssim \frac{1}{nh^d\ln n}. 
\end{align}
In fact, by Lemma \ref{lemma.non-smooth} and the Efron--Stein inequality (cf.~Lemma \ref{lemma:efron-stein}), 
\begin{align*}
\var(\hat{H}_1(x)) \le n\cdot \bE[(\hat{H}_1(x) - \hat{H}_1'(x))^2] \lesssim \frac{1}{n^{2-\varepsilon}h^{2d}}. 
\end{align*}
Then by the bias bound in Lemma \ref{lemma.non-smooth}, 
\begin{align*}
\bE[\hat{H}_1(x)^2] = (\bE[\hat{H}_1(x)])^2 + \var(\hat{H}_1(x)) \lesssim \frac{1}{n^{2-\varepsilon}h^{2d}}. 
\end{align*}
Using $\bE[X\mathbbm{1}(X\ge c)]\le c^{-1}\bE[X^2]$ for any random variable $X$ and $c>0$, we have
\begin{align*}
0\le \bE\hat{H}_1(x) - \bE \min\sth{\hat{H}_1(x), \frac{1}{n^{1-2\varepsilon}h^d}} \le n^{1-2\varepsilon}h^d\cdot \bE[\hat{H}_1(x)^2] \lesssim \frac{1}{n^{1+\varepsilon}h^d},
\end{align*}
which together with the bias bound in Lemma \ref{lemma.non-smooth} establishes \eqref{eq.min}. 

By \eqref{eq.min} and Lemma \ref{lemma.smooth}, we have
\begin{align*}
\bE[\hat{H}(x)|X^{(1)}] + f_h(x)\ln f_h(x)= I_1(x)\mathbbm{1}\left(\hat{f}_{h,1}(x)<\frac{c_1\ln n}{nh^d}\right) + I_2(x)\mathbbm{1}\left(\hat{f}_{h,1}(x)<\frac{c_1\ln n}{nh^d}\right),
\end{align*}
where
\begin{align*}
|I_1(x)| &\triangleq \left| \bE \min\sth{\hat{H}_1(x), \frac{1}{n^{1-2\varepsilon}h^d}} + f_h(x)\ln f_h(x)  \right|  \overset{\prettyref{eq:fhsup}}{\lesssim} \frac{1}{n^{1-2\varepsilon}h^d} + \frac{1}{h^d}\ln \frac{1}{h^d}, \\
|I_2(x)| &\triangleq  |\bE[\hat{H}_2(x)] + f_h(x)\ln f_h(x)| \overset{\prettyref{eq:fhsup},\prettyref{eq:H2_sup}}{\lesssim} 1 + nh^d + \frac{1}{h^d}\ln \frac{n}{h^d}
\end{align*}
for all $x\in\reals^d$. 
Moreover, by \eqref{eq.min} and Lemma \ref{lemma.smooth}, we also have $|I_1(x)|\lesssim (nh^d\ln n)^{-1}$ whenever $f_h(x)\le \frac{2c_1\ln n}{nh^d}$ and $|I_2(x)|\lesssim (nh^d\ln n)^{-1}$ whenever $f_h(x)\ge \frac{c_1\ln n}{2nh^d}$. Note that
\begin{align*}
\bE\left| \bE[\hat{H}(x)|X^{(1)}] + f_h(x)\ln f_h(x) \right|^2 &= I_1(x)^2\cdot \bP\left(\hat{f}_{h,1}(x)<\frac{c_1\ln n}{nh^d}\right) \nonumber \\
&\qquad + I_2(x)^2\cdot \bP\left(\hat{f}_{h,1}(x)\ge \frac{c_1\ln n}{nh^d}\right).
\end{align*}
Hence, by considering the three cases of $f_h(x)\in (0,\frac{c_1\ln n}{2nh^d}), [\frac{c_1\ln n}{2nh^d},\frac{2c_1\ln n}{nh^d}], (\frac{2c_1\ln n}{nh^d},\infty)$ separately and applying Lemma \ref{lemma.regime} in the first and third cases, we arrive at the desired bound \eqref{eq:Ix1}. 
 
The proof of the upper bound on the second moment of the perturbation is similar. If one observation in $X^{(2)}\cup X^{(3)}$ is perturbed, we have
\begin{align*}
\bE|\hat{H}(x) - \hat{H}'(x)|^2 & = J_1(x)^2\cdot  \bP\left(\hat{f}_{h,1}(x)<\frac{c_1\ln n}{nh^d}\right) + J_2(x)^2 \cdot \bP\left(\hat{f}_{h,1}(x)\ge \frac{c_1\ln n}{nh^d}\right),
\end{align*}
where
\begin{align*}
|J_1(x)| &\triangleq \bE \left| \min\sth{\hat{H}_1(x), \frac{1}{n^{1-2\varepsilon}h^d}} - \min\sth{\hat{H}_1'(x), \frac{1}{n^{1-2\varepsilon}h^d}}   \right|^2 \lesssim \frac{1}{(n^{1-2\varepsilon}h^d)^2}, \\
|J_2(x)| &\triangleq  \bE|\hat{H}_2(x) - \hat{H}_2'(x)|^2 \overset{\prettyref{eq:H2_sup}}{\lesssim} \left(1 + nh^d + \frac{1}{h^d}\ln \frac{n}{h^d}\right)^2
\end{align*}
for all $x$. Moreover, Lemma \ref{lemma.non-smooth} and \ref{lemma.smooth} give the improvement $|J_1(x)|\le \bE|\hat{H}_1(x) - \hat{H}_1'(x)|^2\lesssim \frac{1}{n^{3-\varepsilon}h^{2d}}$ whenever $f_h(x)\le \frac{2c_1\ln n}{nh^d}$ and $|J_2(x)|\lesssim \frac{f_h(x)(1+(\ln f_h(x))^2) }{n^2h^{d}}$ whenever $f_h(x)\ge \frac{c_1\ln n}{2nh^d}$. Again, by considering the three regimes $f_h(x)\in (0,\frac{c_1\ln n}{2nh^d}), [\frac{c_1\ln n}{2nh^d},\frac{2c_1\ln n}{nh^d}], (\frac{2c_1\ln n}{nh^d},\infty)$ and applying Lemma \ref{lemma.regime} in the first and third cases, we arrive at the desired bound \eqref{eq:Ix2}. 

\qed

\subsection{Proof of Lemma \ref{lemma:varentropy_lower}}
	We use the following two-point argument. 
	Fix $A>0,\varepsilon>0$ to be specified later.
	Let
	$f \in \text{\rm Lip}_{s,p,d}(L)$ be some fixed density supported on $[1/4,3/4]^d$ which is bounded from below by some constant, say $0.1$.
	Let $g$ be a dilation of $f$ defined by $g(x)=A^d f(A(x-{\bf 1}/2))$ such that $g$ is a density supported on $[1/2-1/(4A), 1/2+1/(4A)]^d$.
	We take $A \leq L^{1/(s+d)}$ so that $g$ is an element of $\text{\rm Lip}_{s,p,d}(L)$.
	
	Consider testing the hypothesis $H_0: X_i \iiddistr f_0=\frac{f+g}{2}$ against $H_1: X\iiddistr f_1 = \frac{1-\varepsilon}{2} f + \frac{1+\varepsilon}{2} g$.
Then
\begin{align*}
\chi^2(f_1\|f_0) &= \chi^2\pth{\frac{f+g}{2} \Big\| \frac{1-\varepsilon}{2} f + \frac{1+\varepsilon}{2} g} \\
&\leq  \chi^2(\Bern((1+\varepsilon)/2) \| \Bern(1/2)) = \varepsilon^2, 
\end{align*}
where the inequality follows from the data processing inequality for $\chi^2$-divergence.
Thus, choosing $\varepsilon=1/\sqrt{n}$ ensures the indistinguishability of the two hypotheses with $n$ observations.

Finally, for the separation of entropy values, by the Taylor expansion of the function $t\mapsto t \ln t$ and that $f(x)\ge 0.1$ everywhere on its support, we have 
\begin{align*}
|H(f_0)-H(f_1)| &= \frac{\varepsilon}{2} \left| \int_{[\frac{1}{4},\frac{3}{4}]^d} (f(x)-g(x)) \ln\frac{f(x)+g(x)}{2}dx  \right| \\
&\qquad + O\pth{\varepsilon^2 \int_{[\frac{1}{4},\frac{3}{4}]^d} (f(x)-g(x))^2 dx  }. 
\end{align*}
Using $f(x)\ge 0.1$ and $\|f\|_\infty + \|g\|_\infty = O(A^d)$, the above inequality further gives $|H(f_0) - H(f_1)| = O(\varepsilon \ln A + \varepsilon^2 A^{2d})$. Finally, in view of $\varepsilon = 1/\sqrt{n}$, taking $A = \min\{L^{1/(s+d)},n^{1/(4d)}\}$ yields the desired lower bound $n^{-1/2}\ln A \asymp n^{-1/2}\ln L$.
\qed

\subsection{Proof of Lemma \ref{lem.risk_multi_poisson}}
For any $\delta>0$, suppose $\hat{H}_n$ nearly attains the minimax risk in the usual sampling model with
\begin{align*}
\sup_{f\in \text{Lip}_{s,p,d}(L_0)} \bE_f(\hat{H}_n-H(f))^2 \le R_n^2 + \delta.
\end{align*} 
Now we use $\hat{H}_n$ to construct an estimator in the Poisson sampling model. Let $N\sim\mathsf{Poi}(n\|f\|_1)$ be the number of observations drawn in the Poisson sampling model, our estimator is just constructed as $\hat{H}=\hat{H}_N$. Note that conditioned on the event that $N=m$, the Poisson sampling model is just the usual sampling model with sample size $m$ and density $\frac{f}{\|f\|_1}$. As a result, for any $P\in \mathcal{P}$, we have
\begin{align*}
&\bE_{f_P}(\hat{H}-H(f_P))^2 \\
&\le 2\bE_{f_P}\pth{\hat{H}-H\left(\frac{f_P}{\|f_P\|_1}\right)}^2 + 2\pth{H(f_P)-H\left(\frac{f_P}{\|f_P\|_1}\right)}^2\\
&= 2\sum_{m=0}^\infty \bE_{f_P}\qth{\left(\hat{H}-H\left(\frac{f_P}{\|f_P\|_1}\right)\right)^2\bigg|N=m}\cdot \bP(N=m) + 2\pth{H(f_P)-H\left(\frac{f_P}{\|f_P\|_1}\right)}^2\\
&\le  2\sum_{n=0}^\infty (R_n^2+\delta) \bP(N=n) + 2\pth{H(f_P)-H\left(\frac{f_P}{\|f_P\|_1}\right)}^2\\
&\le 2c_0\cdot \bP(N\le n/2) + 2R_{\frac{n}{2}}^2\cdot \bP(N>n/2) + 2\delta + 2\pth{H(f_P)-H\left(\frac{f_P}{\|f_P\|_1}\right)}^2\\
&\le 2c_0\exp(-\frac{n}{5}) + 2R_{\frac{n}{2}}^2 + 2\delta + 4(\|f_P\|_1-1)^2c_0+4(\|f_P\|_1\ln \|f_P\|_1)^2
\end{align*}
where $c_0\triangleq\sup_{f\in \text{Lip}_{s,p,d}(L)}H(f)^2\lesssim (\ln L)^2$ in view of Lemma \ref{lemma:varentropy}, and in the last step we have used the Poisson tail bound. Note that by definition of $\mathcal{P}$,
\begin{align*}
| \|f_P\|_1 - 1 |\lesssim \frac{1}{nh^d(\ln n)^3\ln L}, 
\end{align*}
and thus we conclude that
\begin{align*}
R_n^P \lesssim R_{\frac{n}{2}}+\frac{1}{nh^d(\ln n)^3} + \sqrt{\delta}.
\end{align*}
The desired result follows by the arbitrariness of $\delta>0$. \qed

\subsection{Proof of Lemma \ref{lem.approx}}
We introduce some necessary definitions and results from approximation theory first. For functions defined on $[0,1]$, define the $r$-th order Ditzian--Totik modulus of smoothness by \cite{Ditzian--Totik1987}
\begin{align}\label{eq.DT_modulus}
\omega_\varphi^r(f,t)_\infty \triangleq \sup_{0<h\le t} \|\Delta_{h\varphi(x)}^r f(x)\|_\infty
\end{align}
where $\varphi(x)\triangleq \sqrt{x(1-x)}$. This quantity is related to the polynomial approximation error via the following lemma. 
\begin{lemma}[{\cite{Ditzian--Totik1987}}]\label{lem.DT_modulus}
	For any integer $u>0$ and $n>u$, there exists some constant $M_u$ depending only on $u$, such that for all $t\in (0,1)$ and $f$,
	\begin{align*}
	E_{n}(f;[0,1]) &\le M_u\omega_\varphi^u(f,1/n)_\infty\\
	\frac{M_u}{n^u}\sum_{l=0}^n (l+1)^{u-1}E_{l}(f;[0,1]) &\ge \omega_\varphi^u(f,1/n)_\infty.
	\end{align*}
\end{lemma}

Defining $r=q-1\ge 0$, we will make use of the second inequality in Lemma \ref{lem.DT_modulus} to prove that $E_{r,n}(x^{-r}\ln x;[cn^{-2},1])\gtrsim n^{2r}$ for some proper constant $c$. The case $r=0$ has been handled in \cite[Lemma 5]{wu2016minimax}, and we assume that $r>0$. Note that by definition \prettyref{eq:rationalapproxerror}, for any function $f$ and interval $I$, the rational approximation error can be related to the polynomial approximation error as
\begin{align*}
E_{r,n}(f;I) = \inf_{a_1,\ldots,a_r} E_{n}\pth{f(x)+\sum_{l=1}^r a_lx^{-l}; I}.
\end{align*}
Thus it suffices to prove that for any coefficients $a_1,\ldots,a_r\in \mathbb{R}$, $g(x)\triangleq x^{-r}\ln x+\sum_{l=1}^r a_lx^{-l}$ and $\tilde{g}(x)\triangleq g(cn^{-2}+(1-cn^{-2})x)$, the inequality
\begin{align*}
E_{n}(\tilde{g};[0,1]) \ge c'
\end{align*}
holds with proper universal constants $c,c'$ independent of the choices of $a_1,\ldots,a_r$. For the sake of simplicity, in the sequel we use the following equivalent definition of $g(x)$ (by a proper translation of $a_r$):
\begin{align*}
g(x) \triangleq x^{-r}\ln(c^{-1}n^2x) + \sum_{l=1}^r a_lx^{-l}.
\end{align*}

Let $u=1$ and $m=c^{-\frac{1}{2}}n$ with $c<1$ in Lemma \ref{lem.DT_modulus}, we have
\begin{align*}
E_{n}(\tilde{g};[0,1]) &\ge \frac{1}{m}\sum_{l=n}^{m}E_{l}(\tilde{g};[0,1]) \\
&\ge \frac{1}{M_1}\omega_{\varphi}^1(\tilde{g},\frac{1}{m})_\infty - \frac{1}{m}\sum_{l=0}^{n-1}E_{l}(\tilde{g};[0,1]) \\
&\ge \frac{1}{M_1}\omega_{\varphi}^1(\tilde{g},\frac{1}{m})_\infty - \frac{n}{m}E_{0}(\tilde{g};[0,1]).
\end{align*}
We distinguish into two cases. If $E_{0}(\tilde{g};[0,1])\le 2m^{2r}$, by definition of the Ditzian--Totik modulus of smoothness, there exists universal constants $0<A<B$ (independent of $c$) such that
\begin{align*}
\omega_\varphi^1(\tilde{g},\frac{1}{m}) &\ge \max_{z\in [A,B]}\left| \pth{\frac{z+1}{m^2}}^{-r}\ln (z+1) - \pth{\frac{z}{m^2}}^{-r}\ln z + \sum_{l=1}^r a_l\left(\pth{\frac{z+1}{m^2}}^{-l}-\pth{\frac{z}{m^2}}^{-l}\right) \right|\\
&= m^{2r} \cdot \max_{z\in [A,B]} \left| h(z) - \sum_{l=1}^r a_lm^{2(l-r)}h_l(z) \right| \\
&\ge m^{2r}\cdot \inf_{b_1,\ldots,b_r}\max_{z\in [A,B]} \left|h(z)-\sum_{l=1}^r b_lh_l(z)\right|
\end{align*}
where
\begin{align*}
h(z) &\triangleq (z+1)^{-r}\ln(z+1) - z^{-r}\ln z, \\
h_l(z) &\triangleq (z+1)^{-l} - z^{-l}, \qquad l=1,\ldots,r.
\end{align*}
It is clear that the functions $h_l,l=1,\ldots,r$ and $h$ are linearly independent over $[A,B]$, and hence
\begin{align*}
\inf_{b_1,\ldots,b_r}\max_{z\in [A,B]} \left|h(z)-\sum_{l=1}^r b_lh_l(z)\right| \ge C_1
\end{align*}
for some  universal constant $C_1>0$.
Hence, in this case we have
\begin{align}\label{eq.approx_lower_1}
E_{n}(\tilde{g};[0,1]) \ge \left(\frac{C_1}{M_1} - 2\sqrt{c}\right)m^{2r}.
\end{align}

Now we consider the second case where $E_{0}(\tilde{g};[0,1])>2m^{2r}$, which implies that 
\begin{align*}
E_{0}(\tilde{g};[0,1]) &\le \max_{x\in [m^{-2},1]}\left|x^{-r}\ln(m^2x) + \sum_{l=1}^r a_lx^{-l}\right| \\
&\le m^{2r}\cdot \max_{z\in [1,\infty)}z^{-r}\ln z + \sum_{l=1}^r |a_l|m^{2l} \\
&\le m^{2r} + \sum_{l=1}^r |a_l|m^{2l} \\
&\le \frac{1}{2}E_{0}(\tilde{g};[0,1]) + r\cdot \max_{1\le l\le r} |a_l|m^{2l}.
\end{align*}
As a result,
\begin{align*}
\max_{1\le l\le r} |a_l|m^{2l} \ge \frac{E_{0}(\tilde{g};[0,1])}{2r}.
\end{align*}
Suppose the maximum on the LHS is achieved by $l^*$. Then we have
\begin{align*}
\omega_\varphi^1(\tilde{g},\frac{1}{m}) &\ge \max_{z\in [A,B]}\left| \pth{\frac{z+1}{m^2}}^{-r}\ln (z+1) - \pth{\frac{z}{m^2}}^{-r}\ln z + \sum_{l=1}^r a_l\left(\pth{\frac{z+1}{m^2}}^{-l}-\pth{\frac{z}{m^2}}^{-l}\right) \right|\\
&= |a_{l^*}|m^{2l^*} \cdot \max_{z\in [A,B]} \left| h_{l^*}(z) - \frac{m^{2(r-l^*)}}{a_{l^*}}h(z) - \sum_{l\neq l^*} \frac{m^{2(l-l^*)}a_l}{a_{l^*}}h_l(z) \right| \\
&\ge |a_{l^*}|m^{2l^*}\cdot \inf_{b,b_l:l\neq l^*}\max_{z\in [A,B]} \left|h_{l^*}(z)-bh(z)-\sum_{l\neq l^*} b_lh_l(z)\right| \\
&\ge C_2|a_{l^*}|m^{2l^*}\\
&\ge \frac{C_2E_{0}(\tilde{g};[0,1])}{2r}
\end{align*}
where, again by the linear independence,
\begin{align*}
C_2 \triangleq \min_{1\le l^*\le r}\inf_{b,b_l:l\neq l^*}\max_{z\in [A,B]} \left|h_{l^*}(z)-bh(z)-\sum_{l\neq l^*} b_lh_l(z)\right|>0
\end{align*}
is a universal constant.
Hence
\begin{align}
E_{n}(\tilde{g};[0,1]) &\ge \left(\frac{C_2}{2rM_1} - \sqrt{c}\right) E_{0}(\tilde{g};[0,1])\nonumber \\\label{eq.approx_lower_2}
&\ge 2\left(\frac{C_2}{2rM_1} - \sqrt{c}\right)m^r.
\end{align}

Combining \eqref{eq.approx_lower_1} and \eqref{eq.approx_lower_2}, we conclude that
\begin{align*}
E_{r,n}(x^{-r}\ln x;[cn^{-2},1]) &\ge \min\left\{ \frac{C_1}{M_1}-2\sqrt{c}, 2\left(\frac{C_2}{2rM_1} - \sqrt{c}\right)\right\}\cdot m^r \\
&= \min\left\{ \frac{C_1}{M_1}-2\sqrt{c}, 2\left(\frac{C_2}{2rM_1} - \sqrt{c}\right)\right\}\cdot (\frac{n}{\sqrt{c}})^r.
\end{align*}
Now choosing $c>0$ small enough completes the proof of the lemma. \qed

\bibliographystyle{alpha}
\bibliography{di}

\newcommand{\etalchar}[1]{$^{#1}$}
\newcommand{\noopsort}[1]{}
\begin{thebibliography}{BDGVdM97}

\bibitem[BDGVdM97]{beirlant1997nonparametric}
Jan Beirlant, Edward~J Dudewicz, L{\'a}szl{\'o} Gy{\"o}rfi, and Edward~C
  Van~der Meulen.
\newblock Nonparametric entropy estimation: An overview.
\newblock {\em International Journal of Mathematical and Statistical Sciences},
  6(1):17--39, 1997.

\bibitem[BKB{\etalchar{+}}93]{bickel1993efficient}
Peter~J Bickel, Chris~AJ Klaassen, Peter~J Bickel, Y~Ritov, J~Klaassen, Jon~A
  Wellner, and Ya'acov Ritov.
\newblock {\em Efficient and adaptive estimation for semiparametric models}.
\newblock Johns Hopkins University Press Baltimore, 1993.

\bibitem[BLM13]{Boucheron--Lugosi--Massart2013}
St{\'e}phane Boucheron, G{\'a}bor Lugosi, and Pascal Massart.
\newblock {\em Concentration inequalities: A nonasymptotic theory of
  independence}.
\newblock Oxford University Press, 2013.

\bibitem[BM95]{birge1995estimation}
Lucien Birg{\'e} and Pascal Massart.
\newblock Estimation of integral functionals of a density.
\newblock {\em The Annals of Statistics}, pages 11--29, 1995.

\bibitem[BR88]{bickel1988estimating}
Peter~J Bickel and Ya'acov Ritov.
\newblock Estimating integrated squared density derivatives: sharp best order
  of convergence estimates.
\newblock {\em Sankhy{\=a}: The Indian Journal of Statistics, Series A}, pages
  381--393, 1988.

\bibitem[BSY19]{berrett2019efficient}
Thomas~B Berrett, Richard~J Samworth, and Ming Yuan.
\newblock Efficient multivariate entropy estimation via $ k $-nearest neighbour
  distances.
\newblock {\em The Annals of Statistics}, 47(1):288--318, 2019.

\bibitem[BZLV18]{bu2018estimation}
Yuheng Bu, Shaofeng Zou, Yingbin Liang, and Venugopal~V Veeravalli.
\newblock Estimation of kl divergence: Optimal minimax rate.
\newblock {\em IEEE Transactions on Information Theory}, 64(4):2648--2674,
  2018.

\bibitem[CH04]{costa2004geodesic}
Jose~A Costa and Alfred~O Hero.
\newblock Geodesic entropic graphs for dimension and entropy estimation in
  manifold learning.
\newblock {\em IEEE Transactions on Signal Processing}, 52(8):2210--2221, 2004.

\bibitem[CL03]{cai2003note}
T~Tony Cai and Mark~G Low.
\newblock A note on nonparametric estimation of linear functionals.
\newblock {\em Annals of statistics}, 31(4):1140--1153, 2003.

\bibitem[CL05]{cai2005nonquadratic}
T~Tony Cai and Mark~G Low.
\newblock Nonquadratic estimators of a quadratic functional.
\newblock {\em The Annals of Statistics}, 33(6):2930--2956, 2005.

\bibitem[CL11]{Cai--Low2011}
T~Tony Cai and Mark~G Low.
\newblock Testing composite hypotheses, {H}ermite polynomials and optimal
  estimation of a nonsmooth functional.
\newblock {\em The Annals of Statistics}, 39(2):1012--1041, 2011.

\bibitem[DF17]{delattre2017kozachenko}
Sylvain Delattre and Nicolas Fournier.
\newblock On the {Kozachenko--Leonenko} entropy estimator.
\newblock {\em Journal of Statistical Planning and Inference}, 185:69--93,
  2017.

\bibitem[DL93]{Devore--Lorentz1993}
Ronald~A DeVore and George~G Lorentz.
\newblock {\em Constructive approximation}, volume 303.
\newblock Springer, 1993.

\bibitem[DN90]{donoho1990minimax2}
David~L Donoho and Michael Nussbaum.
\newblock Minimax quadratic estimation of a quadratic functional.
\newblock {\em Journal of Complexity}, 6(3):290--323, 1990.

\bibitem[Don97]{donoho1997renormalizing}
David~L Donoho.
\newblock Renormalizing experiments for nonlinear functionals.
\newblock {\em Festschrift for Lucien Le Cam}, pages 167--181, 1997.

\bibitem[DT87]{Ditzian--Totik1987}
Zeev Ditzian and Vilmos Totik.
\newblock {\em Moduli of smoothness}.
\newblock Springer, 1987.

\bibitem[EHHG09]{el2009entropy}
Fidah El~Haje~Hussein and Yu~Golubev.
\newblock On entropy estimation by $m$-spacing method.
\newblock {\em Journal of Mathematical Sciences}, 163(3):290--309, 2009.

\bibitem[Fan91]{fan1991estimation}
Jianqing Fan.
\newblock On the estimation of quadratic functionals.
\newblock {\em The Annals of Statistics}, 19(3):1273--1294, 1991.

\bibitem[GM92]{goldstein1992optimal}
Larry Goldstein and Karen Messer.
\newblock Optimal plug-in estimators for nonparametric functional estimation.
\newblock {\em The Annals of Statistics}, 20(3):1306--1328, 1992.

\bibitem[GOV16]{gao2016breaking}
Weihao Gao, Sewoong Oh, and Pramod Viswanath.
\newblock Breaking the bandwidth barrier: Geometrical adaptive entropy
  estimation.
\newblock In {\em Advances in Neural Information Processing Systems}, pages
  2460--2468, 2016.

\bibitem[GOV18]{gao2018demystifying}
Weihao Gao, Sewoong Oh, and Pramod Viswanath.
\newblock Demystifying fixed $ k $-nearest neighbor information estimators.
\newblock {\em IEEE Transactions on Information Theory}, 64(8):5629--5661,
  2018.

\bibitem[GVdM91]{gyorfi1991nonparametric}
L{\'a}szl{\'o} Gy{\"o}rfi and Edward~C Van~der Meulen.
\newblock On the nonparametric estimation of the entropy functional.
\newblock In {\em Nonparametric functional estimation and related topics},
  pages 81--95. Springer, 1991.

\bibitem[Hal84]{hall1984limit}
Peter Hall.
\newblock Limit theorems for sums of general functions of $m$-spacings.
\newblock In {\em Mathematical Proceedings of the Cambridge Philosophical
  Society}, volume~96, pages 517--532. Cambridge University Press, 1984.

\bibitem[HI80]{hasminskii1980some}
Rafail~Z Hasminskii and Ildar~A Ibragimov.
\newblock Some estimation problems for stochastic differential equations.
\newblock In {\em Stochastic Differential Systems Filtering and Control}, pages
  1--12. Springer, 1980.

\bibitem[HJMW17]{Han--Jiao--Mukherjee--Weissman2017adaptive}
Yanjun Han, Jiantao Jiao, Rajarshi Mukherjee, and Tsachy Weissman.
\newblock On estimation of ${L}_r$-norms in {G}aussian white noise models.
\newblock {\em arXiv preprint arXiv:1710.03863}, 2017.

\bibitem[HJW16]{han2016minimax}
Yanjun Han, Jiantao Jiao, and Tsachy Weissman.
\newblock Minimax rate-optimal estimation of divergences between discrete
  distributions.
\newblock {\em arXiv preprint arXiv:1605.09124}, 2016.

\bibitem[HJWW17]{HJWW17v2}
Yanjun Han, Jiantao Jiao, Tsachy Weissman, and Yihong Wu.
\newblock Optimal rates of entropy estimation over {L}ipschitz balls.
\newblock {\em arXiv preprint}, Nov 2017.
\newblock \url{https://arxiv.org/abs/1711.02141v2}.

\bibitem[HKPT12]{hardle2012wavelets}
Wolfgang H{\"a}rdle, Gerard Kerkyacharian, Dominique Picard, and Alexander
  Tsybakov.
\newblock {\em Wavelets, approximation, and statistical applications}, volume
  129.
\newblock Springer Science \& Business Media, 2012.

\bibitem[HM87]{hall1987estimation}
Peter Hall and James~Stephen Marron.
\newblock Estimation of integrated squared density derivatives.
\newblock {\em Statistics \& Probability Letters}, 6(2):109--115, 1987.

\bibitem[HM93]{hall1993estimation}
Peter Hall and Sally~C Morton.
\newblock On the estimation of entropy.
\newblock {\em Annals of the Institute of Statistical Mathematics},
  45(1):69--88, 1993.

\bibitem[HSPVB07]{hlavavckova2007causality}
Katerina Hlav{\'a}{\v{c}}kov{\'a}-Schindler, Milan Palu{\v{s}}, Martin
  Vejmelka, and Joydeep Bhattacharya.
\newblock Causality detection based on information-theoretic approaches in time
  series analysis.
\newblock {\em Physics Reports}, 441(1):1--46, 2007.

\bibitem[INH87]{Ibragimov--Nemirovskii--Khasminskii1987some}
Ildar~A Ibragimov, Arkadi~S Nemirovskii, and Rafail~Z Hasminskii.
\newblock Some problems on nonparametric estimation in {G}aussian white noise.
\newblock {\em Theory of Probability \& Its Applications}, 31(3):391--406,
  1987.

\bibitem[Jaw77]{jawerth1977some}
Bj{\"o}rn Jawerth.
\newblock Some observations on {Besov} and {Lizorkin--Triebel} spaces.
\newblock {\em Mathematica Scandinavica}, 40(1):94--104, 1977.

\bibitem[JHW17]{jiao2017bias}
Jiantao Jiao, Yanjun Han, and Tsachy Weissman.
\newblock Bias correction with jackknife, bootstrap, and taylor series.
\newblock {\em arXiv preprint arXiv:1709.06183}, 2017.

\bibitem[JHW18]{jiao2018minimax}
Jiantao Jiao, Yanjun Han, and Tsachy Weissman.
\newblock Minimax estimation of the $l_1$ distance.
\newblock {\em IEEE Transactions on Information Theory}, 64(10):6672--6706,
  2018.

\bibitem[Joe89]{joe1989estimation}
Harry Joe.
\newblock Estimation of entropy and other functionals of a multivariate
  density.
\newblock {\em Annals of the Institute of Statistical Mathematics},
  41(4):683--697, 1989.

\bibitem[JS77]{johnen1977equivalence}
Hans Johnen and Karl Scherer.
\newblock On the equivalence of the {$K$}-functional and moduli of continuity
  and some applications.
\newblock In {\em Constructive theory of functions of several variables}, pages
  119--140. Springer, 1977.

\bibitem[JVHW15]{Jiao--Venkat--Han--Weissman2015minimax}
Jiantao Jiao, Kartik Venkat, Yanjun Han, and Tsachy Weissman.
\newblock Minimax estimation of functionals of discrete distributions.
\newblock {\em IEEE Transactions on Information Theory}, 61(5):2835--2885,
  2015.

\bibitem[KKPW14]{krishnamurthy2014nonparametric}
Akshay Krishnamurthy, Kirthevasan Kandasamy, Barnabas Poczos, and Larry
  Wasserman.
\newblock Nonparametric estimation of {R\'enyi} divergence and friends.
\newblock In {\em International Conference on Machine Learning}, pages
  919--927, 2014.

\bibitem[KKPW15]{kandasamy2015nonparametric}
Kirthevasan Kandasamy, Akshay Krishnamurthy, Barnabas Poczos, and Larry
  Wasserman.
\newblock Nonparametric von mises estimators for entropies, divergences and
  mutual informations.
\newblock In {\em Advances in Neural Information Processing Systems}, pages
  397--405, 2015.

\bibitem[KP96]{kerkyacharian1996estimating}
G{\'e}rard Kerkyacharian and Dominique Picard.
\newblock Estimating nonquadratic functionals of a density using haar wavelets.
\newblock {\em The Annals of Statistics}, 24(2):485--507, 1996.

\bibitem[KSG04]{Kraskov2004}
Alexander Kraskov, Harald St{\"o}gbauer, and Peter Grassberger.
\newblock Estimating mutual information.
\newblock {\em Physical Review E}, 69(6):066138, 2004.

\bibitem[KT12]{korostelev2012minimax}
Aleksandr~Petrovich Korostelev and Alexandre~B Tsybakov.
\newblock {\em Minimax theory of image reconstruction}, volume~82.
\newblock Springer Science \& Business Media, 2012.

\bibitem[Lau96]{laurent1996efficient}
B{\'e}atrice Laurent.
\newblock Efficient estimation of integral functionals of a density.
\newblock {\em The Annals of Statistics}, 24(2):659--681, 1996.

\bibitem[Lev78]{levit1978asymptotically}
Boris~Ya Levit.
\newblock Asymptotically efficient estimation of nonlinear functionals.
\newblock {\em Problemy Peredachi Informatsii}, 14(3):65--72, 1978.

\bibitem[LNS99]{lepski1999estimation}
Oleg Lepski, Arkady Nemirovski, and Vladimir Spokoiny.
\newblock On estimation of the {$L_r$} norm of a regression function.
\newblock {\em Probability theory and related fields}, 113(2):221--253, 1999.

\bibitem[MNR17]{mukherjee2017semiparametric}
Rajarshi Mukherjee, Whitney~K Newey, and James~M Robins.
\newblock Semiparametric efficient empirical higher order influence function
  estimators.
\newblock {\em arXiv preprint arXiv:1705.07577}, 2017.

\bibitem[MSGHI16]{moon2016nonparametric}
Kevin~R Moon, Kumar Sricharan, Kristjan Greenewald, and Alfred~O Hero~III.
\newblock Nonparametric ensemble estimation of distributional functionals.
\newblock {\em arXiv preprint arXiv:1601.06884}, 2016.

\bibitem[Nem00]{nemirovski2000topics}
Arkadi Nemirovski.
\newblock Topics in non-parametric statistics.
\newblock {\em Ecole d’Et{\'e} de Probabilit{\'e}s de Saint-Flour}, 28, 2000.

\bibitem[Nus96]{nussbaum1996asymptotic}
Michael Nussbaum.
\newblock Asymptotic equivalence of density estimation and gaussian white
  noise.
\newblock {\em The Annals of Statistics}, pages 2399--2430, 1996.

\bibitem[Pan04]{Paninski2004}
Liam Paninski.
\newblock Estimating entropy on $m$ bins given fewer than $m$ samples.
\newblock {\em Information Theory, IEEE Transactions on}, 50(9):2200--2203,
  2004.

\bibitem[PR16]{patschkowski2016adaptation}
Tim Patschkowski and Angelika Rohde.
\newblock Adaptation to lowest density regions with application to support
  recovery.
\newblock {\em The Annals of Statistics}, 44(1):255--287, 2016.

\bibitem[PY08]{paninski2008undersmoothed}
Liam Paninski and Masanao Yajima.
\newblock Undersmoothed kernel entropy estimators.
\newblock {\em IEEE Transactions on Information Theory}, 54(9):4384--4388,
  2008.

\bibitem[R{\'e}n59]{Renyi59}
Alfr{\'e}d R{\'e}nyi.
\newblock On the dimension and entropy of probability distributions.
\newblock {\em Acta Mathematica Hungarica}, 10(1--2), Mar 1959.

\bibitem[RLM{\etalchar{+}}17]{robins2017higher}
James Robins, Lingling Li, Rajarshi Mukherjee, Eric~Tchetgen Tchetgen, and Aad
  van~der Vaart.
\newblock Higher order estimating equations for high-dimensional models.
\newblock {\em Annals of statistics}, 45(5):1951--1987, 2017.

\bibitem[RLTvdV08]{robins2008higher}
James Robins, Lingling Li, Eric Tchetgen, and Aad van~der Vaart.
\newblock Higher order influence functions and minimax estimation of nonlinear
  functionals.
\newblock In {\em Probability and Statistics: Essays in Honor of David A.
  Freedman}, pages 335--421. Institute of Mathematical Statistics, 2008.

\bibitem[RSH18]{ray2018asymptotic}
Kolyan Ray and Johannes Schmidt-Hieber.
\newblock Asymptotic nonequivalence of density estimation and gaussian white
  noise for small densities.
\newblock {\em arXiv preprint arXiv:1802.03425}, 2018.

\bibitem[SP16]{singh2016finite}
Shashank Singh and Barnab{\'a}s P{\'o}czos.
\newblock Finite-sample analysis of fixed-$k$ nearest neighbor density
  functional estimators.
\newblock In {\em Advances in Neural Information Processing Systems}, pages
  1217--1225, 2016.

\bibitem[SRH12]{sricharan2012estimation}
Kumar Sricharan, Raviv Raich, and Alfred~O Hero.
\newblock Estimation of nonlinear functionals of densities with confidence.
\newblock {\em IEEE Transactions on Information Theory}, 58(7):4135--4159,
  2012.

\bibitem[ST77]{Strukov--Timan1977mathematical}
LI~Strukov and AF~Timan.
\newblock Mathematical expectation of continuous functions of random variables.
  smoothness and variance.
\newblock {\em Siberian Mathematical Journal}, 18(3):469--474, 1977.

\bibitem[Ste79]{stein2016singular}
Elias~M Stein.
\newblock {\em Singular integrals and differentiability properties of functions
  (PMS-30)}.
\newblock Princeton university press, 1979.

\bibitem[Ste86]{steele86}
J~Michael Steele.
\newblock An {E}fron-{S}tein inequality for nonsymmetric statistics.
\newblock {\em The Annals of Statistics}, pages 753--758, 1986.

\bibitem[Tao15]{tao-torus}
Terrence Tao.
\newblock {The Hardy-Littlewood maximal inequality and applications. Lecture
  notes for MATH 247A: Fourier analysis}.
\newblock 2015.
\newblock \url{http://www.math.ucla.edu/~tao/247a.1.06f/notes3.pdf}.

\bibitem[Tim63]{timan63}
Aleksandr~Filippovich Timan.
\newblock {\em Theory of approximation of functions of a real variable}.
\newblock Pergamon Press, 1963.

\bibitem[TLRvdV08]{tchetgen2008minimax}
Eric Tchetgen, Lingling Li, James Robins, and Aad van~der Vaart.
\newblock Minimax estimation of the integral of a power of a density.
\newblock {\em Statistics \& Probability Letters}, 78(18):3307--3311, 2008.

\bibitem[Tot94]{Totik1994approximation}
Vilmos Totik.
\newblock Approximation by {B}ernstein polynomials.
\newblock {\em American Journal of Mathematics}, pages 995--1018, 1994.

\bibitem[Tsy09]{Tsybakov2008}
A.~Tsybakov.
\newblock {\em Introduction to Nonparametric Estimation}.
\newblock Springer-Verlag, 2009.

\bibitem[TVdM96]{tsybakov1996root}
Alexandre~B Tsybakov and EC~Van~der Meulen.
\newblock Root-$n$ consistent estimators of entropy for densities with
  unbounded support.
\newblock {\em Scandinavian Journal of Statistics}, pages 75--83, 1996.

\bibitem[VdV00]{van2000asymptotic}
Aad~W Van~der Vaart.
\newblock {\em Asymptotic statistics}, volume~3.
\newblock Cambridge university press, 2000.

\bibitem[VE92]{van1992estimating}
Bert Van~Es.
\newblock Estimating functionals related to a density by a class of statistics
  based on spacings.
\newblock {\em Scandinavian Journal of Statistics}, pages 61--72, 1992.

\bibitem[VV11]{valiant2011power}
Gregory Valiant and Paul Valiant.
\newblock The power of linear estimators.
\newblock In {\em Foundations of Computer Science (FOCS), 2011 IEEE 52nd Annual
  Symposium on}, pages 403--412. IEEE, 2011.

\bibitem[WKV09]{WKV09b}
Qing Wang, Sanjeev~R Kulkarni, and Sergio Verd{\'u}.
\newblock Universal estimation of information measures for analog sources.
\newblock {\em Foundations and Trends in Communications and Information
  Theory}, 5(3):265--353, 2009.

\bibitem[WY16]{wu2016minimax}
Yihong Wu and Pengkun Yang.
\newblock Minimax rates of entropy estimation on large alphabets via best
  polynomial approximation.
\newblock {\em IEEE Transactions on Information Theory}, 62(6):3702--3720,
  2016.

\bibitem[WY19]{wu2019chebyshev}
Yihong Wu and Pengkun Yang.
\newblock Chebyshev polynomials, moment matching, and optimal estimation of the
  unseen.
\newblock {\em The Annals of Statistics}, 47(2):857--883, 2019.

\end{thebibliography}
\end{document}